%% file: main.tex
\documentclass[11pt,a4paper]{article}
\usepackage[margin=1in]{geometry} 
\usepackage[english]{babel} 
\usepackage{graphicx} 
\usepackage{float} 
\usepackage{subcaption} 
\usepackage{tcolorbox} 
\usepackage{amsmath, amsfonts, bm} 
\usepackage{tikz, tikz-cd, tikz-3dplot} 
\usetikzlibrary{decorations.pathmorphing, decorations.pathreplacing, patterns, angles, quotes}

\usepackage{authblk}
\usepackage{derivative}
\usepackage{hyperref}
\usepackage{xcolor}

\input{macros}

\title{A domain decomposition strategy for natural imposition of mixed boundary conditions in port-Hamiltonian systems}

\author[1]{Sjoerd D.M. de Jong}
\author[2]{Andrea Brugnoli}
\author[3]{Ramy Rashad}
\author[4]{Yi Zhang}
\author[5]{Stefano Stramigioli}

\affil[1]{Department of Microelectronics, Delft University of Technology, Delft, Netherlands}
\affil[2]{ICA, Universit\'e de Toulouse, ISAE–SUPAERO, INSA, CNRS, MINES ALBI, UPS, Toulouse, France}
\affil[3]{Control and Instrumentation Engineering Department, King Fahd University of Petroleum and Minerals, Saudi Arabia}
\affil[4]{School of Mathematics and Computing Science, Guilin University of Electronic Technology, Guilin, China}
\affil[5]{Robotics and Mechatronics Department, University of Twente, The Netherlands}

\date{}

\begin{document}

\maketitle

\input{Text/abstract}

\input{Text/Introduction}

\input{Text/Discretization}

\input{Text/numerical_examples}

\input{Text/pH_intrinsic_nonlinear_beam}

\input{Text/pH_WaveEquation}

\input{Text/pH_elasticity}

\input{Text/pH_Mindlin_plate}

\input{Text/Conclusion}

\section*{Funding}
This research was carried out under project number N21006 in the framework of the Partnership Program of the Material innovation institute M2i (www.m2i.nl) and the Dutch Research Council (www.nwo.nl). This work was supported by the PortWings project funded by the European Research Council [Grant Agreement No. 787675]. The research of Yi Zhang is supported by the Natural Science Foundation of Guangxi under grant number 2024JJB110005.

\section*{Acknowledgements}

The second author would like to thank Flavio Cardoso Ribeiro (Technical Institute of Aeronautics) and Paul Kotyczka (Technical University of Munich) for fruitful discussions on this topic.

\addcontentsline{toc}{section}{References}
\bibliographystyle{elsarticle-num}
\bibliography{References}

\end{document}

%% file: macros.tex
\renewcommand\d{\ensuremath{\mathrm{d}}}

\newcommand{\bbR}{\mathbb{R}}

\newcommand{\bbS}{\mathbb{S}}

\DeclareMathOperator*{\grad}{grad}
\DeclareMathOperator*{\Grad}{Grad}
\renewcommand{\div}{\operatorname{div}}
\DeclareMathOperator*{\Div}{Div}

\DeclareMathOperator*{\curl}{curl}
\DeclareMathOperator*{\rot}{rot}
\DeclareMathOperator*{\Rot}{Rot}

\DeclareMathOperator*{\sym}{sym}

\newcommand{\innerproduct}[3][\Omega]{(#2, #3)_{#1}}
\newcommand{\boundary}[3][\partial\Omega]{\langle #2, #3\rangle_{#1}}

\newtheorem{remark}{Remark}
\newtheorem{assumption}{Assumption}
\newtheorem{example}{Example}
\newtheorem{definition}{Definition}
\newtheorem{proof}{Proof}

\newtheorem{proposition}{Proposition}

\colorlet{rev1}{blue!80!black} 
\colorlet{rev2}{red!80!black} 
\colorlet{rev3}{magenta!80!black} 

\newcommand{\revone}[1]{{\textcolor{black}{#1}}}
\newcommand{\revtwo}[1]{{\textcolor{black}{#1}}}
\newcommand{\revthree}[1]{{\textcolor{black}{#1}}}

%% file: Text/abstract.tex
\begin{abstract}
In this contribution, a finite element scheme to impose mixed boundary conditions without introducing Lagrange multipliers is presented for \revtwo{hyperbolic systems} described as port-Hamiltonian systems. The strategy relies on finite element exterior calculus and domain decomposition to interconnect two systems with \revone{dual input-output behavior}. The spatial domain is split into two parts by introducing an arbitrary interface. Each subdomain is discretized with a mixed finite element formulation that introduces a uniform boundary condition in a natural way as the input. In each subdomain the finite element spaces are selected from a finite element subcomplex to obtain a stable discretization. The two systems are then interconnected together by making use of a feedback interconnection. This is achieved by discretizing the boundary inputs using appropriate spaces that couple the two formulations. The final systems include all boundary conditions explicitly and do not contain any Lagrange multiplier. \revone{Time integration is performed using the implicit midpoint or St\"ormer-Verlet scheme}. \revtwo{The method can also be applied to semilinear systems containing algebraic nonlinearities.} The proposed strategy is tested on different examples: \revtwo{geometrically exact intrinsic beam model}, the wave equation, \revthree{membrane elastodynamics and the Mindlin plate}. Numerical tests assess the conservation properties of the scheme, the effectiveness of the methodology and its robustness against shear locking phenomena.   
\end{abstract}

%% file: Text/Introduction.tex
\section{Introduction}

To simulate, design and analyze modern engineering technologies, modular modeling tools are of great importance, as they allow to simplify validation and verification, speed up prototyping and encapsulate complexity. Paradigms based on a modular description of systems are implemented in many widespread libraries like \textsc{Simulink}\footnote{\url{https://www.mathworks.com/products/simulink.html}} or  \textsc{Dimola}\footnote{\url{https://www.3ds.com/products/catia/dymola}}.  In many cases, a reliable description of a complex technological devices is achieved by using coupled systems of partial differential equations (PDE) where different physics operate together. In recent years, the port-Hamiltonian (pH) formalism  \cite{rashad2020review} has established itself as a sound and powerful mathematical framework for modeling and control of complex multiphysical systems. At the core of this framework lies the idea of composability, i.e. the fact that interconnecting port-Hamiltonian systems (pHs) leads to another system of the same kind.  \\

The theory of port-Hamiltonian systems is built upon a rich geometrical structure based on exterior calculus and issues may arise if this structure is not preserved at the numerical level \cite{bochev2003discourse,arnold2006finite,christiansen2011topics}. Structure preserving techniques attempt to capture as much of the underlying structures as possible. To this aim, many strategies have been proposed throughout the years, such as mimetic finite differences \cite{bochev2006principles,lipnikov2014}, discrete exterior calculus \cite{hirani2003discrete}, finite element exterior calculus \cite{arnold2006finite} and many others. When devising discretization schemes for port-Hamiltonian systems, boundary conditions have a prominent role in the discussion. This is due to the connection of port-Hamiltonian systems to the concept of Stokes-Dirac structure \cite{van2002hamiltonian}. This geometrical structure characterizes all admissible boundary flows into a spatial domain and is agnostic to the actual boundary conditions of the problem. The way boundary conditions are included in the model is related to the numerical method used. In a finite element context, boundary conditions are either imposed strongly by incorporating them in the discrete spaces used to approximate the variables or weakly when they explicitly appear in the weak formulation \cite{quarteroni2008numerical}. Weak imposition of the boundary conditions typically arises from the variational formulation in a natural manner via integration by parts. There is no general consensus on whether it is preferable to use a weak or strong formulation and the best choice is strongly problem and method dependent \cite{bazilevs2007,eriksson2022}. Strong imposition of the boundary conditions in dynamical systems leads to differential-algebraic equations that are more difficult to solve than ordinary differential equations \cite{kunkel2006differential}. In the port-Hamiltonian community a general effort has been made to incorporate mixed boundary conditions in an explicit manner, see for instance \cite{seslija2014explicit} for a discrete exterior calculus formulation, \cite{kotyczka2018weak} for a Galerkin scheme based on Whitney forms, \cite{cardoso2020pfem} for a mixed finite element framework and  \cite{kumar2024dg} for discontinuous Galerkin discretization based on finite element exterior calculus (FEEC). Wave propagation phenomena exhibit a primal dual structure that was first highlighted in \cite{joly2003}. Therein however no connection with differential geometry is established.  In \cite{Brugnoli2022DualCalculus} the authors used a finite element exterior calculus to highlight the fundamental primal-dual structure of pHs. The two formations are related by the Hodge operator and the resulting scheme is called dual-field as each variable is represented in dual finite element bases. The use of a dual field finite element formulation was initially introduced in \cite{Zhang2022ADomains} as a way of handling the convective non-linearity of Navier-Stokes equations in an explicit manner and still obtaining a conservative scheme in terms of mass, helicity and energy. The work of \cite{Zhang2022ADomains} focused on periodic domains only without dealing with boundary conditions. In port-Hamiltonian systems, the dual field representation allows obtaining the topological power balance that characterizes the Dirac structure when inhomogeneous mixed boundary conditions are considered. Furthermore, it clearly shows that the two formulations \revone{treat the boundary conditions in a dual manner, i.e. the natural boundary conditions for one formulation are essential for the other and vice-versa}, which leads to the question of how can this primal-dual structure be exploited for incorporating them.\\
 
 In the present contribution, the dual-field representation is employed to achieve weak imposition of mixed boundary conditions in hyperbolic systems. In particular, this work formalizes previous results discussed in \cite{brugnoli2020mixed,brugnoli2022explicit} using finite element exterior calculus. The spatial domain is decomposed using an internal interface that separates the two boundary subpartitions when a single boundary condition applies. On each subdomain a mixed finite element formulation is used in such a way that the boundary condition is included naturally. Each mixed formulation uses a pair a finite elements that constitute a Hilbert subcomplex and thus is stable and structure preserving. The two formulations are then interconnected together on the shared interface by means of a feedback (or in port-Hamiltonian jargon a gyrator interconnection) that enforces in weak manner the continuity of the finite element spaces. The resulting system incorporates the mixed boundary conditions of the problem in a completely weak manner and does not require Lagrange multipliers. \revone{The proposed methodology is reminiscent of Dirichlet-Neumann alternating Schwarz methods for non-overlapping domain decomposition-based coupling, cf. \cite{rodriguez2025transmission} for an application to linear elasticity models. However, in the classical domain-decomposition method only one primal formulation and the coupling conditions is achieved via an iterative approach. The present contribution is not really concerned with accelerating numerical methods but rather with showing that the employment of a primal-dual formulation makes it possible to avoid the usage of Lagrange multipliers for subdomain coupling.}  
 \revtwo{The strategy can also be applied to semilinear model containing algebraic nonlinearities.} Even if the methodology is discussed for hyperbolic port-Hamiltonian systems, it can be extended to static elliptic problems. \revone{For the time integration the implicit midpoint scheme and the St\"ormer Verlet method are considered. This choice guarantees the preservation of the power balance in each subdomain \cite{kotyczka2019discrete}. The implicit midpoint preserves the overall energy but requires the solution of a monolithic system. St\"ormer-Verlet decouples the two subdomains but does not enforce energy preservation exactly.} 
 The proposed approach will be shown to be accurate, have proper convergence and to be able to preserve certain mathematical, and thereby physical, properties at the discrete level. To demonstrate this, different physical examples are considered: \revtwo{the nonlinear geometrically exact intrinsic beam model}, the wave equation in two dimensions, \revthree{membrane elastodynamics and the Mindlin plate}. The examples chosen showcase the versatility of our approach in different physical-domains as well as different dimensions.  \revthree{Furthermore, they show that the proposed discretization does not suffer from shear locking phenomena.} To summarize main novel results of the paper are the following: 
 \begin{itemize} 
 \item Formalizes a dual-field representation using Finite Element Exterior Calculus (FEEC) to enforce mixed boundary conditions without the need for Lagrange multipliers. The domain is decomposed into subdomains using a Hilbert subcomplex pair for each, ensuring the formulation remains stable and preserves the underlying physical structure. 
 \item Employs a gyrator (port-Hamiltonian feedback) interconnection at the internal interface to enforce continuity between finite element spaces in a weak manner. \item Distinguished from classical Dirichlet-Neumann alternating Schwarz methods by using a primal-dual formulation that avoids iterative coupling for subdomain interconnection. 
 \item Demonstrates that the Implicit Midpoint rule preserves overall energy, while the Störmer-Verlet method allows for subdomain decoupling at the cost of exact energy preservation. 
 \item The discretization approach is shown to be accurate, convergent, and specifically avoids shear locking phenomena in structural mechanics. 
 \item The method is primarily designed for hyperbolic port-Hamiltonian systems, but explicitly extensible to static elliptic problems and semilinear models with algebraic nonlinearities. It can be extended to electromagnetic phenomena and multiphysics coupling.
 \end{itemize}

 The outline of the rest of the paper is as follows. The assumptions of the study and the mixed discretization approach based on finite element subcomplexes are presented in Sec. \ref{sec:Discretization}. The domain-decomposition strategy is presented in Sec.\ref{sec:domain_decomposition} including the choice of the finite element basis for the boundary input made to couple the two formulation on the interface. The time integration schemes are discussed in Sec. \ref{sec:time-integration}. Sec. \ref{sec:examples} presents the numerical examples.

%% file: Text/Discretization.tex
\section{Galerkin discretization of port-Hamiltonian systems}\label{sec:Discretization}

The general class of port-Hamiltonian systems is now presented. \revtwo{A brief introduction on port-Hamiltonian systems is given by means of the wave equation as an example.} Then we recall the mixed finite element Galerkin discretization presented in \cite{joly2003}. This discretization is such to retain the Hamiltonian structure at the discrete level.

\revtwo{
\subsection{An introductory example: the wave equation}
The propagation of acoustic waves in $\Omega\subset \mathbb{R}^d$ is described by the following hyperbolic partial differential equation, that determines the time-dependent field $\phi(t): \Omega \times [0, T_{\mathrm{end}}] \rightarrow \mathbb{R}$
\begin{equation}\label{eq:irr_wave}
    \partial_{tt}^2 \phi - \div\grad \phi = 0,
\end{equation}
together with time-varying Dirichlet boundary condition
\begin{equation}\label{eq:dirbc_wave}
    \phi|_{\partial\Omega} = g(t).
\end{equation}
The total energy is given by the sum of kinetic and potential energy
$$
H = \frac{1}{2} \int_\Omega (\partial_t \phi)^2 + ||\grad \phi||^2 \; \rm{d}\Omega.
$$
To highlight the Hamiltonian structure of the wave equation, consider the variables
\begin{equation}
    \alpha := \partial_t \phi, \qquad \bm{\beta} = \grad \phi.
\end{equation}
Equation \eqref{eq:irr_wave}, together with the boundary condition \eqref{eq:dirbc_wave}, can now be recast into a first order system 
\begin{equation}\label{eq:mix_wave}
    \begin{pmatrix}
        \partial_t \alpha \\
        \partial_t \bm{\beta}
    \end{pmatrix} = \underbrace{\begin{bmatrix}
        0 & \div \\
        \grad & 0 \\
    \end{bmatrix}}_{J}
    \begin{pmatrix}
        \alpha \\
        \bm{\beta}
    \end{pmatrix}, \qquad \alpha|_{\partial \Omega} = \partial_t g:= u.
\end{equation}
Notice that the operator $J$ is formally skew-adjoint, as for compactly supported function the adjoint of the gradient is minus the divergence $\grad^*= -\div$. Here $u$ corresponds to a control input applied to the boundary. Indeed port-Hamiltonian systems are boundary controlled systems and the boundary conditions coincide with inputs that describe interactions with the external environment. Notice that the Hamiltonian is quadratic in the new variables, making the variational derivative easier to compute
$$
H = \frac{1}{2} \int_\Omega \alpha^2 + ||\bm \beta||^2 \mathrm{d}\Omega, \qquad \qquad \delta_{\alpha} H = \alpha, \quad \delta_{\bm\beta} H = \bm\beta,
$$
where $\delta_\alpha$ and $\delta_\beta$ are the variational derivative with respect to the state variables. The power exchanged through the domain boundary corresponds to the time derivative of the energy
\begin{equation}
	\label{eq:wave_power_balance}
\begin{aligned}
\dot{H}(\alpha,\bm\beta) &= \int_\Omega \delta_\alpha H \cdot \partial_t \alpha + \delta_{\bm\beta} H \cdot \partial_t \bm{\beta} \; \mathrm{d}\Omega, \\
&= \int_\Omega \alpha \cdot \div \bm{\beta} + \bm\beta \cdot \grad {\alpha} \; \mathrm{d}\Omega, \\
&= \int_\Omega \div(\alpha \; \bm{\beta}) \; \mathrm{d}\Omega = \int_{\partial\Omega} \alpha \; \bm{\beta} \cdot \bm{n} \; \mathrm{d}\Gamma,
\end{aligned}
\end{equation}
where $\bm{n}$ is the outward normal to the boundary. The final expression pairs the Dirichlet condition with the Neumann boundary condition. It corresponds to a passive power balance of the form
$$
\dot{H} = \int_{\partial\Omega} u y\; \mathrm{d}\Gamma
$$
where variable ${y}$ is the power conjugated output to the input and corresponds to the Neumann boundary condition, 
\begin{equation}
    y := \bm{\beta} \cdot \bm{n}|_{\partial\Omega} = \grad \phi \cdot \bm{n}|_{\partial\Omega}  
\end{equation}
System \eqref{eq:mix_wave} is an example of a port-Hamiltonian system. The underlying geometrical structure is the Stokes-Dirac structure, an infinite dimensional generalization of Dirac manifolds introduced by Courant \cite{courant1990}.
}

\subsection{An abstract setting for linear port-Hamiltonian hyperbolic systems}

Consider a domain $\Omega\subset \mathbb{R}^d$, $d\in\{1,2,3\}$ and a partition of its boundary $\partial\Omega=\overline{\Gamma}_1\cup\overline{\Gamma}_2$, such that $\Gamma_1\cap\Gamma_2=\emptyset$. Let $\bm{x}=\{\bm{\alpha},\bm{\beta}\}$ be the state variables. The states at a given time are expected to be square integrable functions $\bm{\alpha}(t) \in L^2(\Omega; \mathbb{A})$, $\bm{\beta}(t) \in L^2(\Omega; \mathbb{B})$ taking values in the vector spaces $\mathbb{A}, \; \mathbb{B}$.

To define the dynamics of the system, an unbounded differential operator $\mathcal{L} :L^2(\Omega; \mathbb{A})\rightarrow L^2(\Omega; \mathbb{B})$ is introduced. For a given differential it is possible to define its formal adjoint by means of the integration by parts formula.
\begin{definition}[Formal Adjoint]\label{def:adjoint}
    Let $\bm{u}\in C_0^\infty(\Omega,\mathbb{A})$ and $\bm{v}\in C_0^\infty(\Omega,\mathbb{B})$ smooth variables on $\Omega$, and $\mathcal{L}$ be the differential operator $\mathcal{L}:L^2(\Omega; \mathbb{A})\rightarrow L^2(\Omega; \mathbb{B})$. The formal adjoint of $\mathcal{L}$ is than $\mathcal{L^*}:L^2(\Omega; \mathbb{B})\rightarrow L^2(\Omega; \mathbb{A})$ defined by the relation
    \begin{equation}
        \innerproduct{\mathcal{L}\bm{u}}{\bm{v}} = \innerproduct{\bm{u}}{\mathcal{L}^*\bm{v}}.
    \end{equation}
    where the inner product of two functions is denoted by $\innerproduct{f}{g} = \int_\Omega f \cdot g \, d{\Omega}$.
\end{definition}
The differential operator $\mathcal{L}$ and its formal adjoint $\mathcal{L}^*$ give rise to the Hilbert spaces $H^\mathcal{L}$ and $H^{\mathcal{L}^*}$, which are specified as
\begin{equation}
    \begin{split}
        H^\mathcal{L}(\Omega) &= \{\bm{u} \in L^2(\Omega,\mathbb{A})|\; \mathcal{L}\bm{u} \in L^2(\Omega,\mathbb{B})\}, \\
        H^{\mathcal{L}^*}(\Omega) &= \{\bm{v} \in L^2(\Omega; \mathbb{B})|\; \mathcal{L}^*\bm{v} \in L^2(\Omega; \mathbb{A})\}.
    \end{split}
\label{eq:Spaces-H^L-H^L*}
\end{equation}
The formal adjoint definition does not account for boundary terms arising from the integration by parts formula. They are introduced by means of the next assumption. 
\begin{assumption}[Abstract integration by parts]
    Let $\bm{e}_\alpha \in {H}^{\mathcal{L}}(\Omega), \; \bm{e}_\beta \in {H}^{\mathcal{L}^*}(\Omega)$. Then the following integration by parts formula is assumed to hold
    \begin{equation}\label{eq:abstract_intbyparts}
        \innerproduct{\mathcal{L}\bm{e}_\alpha}{\bm{e}_\beta} = \innerproduct{\bm{e}_\alpha}{\mathcal{L}^*\bm{e}_\beta} + \boundary{\mathcal{T}_{\alpha}\bm{e}_\alpha}{\mathcal{T}_{\beta}\bm{e}_\beta},
    \end{equation}  
for appropriate trace operators \revone{$\mathcal{T}_\alpha$ and $\mathcal{T_\beta}$}, where $\boundary{f}{g} = \int_{\partial\Omega} f \cdot g \, d\Gamma$ denotes the inner product over the boundary.
\end{assumption}
\begin{example}[Gradient and divergence operators]
    Let $\mathcal{L}:= \mathrm{grad}$ be the gradient and $\mathcal{L}^* = -\mathrm{div}$ be the negative of the divergence. Let $f \in H^1(\Omega), \; \bm{v} \in H^{\mathrm{div}}(\Omega)$ be a scalar and a vector function. The integration by parts states that the inner product with the gradient can be written as
    \begin{equation}
        \innerproduct{\grad f}{\bm{v}} = - \innerproduct{f}{\div\bm{v}} + \boundary{f}{\bm{v} \cdot \bm{n}}.
    \end{equation}
    In this case the trace operators correspond to the Dirichlet trace and the normal trace.
\end{example}
In this work we focus on conservation laws describing wave propagation phenomena in Hamiltonian form
\begin{equation}
\begin{pmatrix}
    \partial_t{\bm{\alpha}} \\
    \partial_t{\bm{\beta}}
\end{pmatrix} = 
\begin{bmatrix}
    0 & -\mathcal{L}^* \\
    \mathcal{L} & 0
\end{bmatrix}
\begin{pmatrix}
    \delta_{\bm{\alpha}} H\\
    \delta_{\bm{\beta}} H
\end{pmatrix}.
\label{eq:energy_co-energy-system}
\end{equation}
where $H$ is the Hamiltonian and $\delta_{\bm{\alpha}} H,   \delta_{\bm{\beta}} H$ its variational derivative with respect to the state variables. In this work we restrict our attention to linear wave propagation phenomena described by the Hamiltonian formalism. The linearity of the system translates into a quadratic Hamiltonian.
\begin{assumption}[Quadratic Hamiltonian]
The Hamiltonian is assumed to \revthree{take the following form}
\begin{equation}
    H = \frac{1}{2}\innerproduct{\bm{\alpha}}{\mathcal{Q}_\alpha\bm{\alpha}} + \frac{1}{2}\innerproduct{\bm{\beta}}{\mathcal{Q}_\beta\bm{\beta}},
\label{eq:linearHamiltonian}
\end{equation}
The operators $\mathcal{Q}_\alpha$ and $\mathcal{Q}_\beta$ are symmetric and positive (and therefore invertible).    
\end{assumption}
 The variational derivative of the Hamiltonian (also called co-energy variables) is evaluated as follows \cite{marsdenratiu1999}
\begin{equation}
        \bm{e}_\alpha := \frac{\delta H}{\delta \bm{\alpha}} = \mathcal{Q}_\alpha\bm{\alpha}, \qquad\qquad  \bm{e}_\beta := \frac{\delta H}{\delta \bm{\beta}} = \mathcal{Q}_\beta\bm{\beta}.
\label{eq:co-energy-variables}
\end{equation}
Given \eqref{eq:co-energy-variables}, the state variables are related to the co-energy variables by 
$$\bm{\alpha}=\mathcal{M}_\alpha\bm{e}_\alpha, \qquad \bm{\beta} = \mathcal{M}_\beta \bm{e}_\beta,$$
where $\mathcal{M}_\alpha:=\mathcal{Q}_\alpha^{-1}, \; \mathcal{M}_\beta :=\mathcal{Q}_\beta^{-1}$. The Hamiltonian \eqref{eq:linearHamiltonian} can be expressed in terms of co-energy variables as
\begin{equation}
    H = \frac{1}{2}\innerproduct{\bm{e}_\alpha}{\mathcal{M}_\alpha\bm{e}_\alpha} + \frac{1}{2}\innerproduct{\bm{e}_\beta}{\mathcal{M}_\beta\bm{e}_\beta}.
\label{eq:co-energyHamiltonian}
\end{equation}
The system can be equivalently rewritten in terms of the coenergy variables including mixed boundary conditions as follows
\begin{equation}\label{eq:co-energy-system}
    \begin{split}
    \begin{bmatrix}
        \mathcal{M}_\alpha & 0 \\
        0 & \mathcal{M}_\beta
    \end{bmatrix}
    \begin{pmatrix}
        \partial_t \bm{e}_{\alpha}\\
        \partial_t \bm{e}_{\beta}
    \end{pmatrix} &= 
    \begin{bmatrix}
        0 & -\mathcal{L}^* \\
        \mathcal{L} & 0
    \end{bmatrix}
    \begin{pmatrix}
        \bm{e}_\alpha \\
        \bm{e}_\beta
    \end{pmatrix},  \qquad 
    \begin{aligned}
        \bm{e}_\alpha \in H^\mathcal{L}(\Omega), \\
        \bm{e}_\beta \in H^\mathcal{L^*}(\Omega),
    \end{aligned}\\
    \begin{pmatrix}
    \bm{u}_{\partial, 1} \\
    \bm{u}_{\partial, 2} \\    
    \end{pmatrix}
     &= 
     \begin{bmatrix}
    \mathcal{T}_\alpha\vert_{\Gamma_1} & 0 \\
    0 & \mathcal{T}_\beta \vert_{\Gamma_2}
    \end{bmatrix}
    \begin{pmatrix}
        \bm{e}_\alpha \\
        \bm{e}_\beta
    \end{pmatrix}, \\
    \begin{pmatrix}
    \bm{y}_{\partial, 1} \\
    \bm{y}_{\partial, 2} \\    
    \end{pmatrix}
     &= 
     \begin{bmatrix}
    0 & \mathcal{T}_\beta\vert_{\Gamma_1} \\
    \mathcal{T}_\alpha \vert_{\Gamma_2} & 0
    \end{bmatrix}
    \begin{pmatrix}
        \bm{e}_\alpha \\
        \bm{e}_\beta
    \end{pmatrix}.
    \end{split}
\end{equation}
\revone{Variable $\bm{u}_{\partial, i}$ are the inputs. In the terminology of partial differential equations, they correspond to the boundary data. In the following we will use input as a synonym for boundary conditions to stress the fact that these quantities are not fixed a priori but are the result of the system interacting with the environment. The variables $\bm{y}_{\partial, i}$ correspond to the conjugate outputs, i.e. the dual variable to the corresponding input $\bm{u}_{\partial, i}$ with respect to the power balance
\begin{equation*}
    \dot{H} = \boundary[\Gamma_1]{\bm{u}_{\partial, 1}}{\bm{y}_{\partial, 1}} + \boundary[\Gamma_2]{\bm{u}_{\partial, 2}}{\bm{y}_{\partial, 2}}.
\end{equation*}
Examples of this duality are force and velocity in mechanics or voltage and current in electromagnetism.}  The notation $\mathcal{T}_\alpha\vert_{\Gamma_i}, \; \mathcal{T}_\beta \vert_{\Gamma_i}, \; i=\{1, 2\}$ denotes the restriction of the trace operators to a given subpartition of the boundary.

\begin{remark}[Equivalence with Lagrangian dynamics]
    The presented Hamiltonian formulation can be deduced from the least action principle and is equivalent to a Lagrangian formulation \cite{brugnoli2024OnPhenomena}.
\end{remark}

\subsection{\revone{Conforming finite element discretization of port-Hamiltonian systems under mixed boundary conditions}}
The discretization of problem \eqref{eq:co-energy-system} is detailed in \cite{joly2003}, where its primal-dual structure is highlighted. Therein however the point of view of Hilbert complexes is not considered and this mathematical structure is important for port-Hamiltonian systems. \revone{We will here consider a finite element formulation that respects the Hilbert complex structure. Furthermore, we detail the different numerical treatment of input and output variables in the case of mixed boundary conditions. In a classical monolithic discretization, one input variable enters the system via integration by parts and for this reason is called natural. The second input variable has to be enforced in a strong way and it is typically called essential.}

The weak formulation can now be obtained by applying the test function $\bm{v}=\{\bm{v}_\alpha,\bm{v}_\beta\}$ and integrating over $\Omega$ to end up with
\begin{equation}
    \begin{split}
        \innerproduct{\bm{v}_\alpha}{\mathcal{M}_\alpha \partial_t\bm{e}_{\alpha}} &= -\innerproduct{\bm{v}_\alpha}{\mathcal{L}^*\bm{e}_\beta}, \\
        \innerproduct{\bm{v}_\beta}{\mathcal{M}_\beta \partial_t\bm{e}_{\beta}} &= \innerproduct{\bm{v}_\beta}{\mathcal{L}\bm{e}_\alpha}.
    \end{split}
\end{equation}

Given the abstract integration by parts formula \eqref{eq:abstract_intbyparts}, two possibilities arise. One can either integrate by parts the first line or the second. Depending on the choice, two dynamical systems are obtained. \revone{These two systems differ in the way they treat boundary conditions. In the first system $\bm{u}_{\partial, 1}$ is a natural boundary condition, whereas in the second system the natural boundary condition is $\bm{u}_{\partial, 2}$. To explain the classical Galerkin discretization in the case of mixed boundary conditions, Lagrange multipliers $\bm{\lambda}_{\partial,i}$ will be used.
}
\paragraph{System 1: \revone{natural imposition of $\bm{u}_{\partial, 1}$, essential imposition of $\bm{u}_{\partial, 2}$}} 
If the second line is integrated by parts, the weak formulation reads: find $\bm{e}_{\alpha}  \in L^2(\Omega; \mathbb{A})$ and $\bm{e}_{\beta} \in H^{\mathcal{L}^*}(\Omega)$ such that
\begin{equation}
    \begin{aligned}
        \innerproduct{\bm{v}_{\alpha}}{\mathcal{M}_\alpha \partial_t\bm{e}_{\alpha}} &= -\innerproduct{\bm{v}_{\alpha}}{\mathcal{L}^*\bm{e}_{\beta}}, \\
        \innerproduct{\bm{v}_{\beta}}{\mathcal{M}_\beta \partial_t\bm{e}_{\beta}} &= \innerproduct{\mathcal{L}^*\bm{v}_{\beta}}{\bm{e}_{\alpha}} + \boundary[\Gamma_1]{\mathcal{T}_{\beta}\bm{v}_{\beta}}{\bm{u}_{\partial, 1}} + \boundary[\Gamma_2]{\mathcal{T}_{\beta}\bm{v}_{\beta}}{\bm{\lambda}_{\partial, 2}}, \\
        \mathcal{T}_{\beta}|_{\Gamma_2} \bm{e}_{\beta} &= \bm{u}_{\partial, 2}, \\
        \bm{y}_{\partial, 1} &= \mathcal{T}_\beta|_{\Gamma_1} \bm{e}_{\beta}, \\
        \bm{y}_{\partial, 2} &= \bm{\lambda}_{\partial, 2}.
    \end{aligned} \qquad 
    \begin{aligned}
        \text{for all } \bm{v}_{\alpha}  \in L^2(\Omega; \mathbb{A}), \\
        \text{for all } \bm{v}_{\beta} \in H^{\mathcal{L}^*}(\Omega), \\
        {} \\
        {} \\
        {}
    \end{aligned}
    \label{eq:weak_conforming_L}
\end{equation}
\revone{The essential input $\bm{u}_{\partial, 2}$ and the output $\bm{y}_{\partial, 1}$ are not evaluated weakly, but taken to be the trace of the associated state variable.}
\paragraph{System 2: \revone{natural imposition of $\bm{u}_{\partial, 2}$, essential imposition of $\bm{u}_{\partial, 1}$}}
If the first line is integrated by parts, the following system is obtained: find $\bm{e}_{\alpha} \in H^\mathcal{L}(\Omega), \; \bm{e}_\beta \in L^2(\Omega; \mathbb{B})$ such that
\begin{equation}
    \begin{aligned}
        \innerproduct{\bm{v}_{\alpha}}{\mathcal{M}_\alpha \partial_t\bm{e}_{\alpha}} &= -\innerproduct{\mathcal{L}\bm{v}_{\alpha}}{\bm{e}_{\beta}} + \boundary[\Gamma_1]{\mathcal{T}_{\alpha}\bm{v}_{\alpha}}{\bm{\lambda}_{\partial, 1}} + \boundary[\Gamma_2]{\mathcal{T}_{\alpha}\bm{v}_{\alpha}}{\bm{u}_{\partial, 2}}, \\
        \innerproduct{\bm{v}_{\beta}}{\mathcal{M}_\beta \partial_t\bm{e}_{\beta}} &= \innerproduct{\bm{v}_{\beta}}{\mathcal{L}\bm{e}_{\alpha}}, \\
        \mathcal{T}_{\alpha}|_{\Gamma_1}\bm{e}_{\alpha} &= \bm{u}_{\partial, 1}, \\
        \bm{y}_{\partial, 1} &= \bm{\lambda}_{\partial,1}. \\
        \bm{y}_{\partial, 2} &= \mathcal{T}_\alpha|_{\Gamma_2} \bm{e}_{\alpha}. \\
    \end{aligned} \qquad
    \begin{aligned}
        \text{for all } \bm{v}_{\alpha}  \in H^\mathcal{L}(\Omega), \\
        \text{for all } \bm{v}_{\beta} \in L^2(\Omega; \mathbb{B}), \\
        {} \\ {} \\ {}
    \end{aligned}
    \label{eq:weak_conforming_L*}
\end{equation}

\paragraph{Finite dimensional representation of the variables}
The two systems should not be discretized in the same manner, as different differential operators may arise in the weak formulations \eqref{eq:weak_conforming_L}, \eqref{eq:weak_conforming_L*}. \revone{For this reason we now refer to the variables of formulation $i$ using an appropriate subscript $i=\{1, 2\}$. Consider a finite element Galerkin approximation of the test, trial and boundary input functions. For the two systems they are
\begin{equation}
\begin{aligned}
\bm{v}_{\alpha, i} &\approx \sum_{k=1}^{n_{\alpha, i}}\chi_{\alpha, i}^k(\bm{x})v_{\alpha, i}^k, & \bm{v}_{\beta, i} &\approx \sum_{k=1}^{n_{\beta, i}}\chi_{\beta, i}^k(\bm{x})v_{\beta, i}^k, & \bm{u}_{\partial, i} &\approx \sum_{k=1}^{n_{\partial, i}}\chi_{\partial, i}^k(\bm{x})u_{\partial, i}^k(t), \\
\bm{e}_{\alpha, i} &\approx \sum_{k=1}^{n_{\alpha, i}}\chi_{\alpha, i}^k(\bm{x})e_{\alpha, i}^k(t), & \bm{e}_{\beta, i} &\approx \sum_{k=1}^{n_{\beta, i}}\chi_{\beta, i}^k(\bm{x})e_{\beta, i}^k(t), & \bm{\lambda}_{\partial, i} &\approx \sum_{k=1}^{n_{\partial, i}}\chi_{\partial, i}^k(\bm{x})\lambda_{\partial, i}^k(t), 
\end{aligned}
\label{eq:BasisFunction}
\end{equation}
where $\chi_{\alpha, i}, \; \chi_{\beta, i}, \; \chi_{\partial, i}$ the basis functions for the finite element spaces and $n_{\alpha, i}, \; n_{\beta, i}$ are the number of degrees of freedom associated to variable $\bm{e}_{\alpha}, \bm{e}_{\beta}$ on the domain $\Omega_i$.} The finite element spaces associated to the state variables $\bm{e}_{\alpha, i}, \; \bm{e}_{\beta, i}$ are spanned by the basis functions
\begin{equation}
    V_{\alpha, i} = \mathrm{span}\{\chi_{\alpha, i}\}, \qquad V_{\beta, i} = \mathrm{span}\{\chi_{\beta, i}\},  \qquad i=\{1, 2\}.
\end{equation}
\revone{Considering the fact that the weak formulations in Eqs. \eqref{eq:weak_conforming_L} and \eqref{eq:weak_conforming_L*} are conforming, the spaces verify the following inclusions
\begin{equation}
    V_{\alpha, 1} \subset L^2(\Omega; \mathbb{A}), \quad V_{\beta, 1} \subset H^{\mathcal{L}^*}(\Omega), \qquad V_{\alpha, 2} \subset H^{\mathcal{L}}(\Omega), \quad V_{\beta, 1} \subset L^2(\Omega; \mathbb{B}).
\end{equation}
}

The guiding principle behind the choice of $V_{\alpha, i}, \; V_{\beta, i}$ is that of a Hilbert complex. 
\begin{definition}[Hilbert Complex]
A \textit{Hilbert complex} is a sequence $\{H^k, \mathcal{L}^k\}_{k \in \mathbb{Z}}$ where:
\begin{itemize}
    \item $H^k$ are Hilbert spaces,
    \item $\mathcal{L}^k: H^k \to H^{k+1}$ are bounded linear operators,
    \item $\mathcal{L}^{k+1} \circ \mathcal{L}^k = 0$ for all $k \in \mathbb{Z}$.
\end{itemize}
\end{definition}
Given a Hilbert complex it is possible to define the adjoint complex by using the definition of an adjoint operator.

\begin{example}[de-Rham Complex]
    One important example of a Hilbert complex that will be considered in Sec.~\ref{sec:examples} is the de Rham complex
\begin{equation}
\begin{aligned}
H^1(\Omega) \xrightarrow{\grad} H^{\curl}(\Omega) \xrightarrow{\curl} H^{\div}(\Omega) \xrightarrow{\div} L^2(\Omega)
\end{aligned}
\end{equation}
The adjoint complex reads
\begin{equation}
\begin{aligned}
L^2(\Omega) \xleftarrow{\div} \mathring{H}^{\div}(\Omega) \xleftarrow{\curl} \mathring{H}^{\curl}(\Omega) \xleftarrow{\grad} \mathring{H}^1(\Omega)
\end{aligned}
\end{equation}
where the Hilbert spaces in the adjoint complex include homogeneous boundary conditions.   
\end{example}
 
 Finite element spaces $V_{\alpha, i}, \; V_{\beta, i}$ are chosen from a finite dimensional subcomplex.
\begin{definition}[Hilbert Subcomplex]
Given a Hilbert complex $\{H^k, \mathcal{L}^k\}_{k \in \mathbb{Z}}$, a \textit{subcomplex} is a sequence of closed subspaces $\{V^k \subseteq H^k\}_{k \in \mathbb{Z}}$ such that:
\begin{itemize}
    \item $\mathcal{L}^k(V^k) \subseteq V^{k+1}$ for all $k \in \mathbb{Z}$,
    \item $V^k$ is a closed linear subspace of $H^k$,
    \item The restriction of $\mathcal{L}^k$ to $V^k$ maps $V^k$ to $V^{k+1}$.
\end{itemize}
\end{definition}
\revone{In order to obtain a finite element subcomplex, the finite element spaces are selected in such a way that
\begin{equation}\label{eq:spaces_inclusion}
    \mathcal{L}^*(V_{\beta, 1}) \subset V_{\alpha, 1}, \qquad \qquad \mathcal{L}(V_{\alpha, 2}) \subset V_{\beta, 2}.
\end{equation}
This means that the spaces used for the discretization form two complexes
$$
V_{\alpha, 2} \xrightarrow{\mathcal{L}} V_{\beta, 2}, \qquad \qquad V_{\beta, 1} \xrightarrow{\mathcal{L}^*} V_{\alpha, 1}.
$$
Spaces satisfying such an inclusion can be constructed in several manners \cite{arnold2006finite,bochev2006principles,palha2014physics}. The rationale behind the choice for the boundary spaces $V_{\partial, i} = \mathrm{span}\{\chi_{\partial, i}\}$ is important as it establishes a connection between the two formulations. As $V_{\partial, i}$ are trace spaces, their elements can be taken to be the restriction to the boundary subpartitions $\Gamma_i$ of the spaces $V_{\beta, 1} \subset H^{\mathcal{L}^*}(\Omega)$ and $V_{\alpha, 2} \subset H^{\mathcal{L}}(\Omega)$ 
\begin{equation}\label{eq:boundary_basis_boundary} 
\begin{aligned}
    \mathrm{span}\{\chi_{\partial, 1}\} &:= \mathrm{span} \{\mathrm{\mathcal{T}_\alpha|_{\Gamma_1}\chi_{\alpha, 2}}\}, \\
    \mathrm{span}\{\chi_{\partial, 2}\} &:= \mathrm{span} \{\mathcal{T}_\beta|_{\Gamma_2} \chi_{\beta, 1}\}.
\end{aligned}
\end{equation}
}

\begin{remark}[Equivalence with the second order formulation]
Because of the inclusions \eqref{eq:spaces_inclusion}, the mixed formulations coincides therefore with the second order formulation in time and space \cite{joly2003}.
\end{remark}

\revone{
\paragraph{Algebraic realization for System 1}
In this case the formulation \eqref{eq:weak_conforming_L} is converted into the following differential algebraic system
\begin{equation}
    \begin{aligned}
        \begin{bmatrix}
            \mathbf{M}_{\alpha, 1} & 0 & 0\\
            0 & \mathbf{M}_{\beta, 1} & 0 & \\
            0 & 0 & 0
        \end{bmatrix}
        \odv{}{t}
        \begin{pmatrix}
            \mathbf{e}_{\alpha, 1} \\
            \mathbf{e}_{\beta, 1} \\
            \bm{\lambda}_{\partial, 2}
        \end{pmatrix} &= 
        \begin{bmatrix}
            0 & -\mathbf{D}_{\mathcal{L}^*} & 0 \\
             \mathbf{D}_{\mathcal{L}^*}^\top & 0 & (\mathbf{T}_\beta^{\Gamma_2})^\top \mathbf{M}^{\Gamma_2} \\
             0 & -\mathbf{T}_{\beta}^{\Gamma_2} & 0 
        \end{bmatrix}
        \begin{pmatrix}
            \mathbf{e}_{\alpha, 1} \\
            \mathbf{e}_{\beta, 1} \\
            \bm{\lambda}_{\partial, 2}            
        \end{pmatrix} + 
        \begin{bmatrix}
            0 & 0 \\
            \mathbf{B}_\beta^{\Gamma_1} & 0 \\
            0 & \mathbf{I}
        \end{bmatrix} 
        \begin{pmatrix}
            \mathbf{u}_{\partial, 1} \\
            \mathbf{u}_{\partial, 2}
        \end{pmatrix}, \\
        \begin{pmatrix}
            \mathbf{y}_{\partial, 1} \\
            \mathbf{y}_{\partial, 2}
        \end{pmatrix} &= \begin{bmatrix}
            0 & \mathbf{T}_\beta^{\Gamma_1} & 0 \\ 
            0 & 0 & \mathbf{I}
        \end{bmatrix} \begin{pmatrix}
            \mathbf{e}_{\alpha, 1} \\
            \mathbf{e}_{\beta, 1}, \\
            \bm{\lambda}_{\partial, 2}
        \end{pmatrix},
    \end{aligned}
\label{eq:systemL*}
\end{equation}
where the matrices arising from the weak formulation are defined by
\begin{equation}
\begin{aligned}
    &[\mathbf{M}_{\alpha, 1}]_{mn} =\innerproduct{\chi_{\alpha, 1}^m}{\mathcal{M}_\alpha \chi_{\alpha, 1}^n}, \\
    &[\mathbf{M}_{\beta, 1}]_{pq} = \innerproduct{\chi_{\beta, 1}^p}{\mathcal{M}_\beta \chi_{\beta, 1}^q}, \\
\end{aligned}\qquad
    \begin{aligned}
    &[\mathbf{D}_{\mathcal{L}^*}]_{mp} = \innerproduct{\chi_{\alpha, 1}^m}{\mathcal{L}^*\chi_{\beta, 1}^p}, \\
    &[\mathbf{M}^{\Gamma_2}]_{rs} = \boundary[\Gamma_2]{\chi_{\partial, 2}^r}{\chi_{\partial, 2}^s},
    \end{aligned} \qquad 
    \begin{aligned}
        &[\mathbf{B}_\beta^{\Gamma_1}]_{pl} = \boundary[\Gamma_1]{\mathcal{T}_\beta\chi_{\beta, 1}^p}{\chi_{\partial, 1}^l}, \\
        &{}
    \end{aligned}
\label{eq:Operational-matrices_1}
\end{equation}
where $(m,n) \in \{1, \dots, n_{\alpha, 1}\}, \; (p,q) \in\{1, \dots, n_{\beta, 1}\}, \; (r,s) \in \{1,\dots, n_{\partial, 2}\}, \; l \in \{1,\dots, n_{\partial, 1}\}$. The trace matrix is a Boolean matrix that localizes the degrees of freedom lying on the boundary
\begin{equation}
        [\mathbf{T}_{\beta}^{\Gamma_i}]_{kp} = \begin{cases}
            1,  \quad \text{if} \quad \mathcal{T}_\beta\chi_{\beta, 1}^p \not\equiv 0 \quad \text{on } \Gamma_i, \quad i=\{1,2\}, \\
            0,  \quad \text{otherwise},
        \end{cases}
\end{equation}
where $k=1, \dots, \mathrm{dim} \{\mathcal{T}_\beta\chi_{\beta, 1}^i  \not\equiv 0\}_{i=1}^{n_{\beta, 1}}$ counts over the basis function that lie on the boundary subpartition~$\Gamma_i$. The matrix $\mathbf{B}_\beta^{\Gamma_1}$ can be decomposed using the trace matrix as follows
$$
\mathbf{B}_\beta^{\Gamma_1} = (\bm{\Psi}^{\Gamma_1} \mathbf{T}_\beta^{\Gamma_1})^\top, \quad \text{where } \quad [\bm{\Psi}^{\Gamma_1}]_{lj}:= \boundary[\Gamma_1]{\chi_{\partial, 1}^l}{\chi_{\partial, 2}^j}.
$$
\paragraph{Algebraic realization for System 2}
In this case the formulation \eqref{eq:weak_conforming_L*} is converted into the following differential algebraic system
\begin{equation}
\begin{aligned}
        \begin{bmatrix}
            \mathbf{M}_{\alpha, 2} & 0 & 0\\
            0 & \mathbf{M}_{\beta, 2}  & 0 \\
            0 & 0 & 0
        \end{bmatrix} 
        \odv{}{t}
        \begin{pmatrix}
            \mathbf{e}_{\alpha, 2} \\
            \mathbf{e}_{\beta, 2} \\
            \bm{\lambda}_{\partial, 1}
        \end{pmatrix} &=
        \begin{bmatrix}
            0 & -\mathbf{D}_\mathcal{L}^\top & (\mathbf{T}_{\alpha}^{\Gamma_1})^\top \mathbf{M}^{\Gamma_1} \\
            \mathbf{D}_\mathcal{L} & 0 & 0 \\
            -\mathbf{T}_{\alpha}^{\Gamma_1} & 0 & 0
        \end{bmatrix}
        \begin{pmatrix}
            \mathbf{e}_{\alpha, 2} \\
            \mathbf{e}_{\beta, 2} \\
            \bm{\lambda}_{\partial, 1}
        \end{pmatrix} + 
        \begin{bmatrix}
            0 & \mathbf{B}_\alpha^{\Gamma_2} \\
            0 & 0  \\
            \mathbf{I}& 0
        \end{bmatrix}
        \begin{pmatrix}
            \mathbf{u}_{\partial, 1} \\
            \mathbf{u}_{\partial, 2}
        \end{pmatrix}, \\
        \begin{pmatrix}
            \mathbf{y}_{\partial, 1} \\
            \mathbf{y}_{\partial, 2}
        \end{pmatrix} &= 
        \begin{bmatrix}
            0 & 0 & \mathbf{I} \\
            \mathbf{T}_\alpha^{\Gamma_2} & 0 & 0 \\
        \end{bmatrix}
        \begin{pmatrix}
            \mathbf{e}_{\alpha, 2} \\
            \mathbf{e}_{\beta, 2} \\
            \bm{\lambda}_1
        \end{pmatrix}.
\end{aligned}
\label{eq:systemL}
\end{equation}
The matrix components are obtained as follows
\begin{equation}
\begin{aligned}
&[\mathbf{M}_{\alpha, 2}]_{mn} = \innerproduct{\chi_{\alpha, 2}^m}{\mathcal{M}_\alpha \chi_{\alpha, 2}^n}, \\
&[\mathbf{M}_{\beta, 2}]_{pq} = \innerproduct{\chi_{\beta, 2}^p}{\mathcal{M}_\beta \chi_{\beta, 2}^q}, \\
\end{aligned}\qquad
\begin{aligned}
&[\mathbf{D}_\mathcal{L}]_{pm} = \innerproduct{\chi_{\beta, 2}^p}{\mathcal{L}\chi_{\alpha, 2}^m}, \\
&[\mathbf{M}^{\Gamma_1}]_{rs} = \boundary[\Gamma_1]{\chi_{\partial, 1}^i}{\chi_{\partial, 1}^k}, \\
\end{aligned}
\qquad
\begin{aligned}
&[\mathbf{B}_\alpha^{\Gamma_2}]_{pl} = \boundary[\Gamma_2]{\mathcal{T}_{\alpha}\chi_{\alpha, 2}^p}{\chi_{\partial, 2}^l}, \\
{}
\end{aligned}
\label{eq:Operational-matrices_2}
\end{equation}
where $(m,n) \in \{1, \dots, n_{\alpha, 2}\}, \; (p,q) \in\{1, \dots, n_{\beta, 2}\}, \; (r,s) \in \{1,\dots, n_{\partial, 1}\}, \; l \in \{1,\dots, n_{\partial, 2}\}$. For this system the trace matrix selects the degrees of freedom for the variable $\mathbf{e}_{\beta, 2}$
\begin{equation}
        [\mathbf{T}_{\alpha}^{\Gamma_i}]_{ki} = \begin{cases}
            1,  \quad \text{if} \quad \mathcal{T}_\alpha \chi_{\alpha, 2}^i \not\equiv 0, \quad \text{on } \partial\Omega,\\
            0,  \quad \text{otherwise}
        \end{cases} 
\end{equation}
Once again the control input matrix $\mathbf{B}_\alpha^{\Gamma_2}$ can be factorized using the trace matrix as follows
$$
\mathbf{B}_\alpha^{\Gamma_2} = (\mathbf{T}_\alpha^{\Gamma_2})^\top\bm{\Psi}^{\Gamma_2}, \quad \text{where } \quad [\bm{\Psi}^{\Gamma_2}]_{lj}:= \boundary[\Gamma_2]{\chi_{\partial, 1}^l}{\chi_{\partial, 2}^j}.
$$
} 
\section{Domain Decomposition for mixed boundary conditions}\label{sec:domain_decomposition}
\begin{figure}[tbh]
	\begin{minipage}[t]{0.4\linewidth}
		\centering
        \includegraphics[width=0.95\textwidth]{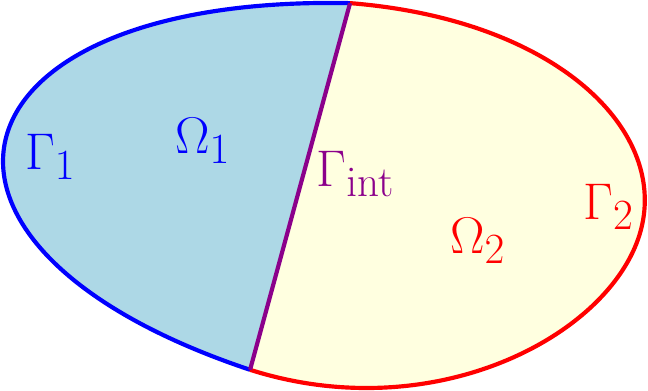}
		\caption{Splitting of the domain.}
		\label{fig:dom_part}
	\end{minipage}
	\hspace{0.5cm}
	\begin{minipage}[t]{0.55\linewidth}
		\centering
\includegraphics[width=0.95\textwidth]{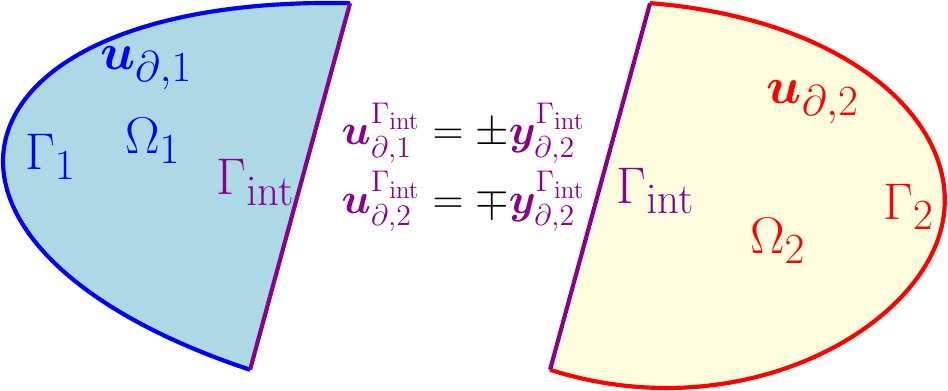}
		\caption{Interconnection at the interface $\Gamma_{12}$.}
		\label{fig:dom_part_int}
	\end{minipage}
\end{figure}

\revone{It has been shown that in order to solve problem \eqref{eq:co-energy-system} numerically, either system \eqref{eq:systemL*} or system \eqref{eq:systemL} can be used. These systems are differential algebraic as the essential imposition of the boundary data leads to constraints imposed on the dynamics. The idea behind the domain decomposition approach is to introduce} an interface boundary $\Gamma_{\rm int}$ to split the domain $\Omega=\Omega_1\cup\Omega_2$, where $\Omega_1\cap\Omega_2=\emptyset$ holds, such that the boundaries of the subdomains are given as $\partial\Omega_1=\overline{\Gamma}_1\cup\overline{\Gamma}_{\rm int}$ and $\partial\Omega_2=\overline{\Gamma}_2\cup\overline{\Gamma}_{\rm int}$ (cf. Fig \ref{fig:dom_part}). This interface boundary is chosen freely. The idea of the discretization is to use both formulations \eqref{eq:weak_conforming_L} and \eqref{eq:weak_conforming_L*} concurrently to achieve natural boundary imposition for both boundary \revone{inputs}. This means applying the formulations \eqref{eq:weak_conforming_L} to the $\Omega_1$ subdomain and \eqref{eq:weak_conforming_L*} to $\Omega_2$. To ensure proper coupling on the interface $\Gamma_{\rm int}$ consider the inputs and outputs from the port-Hamiltonian systems \eqref{eq:pH-system-L*} and \eqref{eq:pH-system-L}. The boundary inputs and outputs for $\Omega_1$ include the boundary condition for the problem and the input along the interconnection boundary 

\begin{equation}\label{eq:input_ouput_1}
    \begin{pmatrix}
        \bm{u}_{\partial, 1} \\
        \bm{u}_{\partial, 1}^{\Gamma_{\rm int}}
    \end{pmatrix} = 
    \begin{bmatrix}
        \mathcal{T}_{\alpha}\vert_{\Gamma_1} & 0 \\
        \mathcal{T}_{\alpha}\vert_{\Gamma_{\rm int}} & 0
    \end{bmatrix}
    \begin{pmatrix}
        \bm{e}_\alpha \\
        \bm{e}_\beta
    \end{pmatrix}, \qquad
    \begin{pmatrix}
        \bm{y}_{\partial,1} \\
        \bm{y}_{\partial,1}^{\Gamma_{\rm int}}
    \end{pmatrix} = 
    \begin{bmatrix}
        0 & \mathcal{T}_{\beta}\vert_{\Gamma_1} \\
        0 & \mathcal{T}_{\beta}\vert_{\Gamma_{\rm int}}
    \end{bmatrix}
    \begin{pmatrix}
        \bm{e}_\alpha \\
        \bm{e}_\beta
    \end{pmatrix}.
\end{equation}
For the $\Omega_2$ domain one input will be the actual boundary condition and a second input represent the exchange of information along the interface
\begin{equation}\label{eq:input_ouput_2}
    \begin{pmatrix}
        \bm{u}_{\partial,2} \\
        \bm{u}_{\partial,2}^{\Gamma_{\rm int}}
    \end{pmatrix} = 
    \begin{bmatrix}
        0 & \mathcal{T}_{\beta}\vert_{\Gamma_2} \\
        0 & \mathcal{T}_{\beta}\vert_{\Gamma_{\rm int}}
    \end{bmatrix}
    \begin{pmatrix}
        \bm{e}_\alpha \\
        \bm{e}_\beta
    \end{pmatrix}, \qquad
    \begin{pmatrix}
        \bm{y}_{\partial,2} \\
        \bm{y}_{\partial,2}^{\Gamma_{\rm int}}
    \end{pmatrix} = 
    \begin{bmatrix}
        \mathcal{T}_{\alpha}\vert_{\Gamma_2} & 0\\
        \mathcal{T}_{\alpha}\vert_{\Gamma_{\rm int}} & 0
    \end{bmatrix}
    \begin{pmatrix}
        \bm{e}_\alpha \\
        \bm{e}_\beta
    \end{pmatrix}.
\end{equation}

The coupling of the two domains takes place on $\Gamma_{\rm int}$ because the inputs and outputs are related by
\begin{equation}\label{eq:coupling}
    \begin{split}
        \bm{u}_{\partial,1}^{\Gamma_{\rm int}} &= \pm\bm{y}_{\partial,2}^{\Gamma_{\rm int}}, \\
        \bm{u}_{\partial,2}^{\Gamma_{\rm  int}} &= \mp\bm{y}_{\partial,1}^{\Gamma_{\rm int}},
    \end{split}
\end{equation}
as shown in Fig.~\ref{fig:dom_part_int}. The $\pm$ and $\mp$ are used due to opposite outward unit normals depending on the domain. \revone{An explicit incorporation of the boundary conditions} is achieved by integrating the $\mathcal{L}$ term by parts on $\Omega_1$, based on the weak form \eqref{eq:weak_conforming_L}, and the $\mathcal{L}^*$ on the $\Omega_2$ subdomain, based on the weak form \eqref{eq:weak_conforming_L*}. Consider the additivity of integral operator, the boundary term $\boundary{\mathcal{T}_{\partial,\beta}\bm{v}_\beta}{\bm{u}_\partial}$ from the $\Omega_1$ domain becomes

\begin{equation}
    \boundary[\partial\Omega_1]{\mathcal{T}_{\beta}\bm{v}_\beta}{\mathcal{T}_{\alpha}\bm{e}_\alpha} = \boundary[\Gamma_1]{\mathcal{T}_{\beta}\bm{v}_\beta}{\bm{u}_{\partial,1}} + \boundary[\Gamma_{\rm int}]{\mathcal{T}_{\beta}\bm{v}_\beta}{\bm{u}_{\partial,1}^{\Gamma_{\rm int}}},
\label{eq:boundary-split-Omega1}
\end{equation}
while for the $\Omega_2$ subdomain 
\begin{equation}
    \boundary[\partial\Omega_2]{\mathcal{T}_{\alpha}\bm{v}_\alpha}{\mathcal{T}_{\beta}\bm{e}_\beta} = \boundary[\Gamma_2]{\mathcal{T}_{\alpha}\bm{v}_\alpha}{\bm{u}_{\partial,2}} + \boundary[\Gamma_{\rm int}]{\mathcal{T}_{\alpha}\bm{v}_\alpha}{\bm{u}_{\partial,2}^{\Gamma_{\rm int}}}.
\label{boundary-split-Omega2}
\end{equation}
The weak formulation for $\Omega_1$ is to find $\bm{e}_\alpha \in L^2(\Omega_1; \mathbb{A}), \; \bm{e}_\beta \in H^{\mathcal{L^*}}(\Omega_1)$ such that $\forall \; \bm{v}_\alpha \in L^2(\Omega_2; \mathbb{A})$ and $\forall \; \bm{v}_\beta \in H^{\mathcal{L^*}}(\Omega_1)$ it holds
\begin{equation}
    \begin{aligned}
        \innerproduct[\Omega_1]{\bm{v}_\alpha}{\mathcal{M}_\alpha\partial_t\bm{e}_\alpha} &= - \innerproduct[\Omega_1]{\bm{v}_\alpha}{\mathcal{L}^*\bm{e}_\beta}, \\
        \innerproduct[\Omega_1]{\bm{v}_\alpha}{\mathcal{M}_\beta\partial_t\bm{e}_\beta} &= \innerproduct[\Omega_1]{\mathcal{L}^*\bm{v}_\beta}{\bm{e}_\alpha} + \boundary[\Gamma_1]{\mathcal{T}_{\beta}\bm{v}_\beta}{\bm{u}_{\partial,1}} + \boundary[\Gamma_{\rm int}]{\mathcal{T}_{\beta}\bm{v}_\beta}{\bm{u}_{\partial,1}^{\Gamma_{\rm int}}},
    \end{aligned}
\label{eq:weak-form-L*}
\end{equation}
where the boundary control and trace matrices are now restricted on the subpartitions  of the boundary $\Gamma_{\rm int}, \; $
 For the $\Omega_2$ subdomain with \ref{boundary-split-Omega2} find $\bm{e}_\alpha\in H^\mathcal{L}(\Omega_2), \;\bm{e}_\beta \in L^2(\Omega_2; \mathbb{B})$ that satisfy $\forall \; \bm{v}_\alpha \in H^{\mathcal{L}}(\Omega_2)$ and $\forall \; \bm{v}_\beta \in L^2(\Omega_2; \mathbb{B})$
\begin{equation}
    \begin{aligned}
        \innerproduct[\Omega_2]{\bm{v}_\alpha}{\mathcal{M}_\alpha\partial_t\bm{e}_\alpha} &= -\innerproduct[\Omega_2]{\mathcal{L}\bm{v}_\alpha}{\bm{e}_\beta} + \boundary[\Gamma_2]{\mathcal{T}_{\alpha}\bm{v}_\alpha}{\bm{u}_{\partial,2}} + \boundary[\Gamma_{\rm int}]{\mathcal{T}_{\alpha}\bm{v}_\alpha}{\bm{u}_{\partial,2}^{\Gamma_{\rm int}}}, \\
        \innerproduct[\Omega_2]{\bm{v}_\beta}{\mathcal{M}_\beta\partial_t\bm{e}_\beta} &= \innerproduct[\Omega_2]{\bm{v}_\beta}{\mathcal{L}\bm{e}_\alpha}.
    \end{aligned}
\label{eq:weak-form-L}
\end{equation}

The weak formulation can be discretized using the basis functions as in \eqref{eq:BasisFunction} where e.g. $\bm{e}_{\alpha,1}$ denotes $\bm{e}_\alpha$ on $\Omega_1$, to include the decomposed domain and interface. Using the basis functions, the formulations for each subdomain can be written into a finite dimensional form. For $\Omega_1$ this becomes
\begin{equation}
    \begin{aligned}
        \begin{bmatrix}
            \mathbf{M}_{\alpha,1} & 0 \\
            0 & \mathbf{M}_{\beta,1}
        \end{bmatrix}
        \odv{}{t}
        \begin{pmatrix}
            \mathbf{e}_{\alpha,1} \\
            \mathbf{e}_{\beta,1}
        \end{pmatrix} &= 
        \begin{bmatrix}
            0 & -\mathbf{D}_{\mathcal{L}^*} \\
            \mathbf{D}_{\mathcal{L}^*}^\top & 0 
        \end{bmatrix}
        \begin{pmatrix}
            \mathbf{e}_{\alpha,1} \\
            \mathbf{e}_{\beta,1}
        \end{pmatrix} + 
        \begin{bmatrix}
            0 & 0 \\
            \mathbf{B}_{\beta}^{\Gamma_1} & \mathbf{B}_{\beta}^{\Gamma_{\rm int}}
        \end{bmatrix}
        \begin{pmatrix}
            \mathbf{u}_{\partial,1} \\
            \mathbf{u}_{\partial,1}^{\Gamma_{\rm int}}
        \end{pmatrix}, \\
        \begin{pmatrix}
            \mathbf{y}_{\partial,1} \\
            \mathbf{y}_{\partial,1}^{\Gamma_{\rm int}}
        \end{pmatrix} &= 
        \begin{bmatrix}
            0 & \mathbf{T}_{\beta}^{\Gamma_1} \\
            0 & \mathbf{T}_{\beta}^{\Gamma_{\rm int}}
        \end{bmatrix}
        \begin{pmatrix}
            \mathbf{e}_{\alpha, 1} \\
            \mathbf{e}_{\beta, 1}
        \end{pmatrix}.
    \end{aligned}
\label{eq:pH-system-L*}
\end{equation}
where the output variables are computed strongly considering discrete trace operators. In an analogous manner for $\Omega_2$ it is obtained
\begin{equation}
    \begin{aligned}
        \begin{bmatrix}
            \mathbf{M}_{\alpha, 2} & 0 \\
            0 & \mathbf{M}_{\beta,2}
        \end{bmatrix}
        \odv{}{t}
        \begin{pmatrix}
            \mathbf{e}_{\alpha,2} \\
            \mathbf{e}_{\beta,2}
        \end{pmatrix} &= 
        \begin{bmatrix}
            0 & -\mathbf{D}_{\mathcal{L}}^\top \\
            \mathbf{D}_{\mathcal{L}} & 0 
        \end{bmatrix}
        \begin{pmatrix}
            \mathbf{e}_{\alpha,2} \\
            \mathbf{e}_{\beta,2}
        \end{pmatrix} + 
        \begin{bmatrix}
            \mathbf{B}_{\alpha}^{\Gamma_2} & \mathbf{B}_{\alpha}^{\Gamma_{\rm int}} \\
            0 & 0
        \end{bmatrix}
        \begin{pmatrix}
            \mathbf{u}_{\partial,2} \\
            \mathbf{u}_{\partial,2}^{\Gamma_{\rm int}}
        \end{pmatrix}, \\
        \begin{pmatrix}
            \mathbf{y}_{\partial,2} \\
            \mathbf{y}_{\partial,2}^{\Gamma_{\rm int}}
        \end{pmatrix} &= 
        \begin{bmatrix}
            \mathbf{T}_{\alpha}^{\Gamma_2} & 0 \\
            \mathbf{T}_{\alpha}^{\Gamma_{\rm int}} & 0
        \end{bmatrix}
        \begin{pmatrix}
            \mathbf{e}_{\alpha,2} \\
            \mathbf{e}_{\beta,2}
        \end{pmatrix}.
    \end{aligned}
\label{eq:pH-system-L}
\end{equation}

\subsection{Choice of the boundary functions} 
\revone{The choice of the boundary spaces follows the same rationale as in Sec.~\ref{sec:Discretization}.}
The boundary shape functions are not chosen in an independent way with respect to the state variables. Given Eqs. \eqref{eq:input_ouput_1}, \eqref{eq:input_ouput_2} and \eqref{eq:coupling}, it is natural to choose the basis functions for the inputs as being the basis function of the associated co-energy variable on the boundary subpartions. \revone{This means leveraging Eq. \eqref{eq:boundary_basis_boundary} also on the interface $\Gamma_{\rm int}$}
\begin{equation}\label{eq:boundary_basis_interface} 
\begin{aligned}
    \mathrm{span}\{\chi_{\partial, 1}\}|_{\partial\Omega_1} &= \mathrm{span} \{\mathrm{\mathcal{T}_\alpha|_{\partial\Omega_1}\chi_{\alpha, 2}}\}, \\
    \mathrm{span}\{\chi_{\partial, 2}\}|_{\partial\Omega_2} &= \mathrm{span} \{\mathcal{T}_\beta|_{\partial\Omega_2} \chi_{\beta, 1}\}.
\end{aligned}
\end{equation}

This choice will couple the two systems and is important for the domain decomposition strategy. The relations in Eq. \eqref{eq:boundary_basis_interface} provide the interconnection of the two system on $\Gamma_{\rm int}$
\begin{equation}\label{eq:interconnection_dofs}
    \begin{split}
        \mathbf{u}_{\partial,1}^{\Gamma_{\rm int}} &= \pm\mathbf{y}_{\partial,2}^{\Gamma_{\rm int}} = \pm\mathbf{T}_\alpha^{\Gamma_{\rm int}} \mathbf{e}_{\alpha, 2}, \\
        \mathbf{u}_{\partial,2}^{\Gamma_{\rm  int}} &= \mp\mathbf{y}_{\partial,1}^{\Gamma_{\rm int}} = \pm\mathbf{T}_\beta^{\Gamma_{\rm int}}
        \mathbf{e}_{\beta, 1}. \\
    \end{split}
\end{equation}
\begin{figure}
    \centering
    \includegraphics[width=0.5\linewidth]{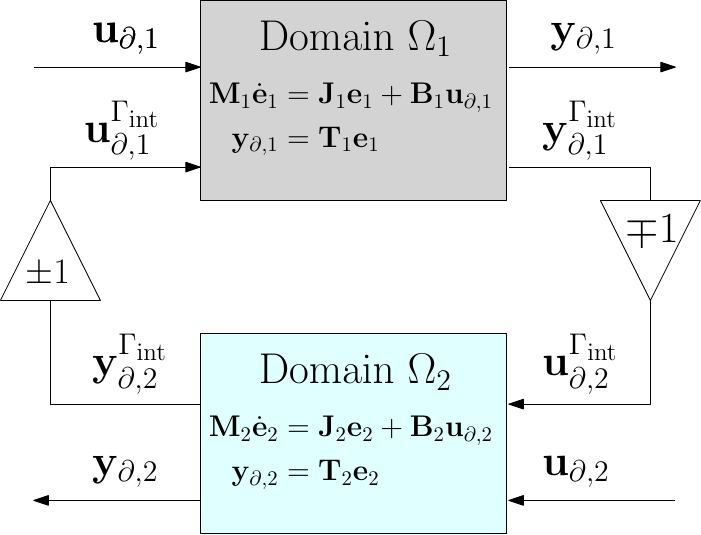}
    \caption{Feedback interconnection of the two systems arising from the domain decomposition.}
    \label{fig:gyrator_interconnection}
\end{figure}

These equations represent a feedback interconnection (cf. Fig. \ref{fig:gyrator_interconnection}) which in port-Hamiltonian systems jargon is also called a gyrator interconnection. Relations \eqref{eq:boundary_basis_interface} are also responsible for a factorization of the $\mathbf{B}$ matrices
\begin{equation}
    \mathbf{B}_\alpha^{\Gamma_{\rm int}} = (\mathbf{T}_\alpha^{\Gamma_{\rm int}})^\top \mathbf{\Psi}^{\Gamma_{\rm int}}, \qquad \mathbf{B}_\beta^{\Gamma_{\rm int}} = (\mathbf{T}_\beta^{\Gamma_{\rm int}})^\top (\mathbf{\Psi}^{\Gamma_{\rm int}})^\top,
\end{equation}
where $[\mathbf{\Psi}^{\Gamma_{\rm int}}]_{lk} = \boundary[\Gamma_{\rm int}]{\chi_{\partial, 1}^l}{\chi_{\partial, 2}^k}$. The systems found for $\Omega_1$ and $\Omega_2$ can be combined into a monolithic interconnected system for the entire domain $\Omega$. The pH-system for the full domain is provided as
\begin{equation}
\begin{aligned}
    \mathrm{Diag}
        \begin{bmatrix}
            \mathbf{M}_{\alpha, 1} \\
            \mathbf{M}_{\beta, 1}\\
            \mathbf{M}_{\alpha, 2}\\
            \mathbf{M}_{\beta, 2}
        \end{bmatrix} \odv{}{t}
        \begin{pmatrix}
            \mathbf{e}_{\alpha, 1} \\
            \mathbf{e}_{\beta, 1} \\
            \mathbf{e}_{\alpha, 2} \\
            \mathbf{e}_{\beta, 2}
        \end{pmatrix} &= 
        \begin{bmatrix}
            0 & -\mathbf{D}_{\mathcal{L}^*} & 0 & 0 \\
            \mathbf{D}_{\mathcal{L}^*}^\top & 0 & \pm \mathbf{L}^{\Gamma_{\rm int}} & 0 \\
            0 & \mp (\mathbf{L}^{\Gamma_{\rm int}})^\top & 0 & -\mathbf{D}_\mathcal{L}^\top \\
            0 & 0 & \mathbf{D}_\mathcal{L} & 0
        \end{bmatrix}
        \begin{pmatrix}
            \mathbf{e}_{\alpha,1} \\
            \mathbf{e}_{\beta,1} \\
            \mathbf{e}_{\alpha,2} \\
            \mathbf{e}_{\beta,2}
        \end{pmatrix} + 
        \begin{bmatrix}
            0 & 0 \\
            \mathbf{B}_\beta^{\Gamma_1} & 0\\
            0 & \mathbf{B}_\alpha^{\Gamma_2} \\
            0 & 0 
        \end{bmatrix}
        \begin{pmatrix}
            \mathbf{u}_{\partial,1} \\
            \mathbf{u}_{\partial,2} \\
        \end{pmatrix}, \\
        \begin{pmatrix}
            \mathbf{y}_{\partial,1} \\
            \mathbf{y}_{\partial,2} \\
        \end{pmatrix} &= 
        \begin{bmatrix}
            0 & \mathbf{T}_\beta^{\Gamma_1} & 0 & 0 \\
            0 & 0 & \mathbf{T}_\alpha^{\Gamma_2} & 0\\
        \end{bmatrix}
        \begin{pmatrix}
            \mathbf{e}_{\alpha,1} \\
            \mathbf{e}_{\beta,1} \\
            \mathbf{e}_{\alpha,2} \\
            \mathbf{e}_{\beta,2}
        \end{pmatrix},
    \end{aligned}
\label{eq:Monolithic-System}
\end{equation}
where $\mathbf{L}^{\Gamma_{\rm int}} =  (\mathbf{\Psi}^{\Gamma_{\rm int}}\mathbf{T}_\beta^{\Gamma_{\rm int}})^\top \mathbf{T}_\alpha^{\Gamma_{\rm int}}$. The structure of the system is again Hamiltonian and can be written compactly as 
\begin{equation}
\begin{aligned}
    \mathbf{M}\dot{\mathbf{e}} &= \mathbf{J} \mathbf{e} + \mathbf{Bu}, \\
    \mathbf{y} &= \mathbf{T} \mathbf{e}, 
\end{aligned}
\end{equation}
where $\mathbf{J} = - \mathbf{J}^\top$ is skew-symmetric.
\revone{
\begin{remark}
    The domain decomposition strategy does not require matching interfaces. The feedback interconnection \eqref{eq:interconnection_dofs} corresponds to a choice of numerical fluxes (as in Discontinuous Galerkin methods).
\end{remark}
}

\revtwo{
\subsection{Extension to the nonlinear case}\label{sec:nonlinear}
The methodology can be applied to nonlinear systems but only if the nonlinearity has to enter the system in an algebraic way. In other words the method is applicable to semi-linear systems of the form
\begin{equation}
    \mathcal{M}\partial_t \bm{e} = \mathcal{J}_{d}\bm{e} + \mathcal{J}_a(\bm{e})\bm{e},
\end{equation}
where $\mathcal{J}_d$ is a differential operator of the form
\begin{equation}
    \mathcal{J}_d = 
    \begin{bmatrix}
        0 & -\mathcal{L}^* \\
        \mathcal{L} & 0
    \end{bmatrix},
\end{equation}
and $\mathcal{J}_a(\bm{e})$ is a nonlinear algebraic operator the depends on the coenergy $\bm{e}$. An example of such a semilinear dynamical systems is the intrinsic geometrically exact beam introduced in \cite{hodges2003intrinsic}. In Sec. \ref{sec:intrinsic_nonlinear_beam} a numerical test for this model will be detailed. \\
If the differential operator appears in a nonlinear term, it might not be possible to apply a primal dual formulation. 
We illustrate this issue on an example.
\paragraph{A nonlinear example that does not fit in the domain decomposition strategy: geometrically nonlinear elasticity}
In geometrically nonlinear elasticity the infinitesimal strain tensor is replaced by the Green-Lagrange tensor 
$$
\bm{E} := \frac{1}{2}(\bm{F}^\top\bm{F} - \bm{I}), \qquad \bm{F} :=  \bm{I} + \nabla \bm{q},
$$
where $\bm{q}$ is the displacement $[\nabla \bm{q}]_{ij} = \partial_j q_i$ is the gradient of a vector defined row-wise, and $\bm{F}$ is the deformation gradient. The kinetic and potential energies are given by 
$$
\begin{aligned}
T = \frac{1}{2} \int_\Omega \rho ||{\partial_t \bm{q}}||^2 \d\Omega, \qquad
V = \frac{1}{2} \int_\Omega \bm{E}:\bm{K}\bm{E} \; \d\Omega,
\end{aligned}
$$
where $\bm{K}$ is the stiffness tensor. For the potential energy a Saint-Venant Kirchhoff material model has been used. The Euler-Lagrange equations are then given by
\begin{equation*}
\rho\, \partial_{tt}\bm{q} = \Div(\bm{F}\bm{S}),
\end{equation*}
where $\Div$ is the row-wise divergence of a tensor and $\bm{S}=\bm{K}\bm{E}$ is the second Piola-Kirchhoff stress tensor.  By introducing the dynamical equation for the second Piola-Kirchhoff stress tensor, the Hamiltonian structure of the equations can be highlighted~\cite{thoma2024velocity}:
\begin{equation*}
\begin{aligned}
\partial_t \bm{q} &= \bm{v}, \\
\begin{bmatrix}
    \rho & 0 \\
    0 & \bm{C}
\end{bmatrix}
\pdv{}{t}
\begin{pmatrix}
    \bm{v} \\
    \bm{S}
\end{pmatrix} &= 
\begin{bmatrix}
    0 & \Div (\bm{F} \; \circ) \\
    \sym( \bm{F}^\top \nabla \; \circ) & 0 \\
\end{bmatrix}
\begin{pmatrix}
    \bm{v} \\
    \bm{S}
\end{pmatrix},    
\end{aligned}
\end{equation*}
where $\bm{C}:=\bm{K}^{-1}$ is the compliance tensor. In this case the  differential operator $\mathcal{L}$ and its adjoint $\mathcal{L}^*$ contains the deformation gradient (that is seen as a parameter for defining the adjoint) as
$$
\begin{aligned}
    \mathcal{L}(\nabla \bm{q}) = \sym( \bm{F}^\top \nabla \; \circ), \qquad
\mathcal{L}^*(\nabla \bm{q}) & =  \Div (\bm{F} \; \circ).
\end{aligned}
$$
Because of the fact the these terms are nonlinear, the discretization can only be performed by integrating by part the $\mathcal{L}^*$ operator. The resulting weak formulation reads
\begin{equation*}
\begin{aligned}
    \partial_t \bm{q}_h &= \bm{v}_h, \\
    \innerproduct[\Omega]{\bm{\psi}}{\rho\,\partial_t \bm{v}_h} &= - \innerproduct[\Omega]{\bm{F}_h^\top\nabla\bm{\psi}}{\bm{S}_h},  \\
    \innerproduct[\Omega]{\bm{\Psi}}{\bm{C}\,\partial_t \bm{S}_h} &= +\innerproduct[\Omega]{\bm{\Psi}}{\bm{F}_h^\top\nabla\bm{v}_h},  \\
\end{aligned} \qquad 
\begin{aligned}
\\
\text{forall } \bm{\psi} \in V_h, \\
\text{forall } \bm{\Psi} \in \Sigma_h.
\end{aligned}
\end{equation*}
For this example it is not clear how a dual system with opposite treatment of the boundary conditions.
}

\revone{
\section{Time integration}\label{sec:time-integration}
We present two different integrators for the system \eqref{eq:Monolithic-System}, the St\"ormer Verlet scheme and the implicit midpoint. The first one allows for decoupling of the two domains, but it is not a Poisson map. The second imposes a monolithic resolution of the problem, but guarantees the preservation of the Poisson structure. To illustrate this method, we assume that the boundary data are homogeneous ($\mathbf{u}_{\partial, 1}=0, \; \mathbf{u}_{\partial, 2}=0$). System \eqref{eq:Monolithic-System} takes the partitioned form
\begin{equation}\label{eq:interconnected_system}
\begin{aligned}
        \begin{bmatrix}
            \mathbf{M}_{1} & 0 \\
            0 & \mathbf{M}_{2}
        \end{bmatrix} 
        \odv{}{t}
        \begin{pmatrix}
            \mathbf{e}_{1} \\
            \mathbf{e}_{2}
        \end{pmatrix} &=
        \begin{bmatrix}
            \mathbf{J}_1 & +\mathbf{G} \\
            -\mathbf{G}^\top & \mathbf{J}_2
        \end{bmatrix}
        \begin{pmatrix}
            \mathbf{e}_{1} \\
            \mathbf{e}_{2}
        \end{pmatrix}.
\end{aligned}
\end{equation}
To simplify the analysis the system can be rewritten by a change of variable $\widehat{\mathbf{e}}_1 = \mathbf{C}_1 \mathbf{e}_1, \; \widehat{\mathbf{e}}_2 = \mathbf{C}_2 \mathbf{e}_2$ where $\mathbf{C}_1, \; \mathbf{C}_2$ are the Cholesky factors of the mass matrices $\mathbf{M}_1 = \mathbf{C}_1^\top \mathbf{C}_1, \; \mathbf{M}_2 = \mathbf{C}_2^\top \mathbf{C}_2$, leading to 
\begin{equation}\label{eq:interconnected_system_hat}
\begin{aligned}
        \odv{}{t}
        \begin{pmatrix}
            \widehat{\mathbf{e}}_{1} \\
            \widehat{\mathbf{e}}_{2}
        \end{pmatrix} &=
        \begin{bmatrix}
            \widehat{\mathbf{J}}_1 & +\widehat{\mathbf{G}} \\
            -\widehat{\mathbf{G}}^\top & \widehat{\mathbf{J}}_2
        \end{bmatrix}
        \begin{pmatrix}
            \widehat{\mathbf{e}}_{1} \\
            \widehat{\mathbf{e}}_{2}
        \end{pmatrix}, \qquad \text{or compactly } \quad \dot{\widehat{\mathbf{e}}} = \widehat{\mathbf{J}}\widehat{\mathbf{e}},
\end{aligned}
\end{equation}
where $\widehat{\mathbf{J}}_1 = \mathbf{C}_1^{-\top} {\mathbf{J}}_1 \mathbf{C}_1^{-1}, \; \widehat{\mathbf{J}}_2 = \mathbf{C}_2^{-\top} {\mathbf{J}}_2 \mathbf{C}_2^{-1}$ and $\widehat{\mathbf{G}} = \mathbf{C}_1^{-\top} {\mathbf{G}} \mathbf{C}_2^{-1}.$
\paragraph{Implicit midpoint scheme}
Consider system \eqref{eq:interconnected_system_hat} $\dot{\mathbf{e}} = {\mathbf{J}}{\mathbf{e}}$, where the hat $\widehat{\cdot}$ is omitted for simplicity. The implicit midpoint rule gives
    $$
     \frac{{\mathbf{e}}^{n+1}-{\mathbf{e}}^n}{\Delta t} = {\mathbf{J}} \left(\frac{{\mathbf{e}}^n + {\mathbf{e}}^{n+1}}{2}\right).
    $$
We now recall a known result that is not easy to find in the literature.
\begin{proposition}
The implicit midpoint scheme applied to a linear Poisson system is a Poisson map.
\end{proposition}
\begin{proof}
    Using a time rescaling, we set $\Delta t/2 = 1$. 
    The application of the midpoint rule leads to the recursion 
    \begin{equation}\label{eq:discrete_flow}
        {\mathbf{e}}^{n+1} = \mathrm{Cay}({\mathbf{J}}) {\mathbf{e}}^{n}, \qquad  \mathrm{Cay}({\mathbf{J}}) := (\mathbf{I} -{\mathbf{J}})^{-1}(\mathbf{I} + {\mathbf{J}}).
    \end{equation}
    For the discrete flow \eqref{eq:discrete_flow} to be a Poisson map, it must hold
    $$
    \mathrm{Cay}({\mathbf{J}}) \; \mathbf{J} \; \mathrm{Cay}({\mathbf{J}})^\top = \mathbf{J}.
    $$
    By exploiting the property $\mathbf{J} = -\mathbf{J}^\top$, the term $\mathrm{Cay}({\mathbf{J}})^\top$ gives
    $$
    \mathrm{Cay}({\mathbf{J}})^\top = (\mathbf{I} -{\mathbf{J}})(\mathbf{I} + {\mathbf{J}})^{-1}.
    $$
    So the discrete flow can be rewritten as
    $$
    \mathrm{Cay}({\mathbf{J}}) \; \mathbf{J} \; \mathrm{Cay}({\mathbf{J}})^\top =  (\mathbf{I} -{\mathbf{J}})^{-1}(\mathbf{I} + {\mathbf{J}}) \mathbf{J} (\mathbf{I} -{\mathbf{J}})(\mathbf{I} + {\mathbf{J}})^{-1}.
    $$
    The following commuting properties holds
    \begin{equation}\label{eq:commutation_matrices}
        \begin{aligned}
            \mathbf{J}(\mathbf{I} + \mathbf{J}) = (\mathbf{I} + \mathbf{J}) \mathbf{J}, \qquad
            \mathbf{J}(\mathbf{I} - \mathbf{J}) = (\mathbf{I} - \mathbf{J}) \mathbf{J}, \qquad
            (\mathbf{I} + \mathbf{J})(\mathbf{I} - \mathbf{J}) = (\mathbf{I} - \mathbf{J})(\mathbf{I} + \mathbf{J}).            
        \end{aligned}
    \end{equation}
    Using these relations, it is obtained
    $$
    (\mathbf{I} -{\mathbf{J}})^{-1}(\mathbf{I} + {\mathbf{J}}) \mathbf{J} (\mathbf{I} -{\mathbf{J}})(\mathbf{I} + {\mathbf{J}})^{-1} = \mathbf{J}.
    $$
\end{proof}
So the implicit midpoint apply to the system given a Poisson map, leading to a symplectic integrator
\paragraph{St\"ormer-Verlet scheme}
The St\"ormer Verlet scheme is a partitioned Runge-Kutta scheme \cite{hairer2006geometric}. In the present case it takes the form
\begin{equation*}
    \begin{aligned}
    \frac{\widehat{\mathbf{e}}_1^{n+1}-\widehat{\mathbf{e}}_1^n}{\Delta t} &= \widehat{\mathbf{J}}_1 \left(\frac{\widehat{\mathbf{e}}_1^n + \widehat{\mathbf{e}}_1^{n+1}}{2}\right) + \widehat{\mathbf{G}}\widehat{\mathbf{e}}_2^{n+\frac{1}{2}}, \\
    \frac{\widehat{\mathbf{e}}_2^{n +\frac{1}{2}}-\widehat{\mathbf{e}}_2^{n -\frac{1}{2}}}{\Delta t} &= \widehat{\mathbf{J}}_2 \left(\frac{\widehat{\mathbf{e}}_2^{n +\frac{1}{2}}+\widehat{\mathbf{e}}_2^{n-\frac{1}{2}}}{2}\right) - \widehat{\mathbf{G}}^\top \widehat{\mathbf{e}}_1^n. \\
    \end{aligned}
\end{equation*}
To start the iterations the St\"ormer-Verlet initial value $\widehat{\mathbf{e}}_2^{\frac{1}{2}}$ is obtained using
$$\widehat{\mathbf{e}}_2^{\frac{1}{2}} =  (\mathbf{I}- \frac{\Delta t}{2}\widehat{\mathbf{J}}_2)^{-1}\widehat{\mathbf{e}}_2^{0} - \frac{\Delta t}{2}\mathbf{G} ^\top\widehat{\mathbf{e}}_1^0.
$$
\begin{remark}
    The St\"ormer-Verlet integrator is not a Poisson map for system \eqref{eq:interconnected_system_hat}. By Lemma 4.9 in \cite{hairer2006geometric}, the St\"ormer-Verlet integrator is not symplectic for the canonical Hamiltonian system obtained via the Darboux-Lie theorem.
\end{remark}
}

%% file: Text/numerical_examples.tex
\section{Numerical examples}\label{sec:examples}
The domain decomposition strategy is applied to four different examples:
\begin{itemize}
    \item \revtwo{the one dimensional nonlinear geometrically exact intrinsic beam model;} 
    \item the two dimensional wave equation;
    \item \revthree{the two dimensional linear elastodynamics problem};
    \item \revthree{the Mindlin plate problem};
\end{itemize}
The decomposition of the mesh has been implemented using \textsc{GMSH} \cite{gmsh2009}. All the investigations will be performed employing the finite element library \textsc{Firedrake} \cite{Rathgeber2016Firedrake:Abstractions}.

%% file: Text/pH_intrinsic_nonlinear_beam.tex
\revtwo{
\subsection{A 1D non linear example: geometrically exact intrinsic geometrically exact beams}\label{sec:intrinsic_nonlinear_beam}
The domain decomposition strategy applied to a semilinear example as described in Sec. \ref{sec:nonlinear}. An example of a semilinear problem is the intrinsic formulation of geometrically exact beams first proposed by Dewey Hodges \cite{hodges2003intrinsic}. This model describes the motion of the beam cross section as a rigid motion and captures the geometric nonlinearity in the deformation without making any additional simplification. The model accounts for shear deformability. The description is in the material reference frame as all variables follow the motion of the cross section. 
In this example a one dimensional beam with length $L$ under a Dirichlet condition at $x=0$ (velocities are set to zero here) and a Neumann boundary condition at $x=L$ (forces and torques are applied at this node), is decomposed into two subdomains $\Omega_1$ and $\Omega_2$ using an interface vertex $\Gamma_{\rm int}$. The results shown in this example use an interface vertex located at $x_{\rm int}=L/2$, but it should be restated that its position is arbitrary. The domain decompostion is plot in Fig. \ref{fig:nonlinear_beam}.}
\begin{figure}[h]
\centering
    \begin{tikzpicture}[scale=3]
	\node[black] at (1,0.2) {$\Gamma_2$};
	\node[black] at (3,0.2) {$\Gamma_1$};
	\node[black] at (2,0.2) {$\Gamma_{\rm int}$};
        \node[black] at (2,-0.2) {$L$};

        \draw[<->, thick] (1,-0.1) -- (3,-0.1);
	\draw[black, thick, fill=gray!40] (1,0) rectangle (2,0.1);
	\draw[black, thick, fill=cyan!40] (2,0) rectangle (3,0.1);
	
    \end{tikzpicture}
\caption{The decomposed beam with mixed boundary conditions.}
\label{fig:nonlinear_beam}
\end{figure}
\revtwo{
The Hamiltonian is given by 
$$
H = \frac{1}{2}\rho A ||\bm{v}||^2 +  \frac{1}{2}\rho \bm{w}^\top\bm{J}\bm{w} +  \frac{1}{2}\bm{n}^\top \bm{C}_t \bm{n} +  \frac{1}{2}\bm{m}^\top \bm{C}_r \bm{m},
$$
where $\bm{v}, \; \bm{w} \in \mathbb{R}^3$ are the {material} linear and the angular velocity, respectively, $\bm{n}, \; \bm{m} \in \mathbb{R}^3$ are the material force and bending moment resultants, respectively. The parameters are the density $\rho$, the cross section area $A$, the moment of area matrix $\bm{J} \in \mathbb{R}^{3\times 3}$, and the translational and rotational compliance $\bm{C}_t, \bm{C}_r \in \mathbb{R}^{3\times 3}$. 
The co-energy variables are given by
$$
\begin{aligned}
\bm{\pi}_v &= \partial_{\bm{v}} H = \rho A \bm{v}, \\
\bm{\pi}_w &= \partial_{\bm{w}} H = \rho \bm{J} \bm{w}, \\
\end{aligned}\qquad 
\begin{aligned}
\bm{\gamma} &= \partial_{\bm{n}} H = \bm{C}_t \bm{n}, \\
\bm{\kappa} &= \partial_{\bm{m}} H = \bm{C}_r \bm{m}, 
\end{aligned}
$$
denoting material momentum and strain quantities. 
In the following the notation $[\bm{v}]_\times$ denotes the skew-symmetric matrix obtained as
\begin{equation}
\bm{v} =
\begin{pmatrix}
v_x \\
v_y \\
v_z
\end{pmatrix}
\rightarrow
[\bm{v}]_\times :=
\begin{bmatrix}
0 & -v_z & v_y \\
v_z & 0 & -v_x \\
-v_y & v_x & 0
\end{bmatrix},
\end{equation}
to rewrite the cross-product as matrix-vector multiplication, i.e. $\bm{v} \times \bm{u} = [\bm{v}]_\times \bm{u}$ for arbitrary $\bm{u} \in \mathbb{R}^3.$ Denoting with variable ${s} \in [0, L]$ the material arc length coordinate, the dynamics of the system over an interval $\Omega=[0, L]$ is given by
\begin{equation}\label{eq:intrinsic_beam}
\mathrm{Diag}
\begin{bmatrix}
\rho A \\
\rho \bm{J} \\
\bm{C}_t\\
\bm{C}_r
\end{bmatrix}
\partial_t
\begin{pmatrix}
\bm{v} \\
\bm{w} \\
\bm{n} \\
\bm{m}
\end{pmatrix}
=
\left(
\begin{bmatrix}
0 &  0 & \partial_s & 0 \\
0 & 0 & 0  & \partial_s \\
\partial_s & 0 & 0 & 0 \\
0 & \partial_s  & 0 & 0
\end{bmatrix}
+ 
\begin{bmatrix}
0 &  [\bm{\pi}_V]_{\times} & [\bm{\kappa}]_{\times} & 0 \\
[\bm{\pi}_V]_{\times} &  [\bm{\pi}_W]_{\times} & [\bm{\gamma} + \bm{e}_1]_{\times}  & [\bm{\kappa}]_{\times} \\
[\bm{\kappa}]_{\times} & [\bm{\gamma} + \bm{e}_1]_{\times}  & 0 & 0 \\
0 & [\bm{\kappa}]_{\times} & 0 & 0
\end{bmatrix}
\right)
\begin{pmatrix}
\bm{v} \\
\bm{w} \\
\bm{n} \\
\bm{m}
\end{pmatrix},
\end{equation}
where $\bm{e}_1 = [1 \; 0\; 0]^\top$. For this examples of operators $\mathcal{L}, \; \mathcal{L}^*$ and the variables $\bm{e}_{\alpha}, \; \bm{e}_\beta$ are given by
$$
\mathcal{L} = \begin{bmatrix}
\partial_s & 0 \\
0 & \partial_s 
\end{bmatrix}, \qquad  \mathcal{L}^* = -\begin{bmatrix}
\partial_s & 0 \\
0 & \partial_s 
\end{bmatrix}, \qquad
\bm{e}_{\alpha} = \begin{pmatrix}
    \bm{v} \\ \bm{w}
\end{pmatrix} \qquad 
\bm{e}_{\beta} = \begin{pmatrix}
    \bm{n} \\ \bm{m}
\end{pmatrix}.
$$
Notice that the linearization of System \eqref{eq:intrinsic_beam} gives the port-Hamiltonian formulation of the Timoshenko beam
\begin{equation}\label{eq:timoshenkp_beam}
\mathrm{Diag}
\begin{bmatrix}
\rho A \\
\rho \bm{J} \\
\bm{C}_t\\
\bm{C}_r
\end{bmatrix}
\partial_t
\begin{pmatrix}
\bm{v} \\
\bm{w} \\
\bm{n} \\
\bm{m}
\end{pmatrix}
=
\begin{bmatrix}
0 &  0 & \partial_s & 0 \\
0 & 0 & [\bm{e}_1]_{\times}  & \partial_s \\
\partial_s & [\bm{e}_1]_{\times} & 0 & 0 \\
0 & \partial_s  & 0 & 0
\end{bmatrix}
\begin{pmatrix}
\bm{v} \\
\bm{w} \\
\bm{n} \\
\bm{m}
\end{pmatrix}.
\end{equation}
The discretization is explained on the linear part only as the nonlinearity is simply projected on finite element spaces. If the last two lines are integrated by parts then one obtains the weak formulation: find $\bm{v}, \bm{w} \in L^2(\Omega_1; \mathbb{R}^3), \; \bm{n}, \bm{m} \in H^1(\Omega_1; \mathbb{R}^3)$ such that forall $\bm{\psi}_v, \bm{\psi}_w \in L^2(\Omega_1; \mathbb{R}^3), \; \bm{\psi}_n, \bm{\psi}_m \in H^1(\Omega_1; \mathbb{R}^3)$
\begin{equation}
\begin{aligned}
   \innerproduct[\Omega_1]{\bm{\psi}_v}{\rho A \partial_t \bm{v}} &= \innerproduct[\Omega_1]{\bm{\psi}_v}{\partial_s \bm{n}}, \\
    \innerproduct[\Omega_1]{\bm{\psi}_w}{\rho \bm{J} \partial_t \bm{w}} &= \innerproduct[\Omega_1]{\bm{\psi}_w}{[\bm{e}_1]_{\times} \bm{n}} + \innerproduct[\Omega_1]{\bm{\psi}_w}{\partial_s \bm{m}}, \\
    \innerproduct[\Omega_1]{\bm{\psi}_n}{\bm{C}_t \partial_t \bm{n}} &= \innerproduct[\Omega_1]{\bm{\psi}_n}{[\bm{e}_1]_{\times} \bm{w}} - \innerproduct[\Omega_1]{\partial_s \bm{\psi}_n}{\bm{v}} + \boundary[\partial \Omega_1]{\bm{\psi}_n }{\bm{v}} + \boundary[\Gamma_{\rm int}]{\bm{\psi}_n }{\bm{v}}, \\
    \innerproduct[\Omega_1]{\bm{\psi}_m}{\bm{C}_r \partial_t \bm{m}} &= - \innerproduct[\Omega_1]{\partial_s \bm{\psi}_m}{ \bm{w}} + \boundary[\partial \Omega_1]{\bm{\psi}_m }{\bm{w}} + \boundary[\Gamma_{\rm int}]{\bm{\psi}_m }{\bm{w}}, \\  
\end{aligned}
\end{equation}
where the trace operator $\mathcal{T}_\beta : H^1([a, b]; \, \mathbb{R}^6) \rightarrow \mathbb{R}^{12}$ is given by
$$
\mathcal{T}_\beta\begin{pmatrix}
    \bm{f} \\ \bm{g}
\end{pmatrix} = \begin{pmatrix}
    +\bm{f}(b) \\ -\bm{f}(a) \\ +\bm{g}(b) \\ -\bm{g}(a)
\end{pmatrix}.
$$
If the first two lines are integrated by parts then one obtains the following formulation: find $\bm{v}, \bm{w} \in H^1(\Omega_2; \mathbb{R}^3), \; \bm{n}, \bm{m} \in L^2(\Omega_2; \mathbb{R}^3)$ such that forall $\bm{\psi}_v, \bm{\psi}_w \in H^1(\Omega_2; \mathbb{R}^3), \; \bm{\psi}_n, \bm{\psi}_m \in L^2(\Omega_2; \mathbb{R}^3)$
\begin{equation}
\begin{aligned}
    \innerproduct[\Omega_2]{\bm{\psi}_v}{\rho A \partial_t \bm{v}} &= - \innerproduct[\Omega_2]{\partial_s  \bm{\psi}_v}{\bm{n}} + \boundary[\partial \Omega_2]{\bm{\psi}_v}{\bm{n}} + \boundary[\Gamma_{\rm int}]{\bm{\psi}_v}{\bm{n}}, \\
    \innerproduct[\Omega_2]{\bm{\psi}_w}{\rho \bm{J} \partial_t \bm{w}} &= \innerproduct[\Omega_2]{\bm{\psi}_w}{[\bm{e}_1]_{\times} \bm{n}} - \innerproduct[\Omega_2]{\partial_s \bm{\psi}_w}{\bm{m}} + \boundary[\partial \Omega_2]{\bm{\psi}_w}{\bm{m}} + \boundary[\Gamma_{\rm int}]{\bm{\psi}_w}{\bm{m}}, \\
    \innerproduct[\Omega_2]{\bm{\psi}_n}{\bm{C}_t \partial_t \bm{n}} &= \innerproduct[\Omega_2]{\bm{\psi}_n}{[\bm{e}_1]_{\times} \bm{w}} + \innerproduct[\Omega_2]{\bm{\psi}_n}{\partial_s \bm{v}}, \\
    \innerproduct[\Omega_2]{\bm{\psi}_m}{\bm{C}_r \partial_t \bm{m}} &= \innerproduct[\Omega_2]{\bm{\psi}_m}{\partial_s \bm{w}}, \\    
\end{aligned}
\end{equation}
and the trace operator $\mathcal{T}_\alpha : H^1([a, b]; \, \mathbb{R}^6) \rightarrow \mathbb{R}^{12}$ is given by
$$
\mathcal{T}_\alpha\begin{pmatrix}
    \bm{f} \\ \bm{g}
\end{pmatrix} = \begin{pmatrix}
    \bm{f}(b) \\ \bm{f}(a) \\ \bm{g}(b) \\ \bm{g}(a)
\end{pmatrix}.
$$
\paragraph{Finite element basis}
The finite element family used to solve this problem is the Discontinuous Galerkin elements of order 1 (DG$_0$) to discretize the $L^2$ space, and linear Lagrange finite elements (CG$_1$) to discretize the $H^1$ space, though the mixed finite element spaces are different on each subdomain. This choice is justified by the de-Rham complex. This example constitutes the simplest example of discrete de~Rham subcomplex
}
\begin{center}
\begin{tikzcd}
H^1 \arrow[r, "\partial_{x}"] \arrow[d, "\Pi"] & L^2 \arrow[d, "\Pi"] \\
\mathrm{CG}_1 \arrow[r, "\partial_{x}"] & \mathrm{DG}_0
\end{tikzcd} 
\end{center}
\revtwo{
Therefore, for the $\mathcal{I}_h^{\Omega_1}$ mesh the finite dimensional spaces are
\begin{equation*}
    \begin{aligned}
        V_{\alpha, 1} &= \{u_h \in L^2(\Omega_1; \mathbb{R}^6)|\; \forall E \in \mathcal{I}_h^{\Omega_1}, \; u_h|_E \in \mathrm{DG}_0(\mathbb{R}^6)\}, \\
       V_{\beta, 1} &= \{u_h \in H^1(\Omega_1; \mathbb{R}^6)|\; \forall E \in \mathcal{I}_h^{\Omega_1}, \; u_h|_E \in \mathrm{CG}_1(\mathbb{R}^6)\},
    \end{aligned}
\end{equation*}
where the notation $L^2(\Omega_1; \mathbb{R}^6), \; H^1(\Omega_1; \mathbb{R}^6), \; \mathrm{DG}_0(\mathbb{R}^6), \; \mathrm{CG}_1(\mathbb{R}^6)$ indicates the functions are six dimensional vectors. For the $\mathcal{I}_h^{\Omega_2}$ domain the finite dimensional spaces are
\begin{equation*}
    \begin{aligned}
        V_{\alpha, 2} &= \{u_h \in H^1(\Omega_2; \mathbb{R}^6)|\; \forall E \in \mathcal{I}_h^{\Omega_2}, \; u_h|_E \in \mathrm{CG}_1(\mathbb{R}^6) \}, \\
        V_{\beta, 2} &= \{u_h \in L^2(\Omega_2; \mathbb{R}^6)|\; \forall E \in \mathcal{I}_h^{\Omega_2}, \; u_h|_E \in \mathrm{DG}_0(\mathbb{R}^6) \}.
    \end{aligned}
\end{equation*}
Introducing the finite element approximation, the following ODE is obtained for subdomain $\Omega_1$  
\begin{equation}\label{eq:discrete_clamped_intrinsic}
\mathrm{Diag}
\begin{bmatrix}
\mathbf{M}_{v, 1} \\ 
\mathbf{M}_{w, 1} \\
\mathbf{M}_{n, 1} \\
\mathbf{M}_{m, 1}
\end{bmatrix}
    \odv{}{t}
    \begin{pmatrix}
        \mathbf{v}_1 \\
        \mathbf{w}_1 \\
        \mathbf{n}_1 \\
        \mathbf{m}_1
    \end{pmatrix} = 
    \begin{bmatrix}
        0 & 0 & \mathbf{D}_{\partial_s} & 0 \\
        0 & 0 & [\mathbf{e}_{1}]_\times & \mathbf{D}_{\partial_s} \\
        -\mathbf{D}_{\partial_s}^\top & [\mathbf{e}_{1}]_\times & 0 & 0 \\
        0 & -\mathbf{D}_{\partial_s}^\top & 0 & 0
    \end{bmatrix}
    \begin{pmatrix}
        \mathbf{v}_1 \\
        \mathbf{w}_1 \\
        \mathbf{n}_1 \\
        \mathbf{m}_1
    \end{pmatrix} + \begin{bmatrix}
        0 & 0 \\
        0 & 0 \\
         \mathbf{T}_{\beta}^\top & 0 \\
        0 & \mathbf{T}_{\beta}^\top \\
    \end{bmatrix}\begin{pmatrix}
        \mathbf{u}_v \\
        \mathbf{u}_w
    \end{pmatrix},
\end{equation}
where $\mathbf{T}_{\beta}$ is the normal trace matrix, taking values $1$ or $-1$ for the right and left extremity degree of freedom. The discrete system for the domain $\Omega_2$ is given by
\begin{equation}\label{eq:discrete_free_intrinsic}
\mathrm{Diag}
\begin{bmatrix}
\mathbf{M}_{v, 2} \\ 
\mathbf{M}_{w, 2} \\
\mathbf{M}_{n, 2} \\
\mathbf{M}_{m, 2}
\end{bmatrix}
    \odv{}{t}
    \begin{pmatrix}
        \mathbf{v}_2 \\
        \mathbf{w}_2 \\
        \mathbf{n}_2 \\
        \mathbf{m}_2
    \end{pmatrix} = 
    \begin{bmatrix}
        0 & 0 & -\mathbf{D}_{\partial_s}^\top & 0 \\
        0 & 0 & [\mathbf{e}_{1}]_\times & -\mathbf{D}_{\partial_s}^\top \\
        \mathbf{D}_{\partial_s} & [\mathbf{e}_{1}]_\times & 0 & 0 \\
        0 & \mathbf{D}_{\partial_s} & 0 & 0
    \end{bmatrix}
    \begin{pmatrix}
        \mathbf{v}_2 \\
        \mathbf{w}_2 \\
        \mathbf{n}_2 \\
        \mathbf{m}_2
    \end{pmatrix} + \begin{bmatrix}
        \mathbf{T}_\alpha^\top & 0 \\
        0 & \mathbf{T}_\alpha^\top \\
        0 & 0 \\
        0 & 0 \\
    \end{bmatrix}\begin{pmatrix}
        \mathbf{u}_n \\
        \mathbf{u}_m
    \end{pmatrix},
\end{equation}
where $\mathbf{T}_\alpha$ is a localization matrix that picks the degrees of freedom at the boundary. 
The nonlinear terms are then discretized by simply projecting on the finite element basis, leading to two finite dimensional ODEs of the form
$$
\begin{aligned}
   \mathbf{M}_i \dot{\mathbf{e}}_i &= \mathbf{J}_{d, i}\mathbf{e}_i + \mathbf{J}_{a, i}(\mathbf{e}_i)\mathbf{e}_i + \mathbf{B}_i \mathbf{u}_i, \\
   \mathbf{y}_i &= \mathbf{B}_i^\top \mathbf{e}_i, \qquad  i=\{1,2\}, 
\end{aligned}
$$
where $\mathbf{J}_{d, i}= - \mathbf{J}_{d, i}^\top$ is the matrix associated to the discretization of the differential operator and $\mathbf{J}_{a, i}(\mathbf{e}_i)= - \mathbf{J}_{d, i}^\top(\mathbf{e}_i)$ is the matrix discretization of the nonlinear terms, modulated by the state variable. The interconnection is then performed using a feedback interconnection leading to a final system of the form
$$
\begin{bmatrix}
    \mathbf{M}_1 & 0 \\
    0 & \mathbf{M}_2
\end{bmatrix}
\odv{}{t}\begin{pmatrix}
    \mathbf{e}_1 \\
    \mathbf{e}_2 \\
\end{pmatrix} = 
\begin{bmatrix}
    \mathbf{J}_1(\mathbf{e}_1) & \mathbf{B}_1^{\Gamma_{\rm int}} (\mathbf{B}_2^{\Gamma_{\rm int}})^\top \\
    - \mathbf{B}_2^{\Gamma_{\rm int}} (\mathbf{B}_1^{\Gamma_{\rm int}})^\top & \mathbf{J}_2(\mathbf{e}_2)
\end{bmatrix}
\begin{pmatrix}
    \mathbf{e}_1 \\
    \mathbf{e}_2 \\
\end{pmatrix} + \begin{bmatrix}
    \mathbf{B}_{1}^{\Gamma_1} & 0 \\
    0 & \mathbf{B}_{2}^{\Gamma_2} \\
\end{bmatrix}
\begin{pmatrix}
    \mathbf{u}_{\partial, 1} \\
    \mathbf{u}_{\partial, 2} \\
\end{pmatrix}.
$$
\paragraph{Time domain simulation: Quasi static roll up of a cantilever beam}
The benchmark problem of rolling up a cantilever beam is used to demonstrate that the proposed formulation effectively avoids shear locking and is suitable for quasi-static simulations. An initially straight cantilever beam is clamped at \(s=0\), enforcing the Dirichlet boundary conditions
\[
\bm{v}(0,t)=0,\qquad \bm{w}(0,t)=0, \qquad \bm{n}(0,t)=0, \qquad \bm{m}(0,t)=0.
\]
A quasi-static torque is applied at the free end \(s=L\), and according to reference results, a torque of \(m_{\text{rollup}} = 2\pi EI / L\) rolls the beam into a complete circular arc. The parameters from~\cite{hante2022} are used, cf. Table \ref{tab:parameters_rollup}. 
To enable quasi-static behavior within the intrinsically dynamical framework, inertial terms are neglected by setting \(\rho = 0\). In this setting, velocity-type variables act as incremental coordinates, and the simulation time \(t \in [0,1]\) serves as a loading factor such that
\[
m_0(t) = t\, m_{\text{rollup}}.
\]
\revone{The time integration is performed using the implicit midpoint method.} The results for the free-free (corresponding to System \eqref{eq:discrete_free_intrinsic}), clamped-clamped (corresponding to System \eqref{eq:discrete_clamped_intrinsic}) and interconnected formulation are shown in Figure \ref{fig:configurations_rollup}. The free-free case requires a strong imposition of the clamp boundary condition. The clamped-clamped case required the strong imposition of the end bending moment. The interconnected system does not require any Lagrange multiplier. The proposed formulation does not exhibit locking, even under coarse spatial discretization.
}
\begin{table}[htb]
\centering
\begin{tabular}{cccccccc}
\hline
$\Delta t$ & $T$ & $N_{\rm elements} = 1/h$ & $L$ & $\rho A = \rho I$ & $EA=GA$ & $EI$ \\
\hline
0.01 & 1 & 8 & 10 & 0 & $10^4$ & 500 \\
\hline
\end{tabular}
\caption{Parameters for roll up example}
\label{tab:parameters_rollup}
\end{table}
\begin{figure}[htb]
    \centering
    \begin{subfigure}[b]{0.32\textwidth}
        \centering
        \includegraphics[width=\textwidth]{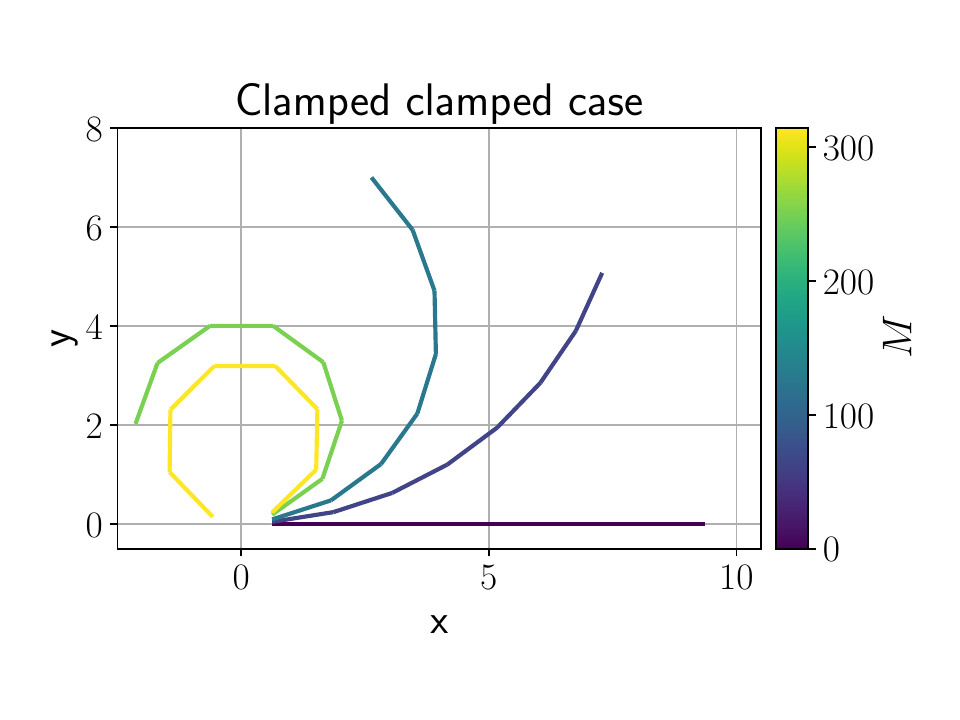}
        \caption{Clamped-clamped}
    \label{fig:clamped_nonlinear_beam}
    \end{subfigure}
    \hfill
    \begin{subfigure}[b]{0.32\textwidth}
        \centering
        \includegraphics[width=\textwidth]{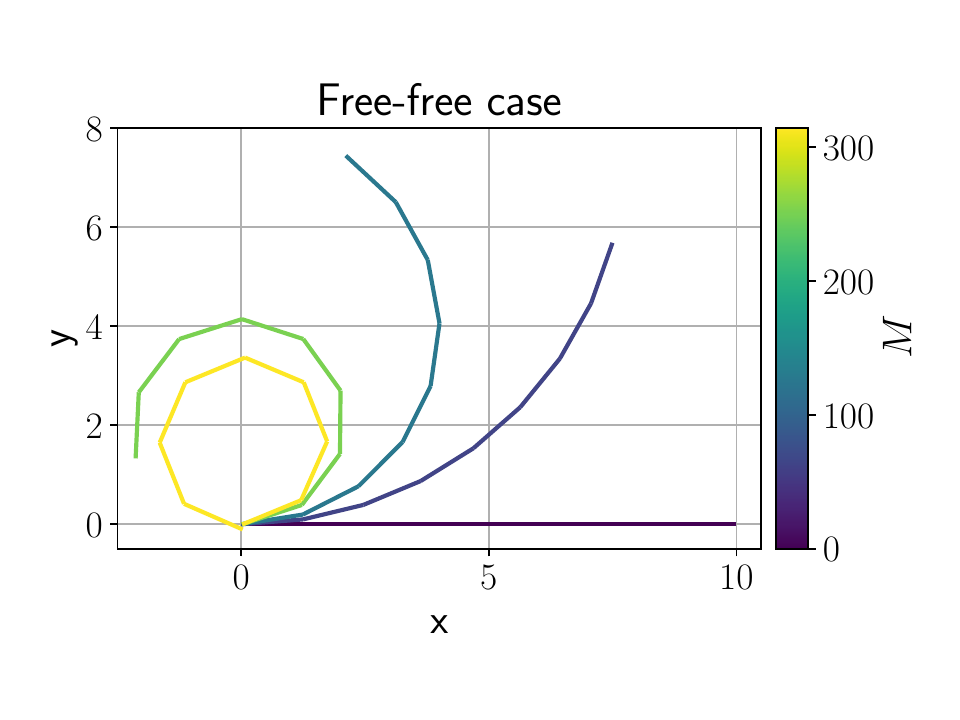}
        \caption{Free-free}
        \label{fig:free_nonlinear_beam}
    \end{subfigure}
    \hfill
    \begin{subfigure}[b]{0.32\textwidth}
        \centering
        \includegraphics[width=\textwidth]{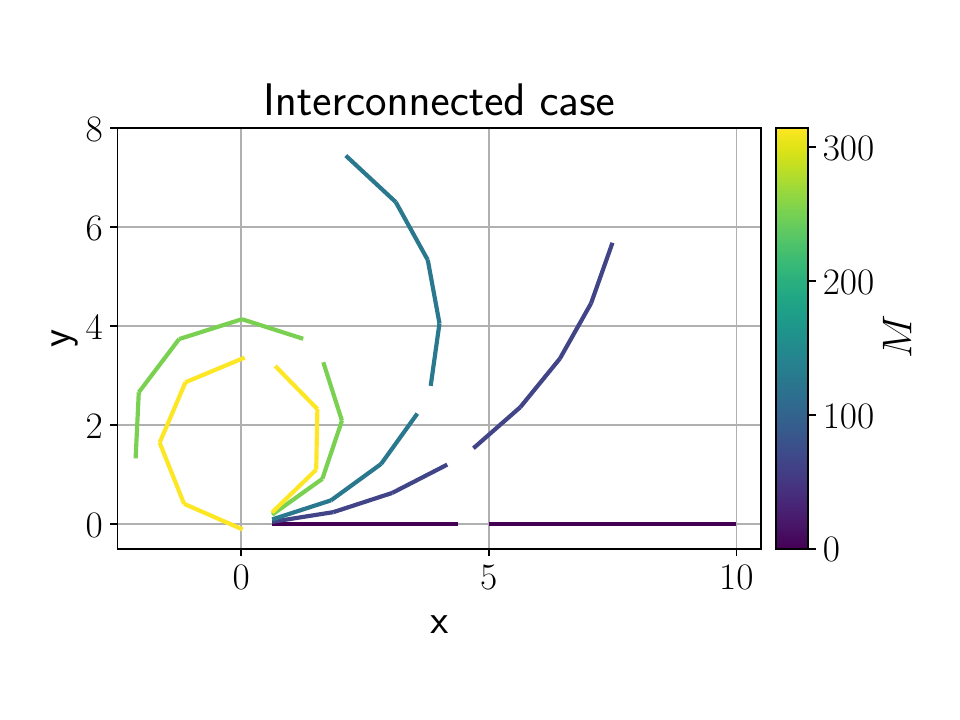}
        \caption{Interconnection}
    \label{fig:interonnected_nonlinear_beam}
    \end{subfigure}
    \caption{Configurations for the quasi static roll up using the free-free, clamped-clamped and interconnected model. System \eqref{eq:discrete_clamped_intrinsic} is used for Fig. \ref{fig:clamped_nonlinear_beam}. System \eqref{eq:discrete_free_intrinsic} is used for Fig. \ref{fig:free_nonlinear_beam}. Fig. \ref{fig:interonnected_nonlinear_beam} shows the results when using the interconnection of the two. The first two cases require a Lagrange multiplier to impose the boundary conditions while the last does not.}
    \label{fig:configurations_rollup}
\end{figure}

%% file: Text/pH_WaveEquation.tex
\subsection{The Wave Equation in 2D}\label{sec:Wave}
The method is now applied to the two dimensional wave equation on a unit square domain
\begin{equation}\label{eq:wave_primal}
\partial_{tt}\phi-\div\grad\phi=0, \qquad \Omega = [0, 1]^ 2,    
\end{equation}
split into subdomains $\Omega_1$ and $\Omega_2$ with a Dirichlet boundary condition on $\Gamma_1$ and a Neumann boundary condition on $\Gamma_2$. The discretization of the wave equation starts again by expressing via the variables
\begin{equation}
    {e}_{\alpha} = \partial_t\phi, \qquad \bm{e}_{\beta} = \grad \phi.
\end{equation}
The system equivalent to \eqref{eq:energy_co-energy-system} is then written as
\begin{equation}
    \begin{pmatrix}
        \partial_t {e}_{\alpha} \\
        \partial_t \bm{e}_{\beta}
    \end{pmatrix} = 
    \begin{bmatrix}
        0 & \div \\
        \grad & 0
    \end{bmatrix}
    \begin{pmatrix}
        {e}_{\alpha} \\
        \bm{e}_{\beta}
    \end{pmatrix},
\end{equation}
so the differential operator for the wave equations is $\mathcal{L}=\grad$ and its formal adjoint is now $\mathcal{L}^*=-\div$. The discretization is obtained by multiplying by the test functions and applying integration by parts as in \eqref{eq:weak_conforming_L} or \eqref{eq:weak_conforming_L*}. The resulting weak formulation for $\Omega_1$ is to find $\bm{e}_\alpha\in L^2(\Omega_1)$ and $\bm{e}_\beta\in H^{\div}(\Omega_1)$ such that $\forall{v}_\alpha\in L^2(\Omega_1)$ and $\forall\bm{v}_\beta\in H^{\div}(\Omega_1)$ it holds
\begin{equation}
    \begin{split}
    \begin{alignedat}{2}
        \innerproduct[\Omega_1]{{v}_\alpha}{\partial_t{e}_\alpha} &= + \innerproduct[\Omega_1]{{v}_\alpha}{\div\bm{e}_\beta}, \\
        \innerproduct[\Omega_1]{\bm{v}_\beta}{\partial_t\bm{e}_\beta} &= - \innerproduct[\Omega_1]{\div\bm{v}_\beta}{{e}_\alpha} + \boundary[\Gamma_1]{\mathcal{T}_{\beta}\bm{v}_\beta}{{u}_{\partial,1 }} + \boundary[\Gamma_{\rm int}]{\mathcal{T}_{\beta}\bm{v}_\beta}{u_{\partial, 1}^{\Gamma_{\rm int}}},
    \end{alignedat}
    \end{split}
\label{eq:Wave-weak-D}
\end{equation}
where $\mathcal{T}_\beta \bm{g} = \bm{g} \cdot \bm{n}|_{\partial\Omega_1}$ is the normal trace. \revtwo{For the subdomain with the Neumann boundary condition $\Omega_2$, seek $\bm{e}_\alpha\in H^1(\Omega_2)$ and $\bm{e}_\beta\in H^{\curl}(\Omega_2)$ to satisfy $\forall{v}_\alpha\in H^1(\Omega_2)$ and $\forall\bm{v}_\beta\in H^{\curl}(\Omega_2)$}
\begin{equation}
    \begin{split}
    \begin{alignedat}{2}
        \innerproduct[\Omega_2]{{v}_\alpha}{\partial_t{e}_\alpha} &= -\innerproduct[\Omega_2]{\grad{v}_\alpha}{\bm{e}_\beta} + \boundary[\Gamma_2]{\mathcal{T}_{\alpha}v_\alpha}{{u}_{\partial,2}} + \boundary[\Gamma_{\rm int}]{\mathcal{T}_{\alpha}{v}_\alpha}{{u}_{\partial,2}^{\Gamma_{\rm int}}}, \\
        \innerproduct[\Omega_2]{\bm{v}_\beta}{\partial_t\bm{e}_\beta} &= \innerproduct[\Omega_2]{\bm{v}_\beta}{\grad {e}_\alpha},
    \end{alignedat}
    \end{split}
\label{eq:Wave-weak-N}
\end{equation}
where $\mathcal{T}_\alpha f = f|_{\partial\Omega_2}$ is the Dirichlet trace. 
\paragraph{Finite element spaces} The mesh consists of a structured triangular mesh. Discontinuous Galerkin of order $k-1$ and Raviart-Thomas of order $k$ (RT$_k$) elements are used for $e_{\alpha, 1}$ and $\bm{e}_{\beta, 1}$ respectively on the $\Omega_1$ subdomain. Continuous Galerkin of order $k$ (CG$_k$) element for $e_{\alpha, 2}$, the N\'ed\'elec first kind of order $k$ (NED$_k$) for $\bm{e}_{\beta, 2}$ on the $\Omega_2$ subdomain.  The justification for this choice comes from de Rham complex and the subcomplex obtained using finite element differential forms of the trimmed polynomial family. The corresponding complex and subcomplex are given by

\begin{minipage}[t]{0.45\textwidth}
\begin{center}
\begin{tikzcd}
H^{\div} \arrow[r, "\div"] \arrow[d, "\Pi"] & L^2 \arrow[d, "\Pi"] \\
\mathrm{RT}_k \arrow[r, "\div"] & \mathrm{DG}_{k-1}
\end{tikzcd}  
\end{center}
\end{minipage} 
\begin{minipage}[t]{0.45\textwidth}
\begin{center}
\begin{tikzcd}
H^{1} \arrow[r, "\grad"] \arrow[d, "\Pi"] & H^{\curl} \arrow[d, "\Pi"] \\
\mathrm{CG}_k \arrow[r, "\grad"] & \mathrm{NED}_k
\end{tikzcd}  
\end{center}
\end{minipage}

\noindent The solution is again found on union of meshes, that is $\mathfrak{T}_h=\mathfrak{T}_h^{\Omega_1}\cup\mathfrak{T}_h^{\Omega_2}$, with finite dimensional spaces for the $\Omega_1$ subdomain given by
\begin{equation}
    \begin{split}
         V_{\alpha, 1} &= \{u_h\in L^2(\Omega_1)|\; \forall T\in\mathfrak{T}_h^{\Omega_1}, \; u_h|_T\in \mathrm{DG}\}, \\
         V_{\beta, 1} &= \{\bm{u}_h\in H^{\div}(\Omega_1)|\; \forall T\in\mathfrak{T}_h^{\Omega_1}, \; \bm{u}_h|_T\in \mathrm{RT}\}, 
    \end{split}
\end{equation}
where $T$ now denotes a triangular mesh element of $\mathfrak{T}_h$. For the $\Omega_2$ subdomain the mixed finite element spaces are
\begin{equation}
    \begin{split}
        V_{\alpha, 2} &= \{u_h\in H^1(\Omega_2)|\; \forall T\in\mathfrak{T}_h^{\Omega_2},\; u_h|_T\in \mathrm{CG}\}, \\
        V_{\beta, 2} &= \{\bm{u}_h\in H^{\curl}(\Omega_2)|\; \forall T\in\mathfrak{T}_h^{\Omega_2}, \bm{u}_h|_T\in \mathrm{Ned}\}.
    \end{split}
\end{equation}
The finite dimensional system for the $\Omega_1$ subdomain becomes
\begin{equation}
    \begin{split}
        \begin{bmatrix}
            \mathbf{M}_{\alpha, 1} & 0 \\
            0 & \mathbf{M}_{\beta, 1}
        \end{bmatrix}\odv{}{t}
        \begin{pmatrix}
            \mathbf{e}_{\alpha,1} \\
            \mathbf{e}_{\beta,1}
        \end{pmatrix} &= 
        \begin{bmatrix}
            0 & \mathbf{D}_{\div} \\
            -\mathbf{D}_{\div}^\top & 0 
        \end{bmatrix}
        \begin{pmatrix}
            \mathbf{e}_{\alpha,1} \\
            \mathbf{e}_{\beta,1}
        \end{pmatrix} + 
        \begin{bmatrix}
            0 & 0 \\
            \mathbf{B}_{\beta}^{\Gamma_1} & \mathbf{B}_{\beta}^{\Gamma_{\rm int}}
        \end{bmatrix}
        \begin{pmatrix}
            \mathbf{u}_{\partial,1} \\
            \mathbf{u}_{\partial,1}^{\rm int}
        \end{pmatrix}, \\
        \begin{pmatrix}
            \mathbf{y}_{\partial,1} \\
            \mathbf{y}_{\partial,1}^{\rm int}
        \end{pmatrix} &= 
        \begin{bmatrix}
            0 & \mathbf{T}_\beta^{\Gamma_1} \\
            0 & \mathbf{T}_\beta^{\Gamma_{\rm int}}
        \end{bmatrix}
        \begin{pmatrix}
            \mathbf{e}_{\alpha,1} \\
            \mathbf{e}_{\beta,1}
        \end{pmatrix},
    \end{split}
\end{equation}
while for the $\Omega_2$ subdomain it becomes
\begin{equation}
    \begin{split}
        \begin{bmatrix}
            \mathbf{M}_{\alpha, 2} & 0 \\
            0 & \mathbf{M}_{\beta, 2}
        \end{bmatrix}\odv{}{t}
        \begin{pmatrix}
            \mathbf{e}_{\alpha,2} \\
            \mathbf{e}_{\beta,2}
        \end{pmatrix} &= 
        \begin{bmatrix}
            0 & -\mathbf{D}_{\grad} \\
            \mathbf{D}_{\grad}^\top & 0 
        \end{bmatrix}
        \begin{pmatrix}
            \mathbf{e}_{\alpha,2} \\
            \mathbf{e}_{\beta,2}
        \end{pmatrix} + 
        \begin{bmatrix}
            \mathbf{B}_\alpha^{\Gamma_2} & \mathbf{B}_\alpha^{\Gamma_{\rm int}} \\
            0 & 0
        \end{bmatrix}
        \begin{pmatrix}
            \mathbf{u}_{\partial,2} \\
            \mathbf{u}_{\partial,2}^{\rm int}
        \end{pmatrix}, \\
        \begin{pmatrix}
            \mathbf{y}_{\partial,2} \\
            \mathbf{y}_{\partial,2}^{\rm int}
        \end{pmatrix} &= 
        \begin{bmatrix}
            \mathbf{T}_\alpha^{\Gamma_2} & 0 \\
            \mathbf{T}_\alpha^{\Gamma_{\rm int}} & 0
        \end{bmatrix}
        \begin{pmatrix}
            \mathbf{e}_{\alpha,2} \\
            \mathbf{e}_{\beta,2}
        \end{pmatrix}.
        \end{split}
\end{equation}

\paragraph{Numerical experiments}
\begin{figure}
    \centering
    \includegraphics[width=0.3\linewidth]{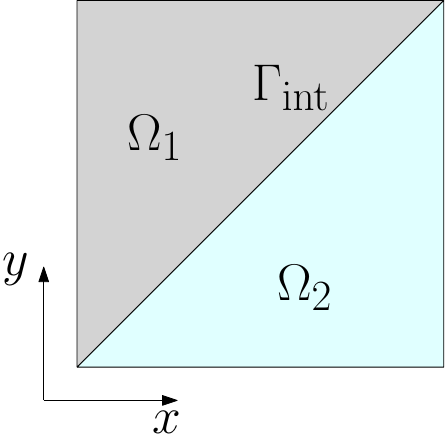}
    \caption{Domain splitting for the wave equation }
    \label{fig:mesh_wave}
\end{figure}

For this example we consider three different numerical analyses:
\begin{itemize}
    \item A convergence study;
    \item A conservation properties study. In particular, a curl free condition of the two subdomains and the power balance will be investigated.
    \item A spectral analysis.
\end{itemize}
The simulations take place on a unit square $\Omega=[0, 1]^2$, decomposed by an interface placed diagonally between the lower left and upper right vertex as shown in Fig \ref{fig:mesh_wave}. \revone{For the first two analysis the St\"ormer-Verlet method has been used with time step time step of $\Delta t=0.001 \; \mathrm{[s]}$.}

\paragraph{Time domain simulations}

\begin{figure}[htbp]
	\centering
	\begin{subfigure}{0.48\textwidth}
		\centering
		\includegraphics[width=\textwidth]{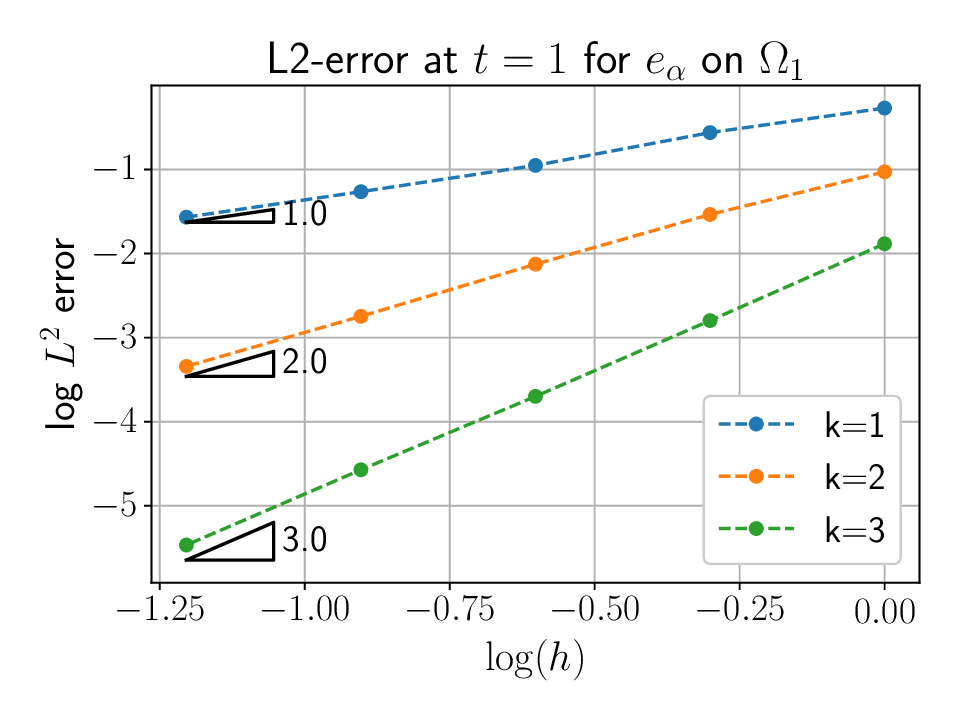}
	\end{subfigure}
	\begin{subfigure}{0.48\textwidth}
		\centering
		\includegraphics[width=\textwidth]{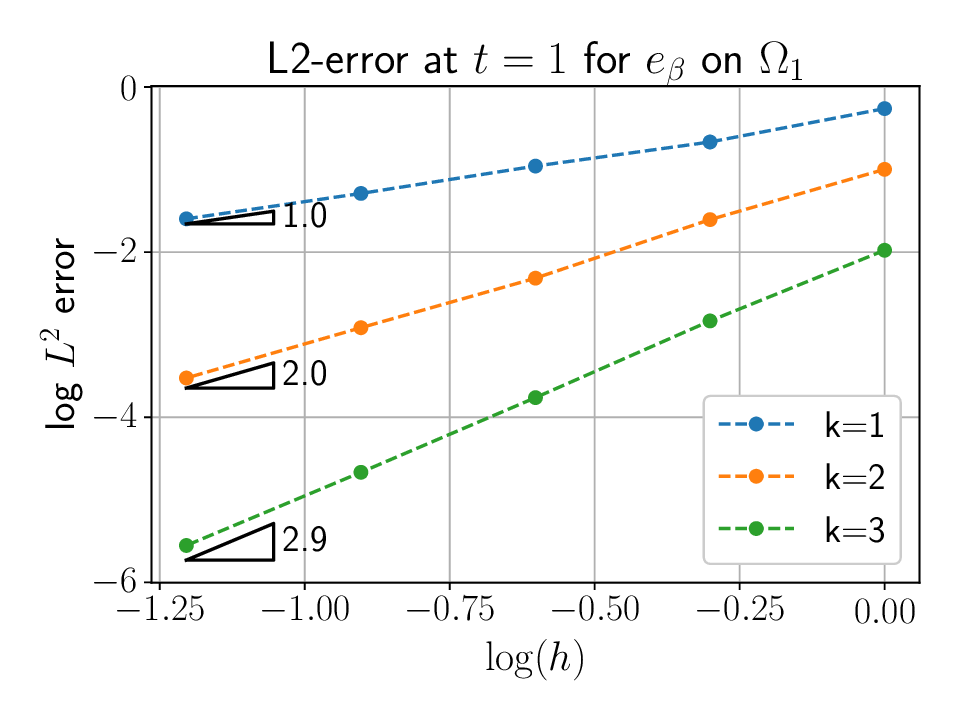}
	\end{subfigure}
        \begin{subfigure}{0.48\textwidth}
		\centering
		\includegraphics[width=\textwidth]{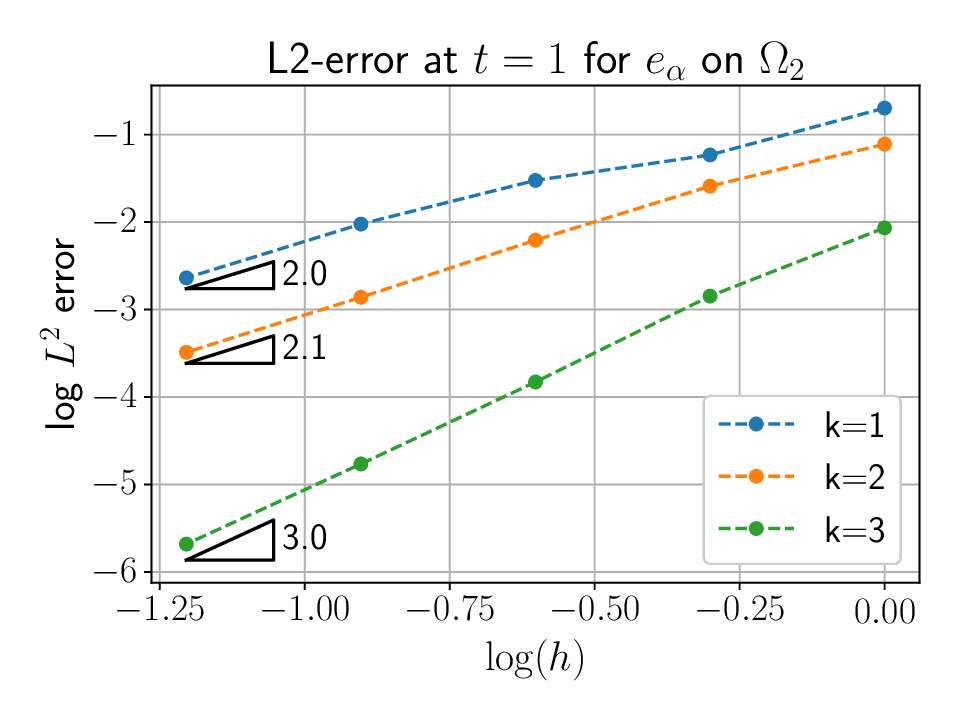}
	\end{subfigure}
	\begin{subfigure}{0.48\textwidth}
		\centering
		\includegraphics[width=\textwidth]{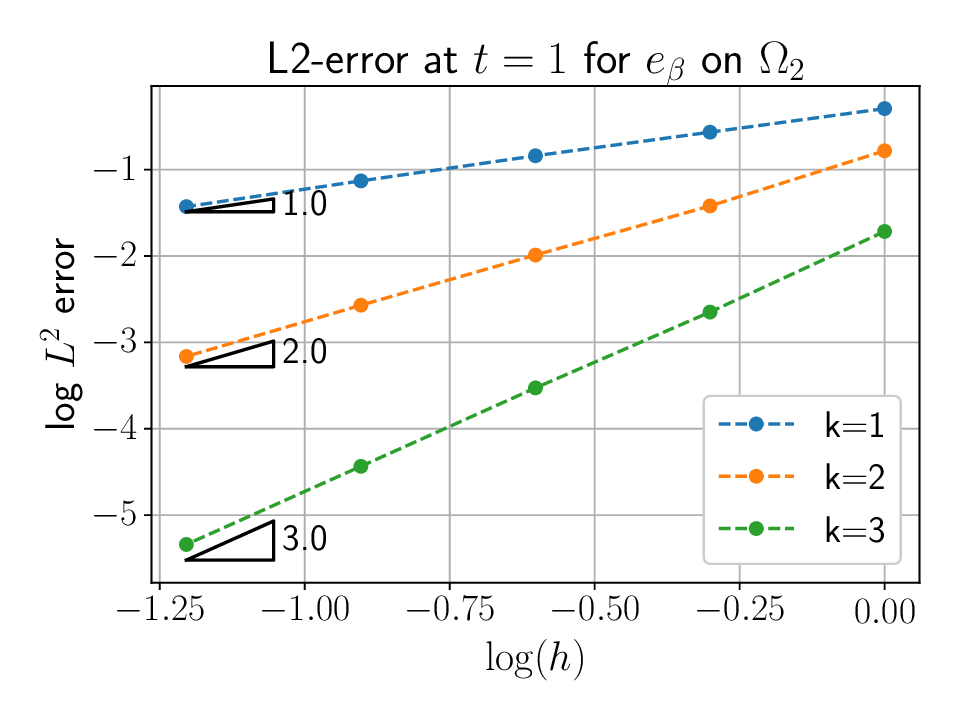}
	\end{subfigure}
	\caption{Convergence rates for wave equation.}
	\label{fig:staggered_convergence}
\end{figure}

\begin{figure}[htb]
    \centering

    \begin{minipage}{0.32\textwidth}
        \centering
        \includegraphics[width=\textwidth]{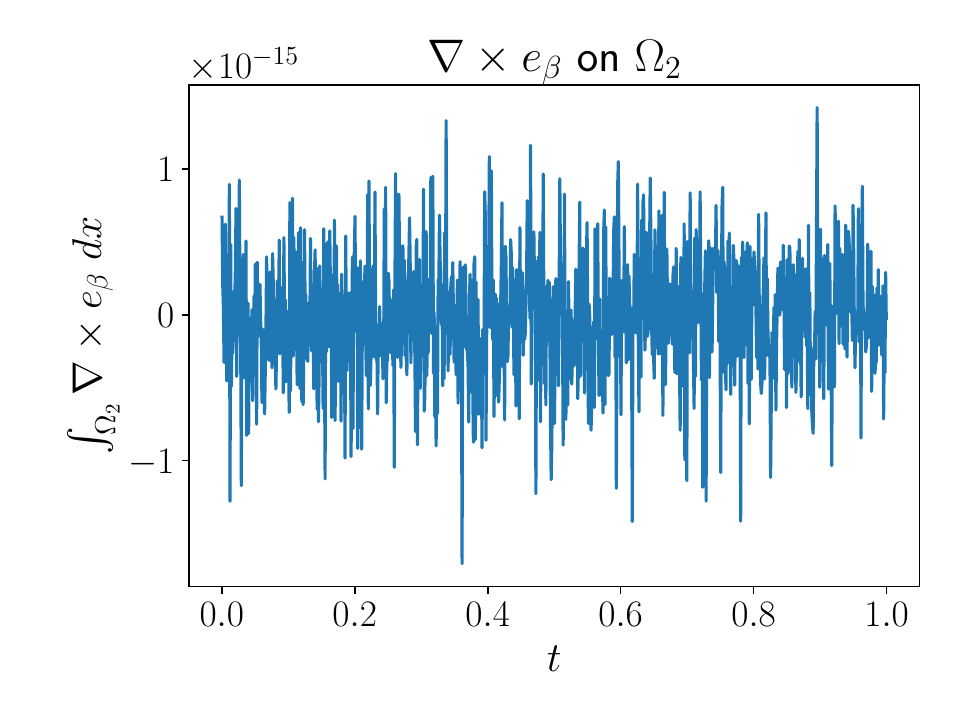}
        \captionof{figure}{$\curl e_{\beta, 2}$ on $\Omega_2$.}
        \label{fig:curl10_staggered}
    \end{minipage}
    \hfill
    \begin{minipage}{0.65\textwidth}
        \centering
        \begin{subfigure}{0.49\textwidth}
            \centering
            \includegraphics[width=\textwidth]{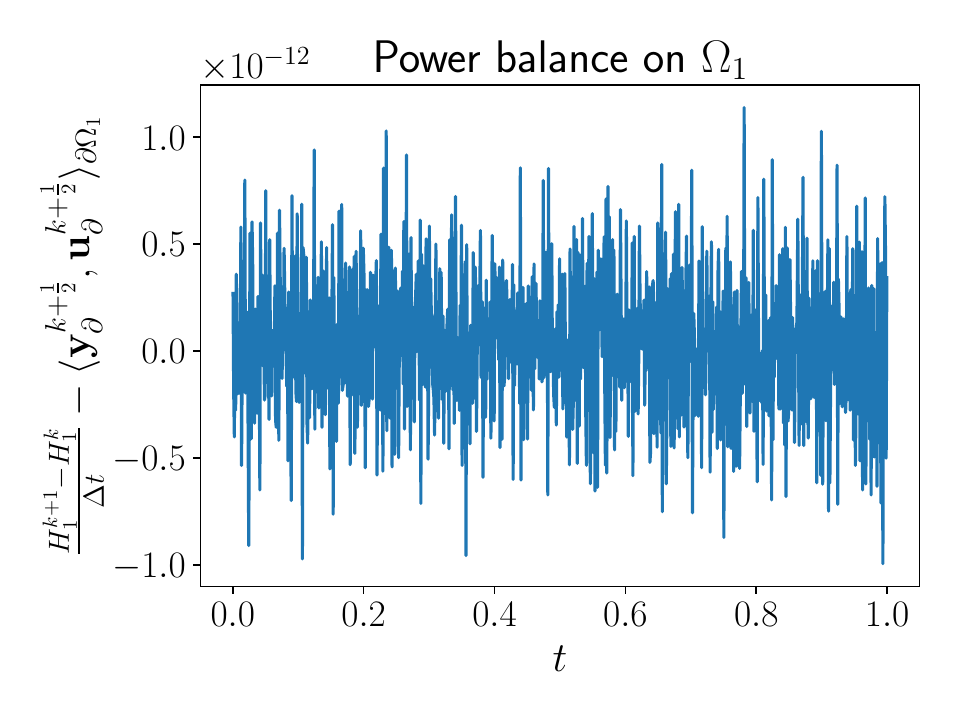}
        \end{subfigure}
        \hfill
        \begin{subfigure}{0.49\textwidth}
            \centering
            \includegraphics[width=\textwidth]{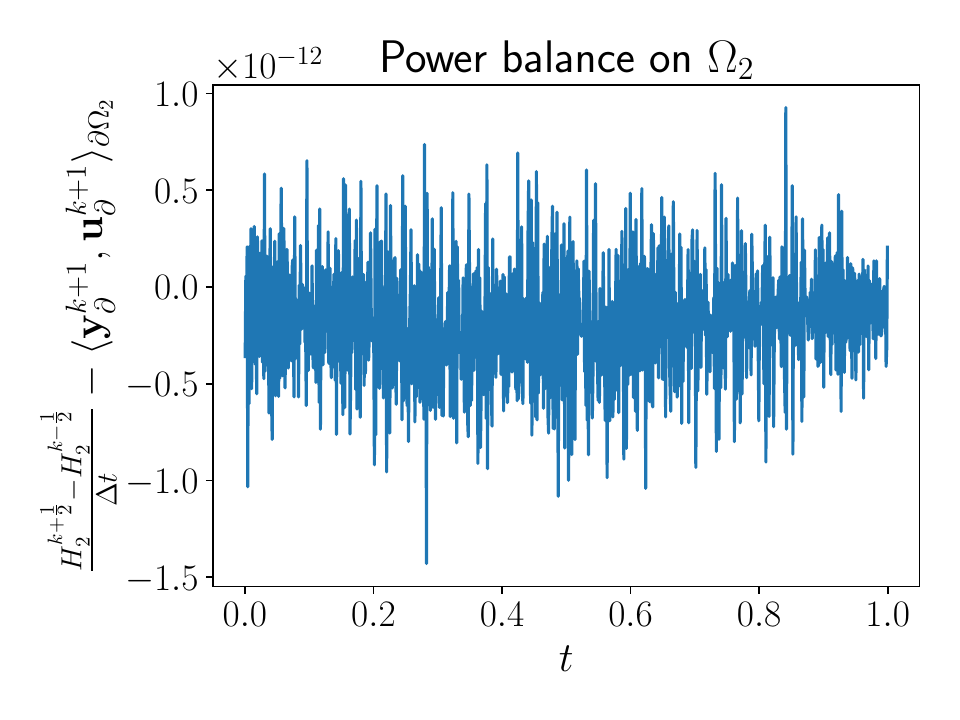}
        \end{subfigure}

        \captionof{figure}{The power balance for the wave equation in two dimensions.}
        \label{fig:Powerbalance_Staggered}
    \end{minipage}

\end{figure}

An analytical solution has been used for the boundary inputs and the verification of the simulations. The exact solution consists of a temporal and spatial part given by
\begin{align}
    f(t) &= 2\sin(\sqrt{2} t)+3\cos(\sqrt{2}t), \\
    g(x,y) &= \cos(x)\sin(y).
\end{align}
The exact solutions are given as
\begin{equation}
        e_\alpha^{ex} = g\odv{f}{t}, \qquad \bm{e}_\beta^{ex} = f\grad g,
\end{equation}
The boundary conditions have been obtained from the exact solutions
$$e_\alpha|_{\Gamma_1}=g\odv{f}{t}, \qquad \bm{e}_\beta\cdot \bm{n}|_{\Gamma_2}=f\, \nabla_{\bm{n}} g|_{\Gamma_2},$$
where $\bm{n}$ denotes the outward unit normal. The spatial convergence has been investigated by performing simulations for five different spatial step sizes $h$ and three polynomial degrees $k=1,2,3$ for a total of 15 simulations. The convergence rates for the mixed finite element formulation are well-known \cite{badia2014stability,grote2006discontinuous}, and are thus expected to have a theoretical convergence rate of $h^k$, apart from $e_{\alpha,2}$. This is due to the fact that $e_{\alpha,2}$ is discretized with a Lagrange element (whose convergence is given by $h^{k+1}$ in the $L^2$ norm) whereas $e_{\alpha,2}$ it is discretized by a discontinuous Galerkin element (that convergences with a rate $h^{k}$). Fig. \ref{fig:staggered_convergence} shows the L2-error of $e_\alpha$ and $e_\beta$ for both the $\Omega_1$ and $\Omega_2$ subdomains. The error is lower with smaller values of $h$ and decreases faster with higher polynomial degrees. The numerical solution is approaching the exact solution with a regular rate. This rate of convergence matches well with the theoretical convergence rate $h^k$, but, as expected, the convergence of $e_\alpha$ on ${\Omega_2}$ behaves slightly differently. For a polynomial order of $k=1$, it converges with $h^{k+1}$, while for higher polynomial orders it converges with $h^k$.
The second equation in \eqref{eq:Wave-weak-N} is satisfied strongly because of the inclusion $V_{\beta, 2} \subset \grad V_{\alpha, 2}$. This means that the following holds
$$\curl\partial_t \bm{e}_\beta=\curl\grad \bm{e}_\alpha=0.$$ 
The curl free condition is instead only satisfied weakly in Eq. \eqref{eq:Wave-weak-D}. Indeed, suppose $\bm{v}_\beta = \curl \bm{v}$, where $\bm{v}$ is chosen in a N\'ed\'elec space (recall that $\curl \mathrm{Ned}_k \subset \mathrm{RT}_k$ so this is a valid choice of test function), then it holds that
$$
\innerproduct[\Omega_1]{\curl \bm{v}}{\partial_t\bm{e}_\beta} = \innerproduct[\Omega_1]{\div\curl \bm{v}}{{e}_\alpha} = 0,
$$
for vanishing boundary conditions. The $L^2$ norm of $\curl \bm{e}_\beta$ is plotted in Fig.~\ref{fig:curl10_staggered} and it is zero within machine precision. 
\revone{For the time integration, the St\"ormer-Verlet scheme is used and each subdomain satisfies a power balance \cite{kotyczka2019discrete}}. Since the solution for $\Omega_1$ is computed at integer time steps, it holds
$$
\frac{{H}_1^{n+1} - H_1^n}{\Delta t} - \boundary[\partial\Omega_1]{\mathbf{y}_\partial^{n+\frac{1}{2}}}{\mathbf{u}_\partial^{n+\frac{1}{2}}} = 0.
$$
For $\Omega_2$ the solution is advanced at half-integer time steps so
$$
\frac{{H}_2^{n+\frac{1}{2}} - H_2^{n-\frac{1}{2}}}{\Delta t} - \boundary[\partial\Omega_2]{\mathbf{y}_\partial^{n+1}}{\mathbf{u}_\partial^{n+1}} = 0.
$$
The power balance is determined for both the $\Omega_1$ and $\Omega_2$ subdomains in Figs.~\ref{fig:Powerbalance_Staggered}. For both parts of the domain the power balance is observed to be in the order of $10^{-12}$, hence zero within machine precision. If the entirer domain $\Omega= \Omega_1 \cup \Omega_2$ was considered than the power balance would not be preserved to machine precision. 

\begin{figure}
	\centering
	\begin{subfigure}{0.32\textwidth}
		\centering
		\includegraphics[width=\textwidth]{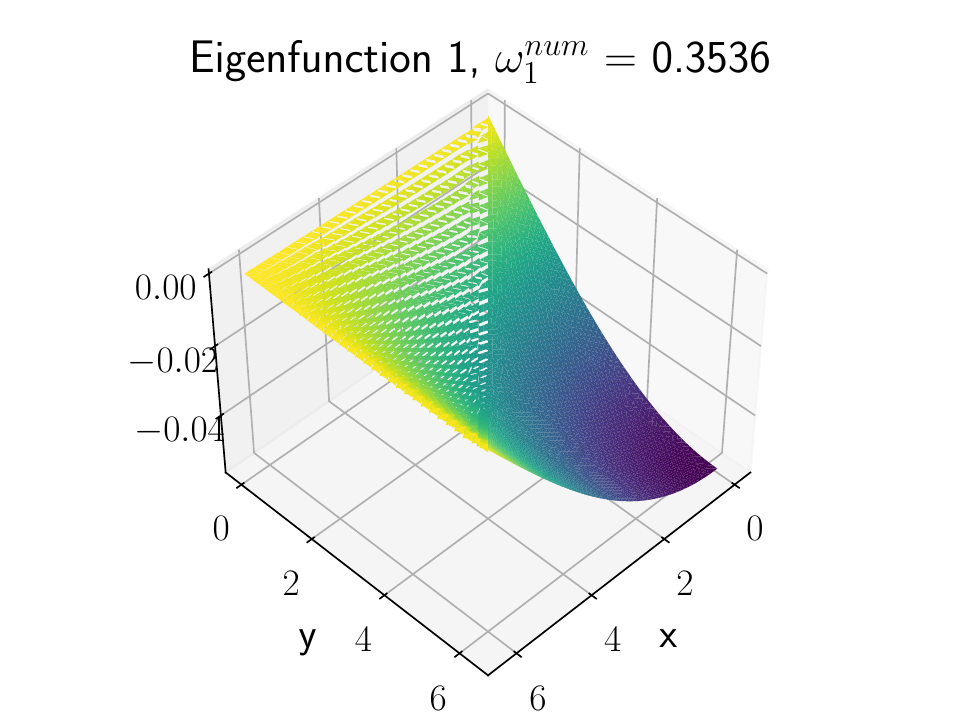}
	\end{subfigure}
	\begin{subfigure}{0.32\textwidth}
		\centering
		\includegraphics[width=\textwidth]{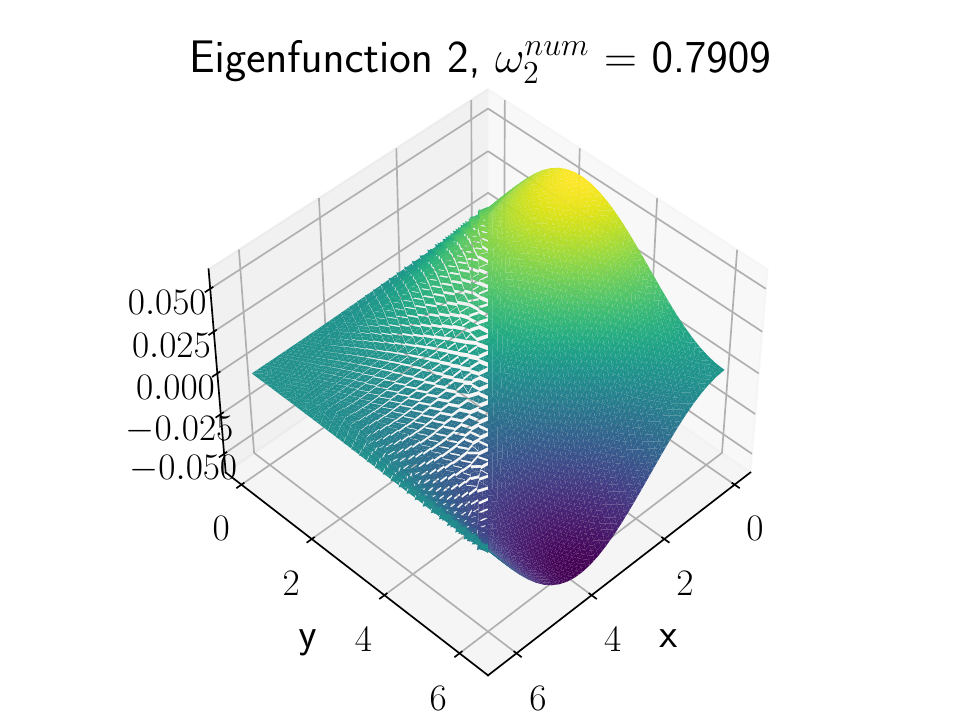}
	\end{subfigure}
        \begin{subfigure}{0.32\textwidth}
		\centering
		\includegraphics[width=\textwidth]{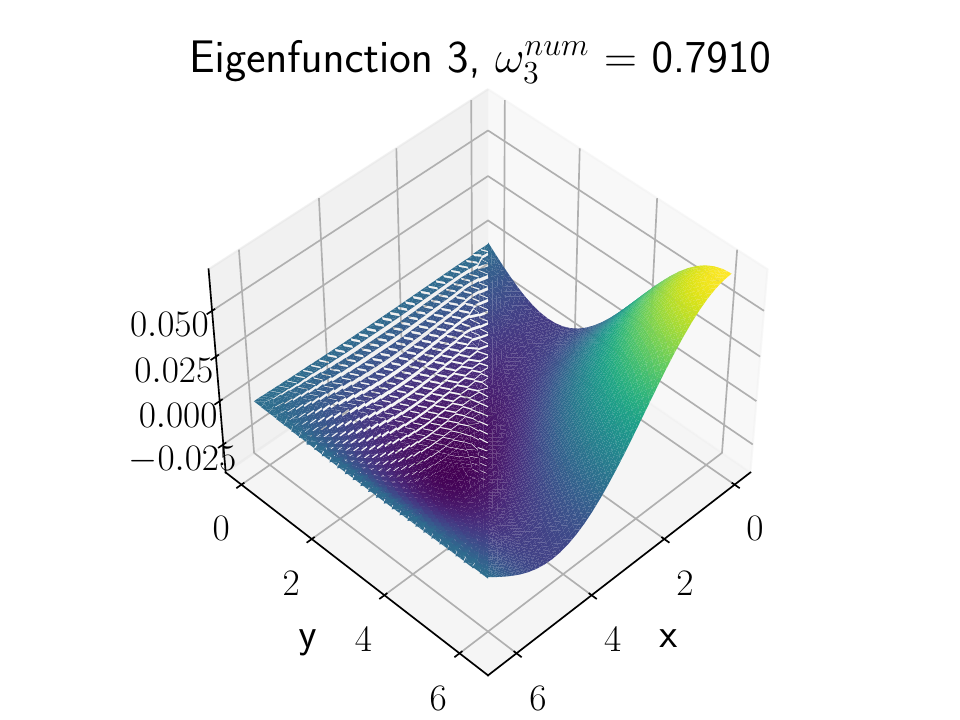}
	\end{subfigure}
	\begin{subfigure}{0.32\textwidth}
		\centering
		\includegraphics[width=\textwidth]{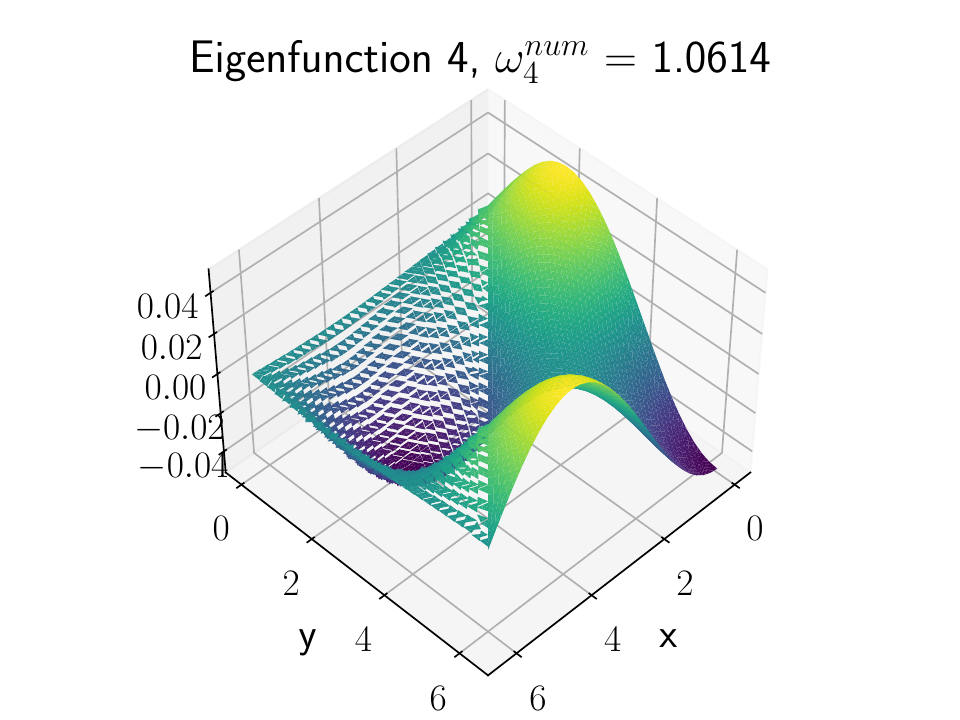}
	\end{subfigure}
        \begin{subfigure}{0.32\textwidth}
		\centering
		\includegraphics[width=\textwidth]{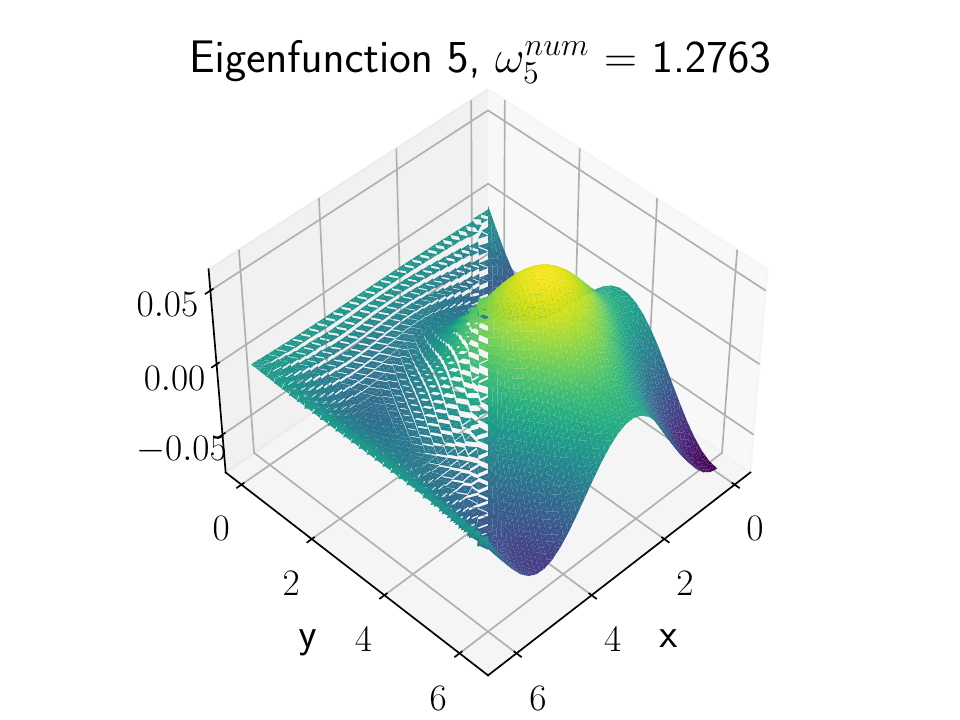}
	\end{subfigure}
	\begin{subfigure}{0.32\textwidth}
		\centering
		\includegraphics[width=\textwidth]{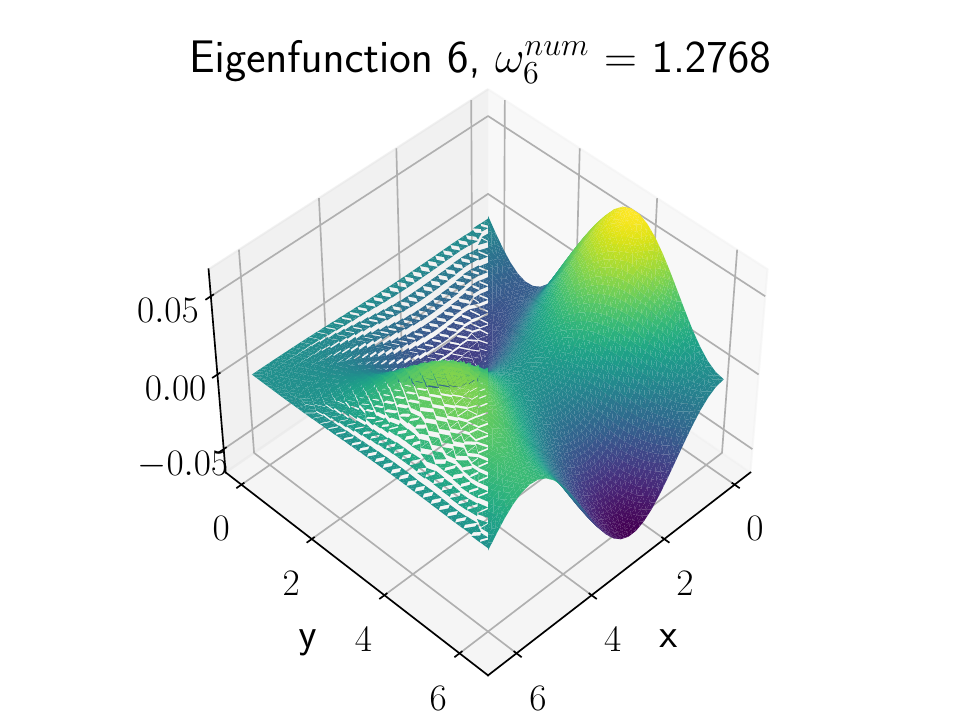}
	\end{subfigure}
	\caption{Numerical eigenvectors for variable $e_\alpha$ wave equation under mixed boundary conditions}
	\label{fig:wave_eigenvectors}
\end{figure}

\paragraph{Modal analysis}
The analytical eigenvalues for the problem under examination are obtained via separation of variables, leading to the following analytical eigenvalues 
$$
\omega_{mn}^{ana} = \frac{\pi}{2 L}\sqrt{(2m-1)^2 + (2n-1)^2}.
$$
The numerical eigenvalues are obtained via the generalized eigenproblem
$$
i \omega_{mn}^{num} \mathbf{M}\bm{\psi}_{mn} = \mathbf{J}\bm{\psi}_{mn},
$$
where $i= \sqrt{-1}$ is the imaginary unit.
For the discretization $30$ finite elements per side are considered. The resulting eigenvalues are shown in Table \ref{tab:eigenvalues_wave} whereas the eigenvectors are plotted in Fig. \ref{fig:wave_eigenvectors}. The obtained eigenfrequencies match the analytical solution and have comparable accuracy with respect to a classical finite element discretization using linear Lagrange elements CG$_1$, with smaller error overall.
\begin{table}[htb]
\centering
\begin{tabular}{c c c c c c}
\hline
Mode 
& Proposed 
& Classical 
& Analytical
& Rel. Err. Prop. 
& Rel. Err. Class. \\
\hline
1  & 0.3536 & 0.3536 & 0.3536 & 0.002 & 0.017 \\
2  & 0.7910 & 0.7912 & 0.7906 & 0.058 & 0.083 \\
3  & 0.7908 & 0.7912 & 0.7906 & 0.035 & 0.085 \\
4  & 1.0613 & 1.0622 & 1.0607 & 0.068 & 0.153 \\
5  & 1.2763 & 1.2775 & 1.2748 & 0.124 & 0.216 \\
6  & 1.2767 & 1.2775 & 1.2748 & 0.158 & 0.220 \\
\hline
\end{tabular}
\caption{Comparison of numerical and analytical eigenvalues of the 2D wave equation. Relative error is computed as $|\omega_{nm}^{\text{num}} - \omega_{nm}^{\text{ana}}|/{\omega_{nm}^{\text{ana}}} \times 100$.}
\label{tab:eigenvalues_wave}
\end{table}

%% file: Text/pH_elasticity.tex
\subsection{Linear Elastodynamics}\label{sec:elasticity}
We now consider the elastodynamics problem in a unit square domain
\begin{equation}\label{eq:elasticity_primal}
\rho\partial_{tt}\bm{u}-\Div\bm{\sigma}=0, \qquad \Omega = [0, 1]^ 2.
\end{equation}
Here $\rho$ is the density, $\bm{u}$ is the displacement field, $\Div$ is row-wise divergence of a tensor and $\bm{\sigma}$ the Cauchy stress tensor
$$
\bm{\sigma} = \bm{K} \Grad \bm{u},
$$
where $\bm{K}$ is the stiffness tensor and $\Grad \bm{u} = \bm{\varepsilon}:=\frac{1}{2}(\nabla\bm{u} + (\nabla\bm{u})^\top)$ is the infinitesimal strain. An isotropic material under plane stress is considered
$$
\bm{K}(\circ) = \frac{E}{(1 - \nu^2)} ((1 - \nu)\circ + \nu\mathrm{tr}(\circ)\bm{I}_2),
$$
where $E$ is the Young modulus and $\nu$ the Poisson ratio. The domain is split into subdomains $\Omega_1$ and $\Omega_2$ with a Dirichlet boundary condition on $\Gamma_1$ and a Neumann boundary condition on $\Gamma_2$. The discretization of the elastodynamics problem starts via the change of variables
\begin{equation}
     \bm{e}_{\alpha} = \partial_t\bm{u}, \qquad \bm{E}_{\beta} = \bm{\sigma},
\end{equation}
where uppercase is used to denote tensor variables. The system equivalent to \eqref{eq:energy_co-energy-system} is then written as
\begin{equation}
    \begin{bmatrix}
        \rho & 0 \\
        0 & \bm{C}
    \end{bmatrix}
    \begin{pmatrix}
        \partial_t\bm{e}_{\alpha} \\
        \partial_t\bm{E}_{\beta}
    \end{pmatrix} = 
    \begin{bmatrix}
        0 & \Div \\
        \Grad & 0
    \end{bmatrix}
    \begin{pmatrix}
        \bm{e}_{\alpha} \\
        \bm{E}_{\beta}
    \end{pmatrix},
\end{equation}
where $\bm{C} := \bm{K}^{-1}$ is the compliance tensor. The differential operator for the elastodynamics problem is $\mathcal{L}=\Grad$ and its formal adjoint is now $\mathcal{L}^*=-\Div$ applied to symmetric tensors. The discretization is obtained by multiplying by the test functions and applying integration by parts as in \eqref{eq:weak_conforming_L} or \eqref{eq:weak_conforming_L*}. The resulting weak formulation for $\Omega_1$ is to find $\bm{e}_\alpha\in L^2(\Omega_1; \bbR^2)$ (the $L^2$ space of two dimensional vectors) and $\bm{E}_\beta\in H^{\Div}(\Omega_1; \bbS)$, where $\bbS = \mathbb{R}^{2\times 2}_{\rm sym}$ (the space of Div conforming symmetric tensors) such that $\forall{v}_\alpha\in L^2(\Omega_1; \bbR^2)$ and $\forall\bm{v}_\beta\in H^{\Div}(\Omega_1; \bbS)$ it holds
\begin{equation}
    \begin{split}
    \begin{alignedat}{2}
        \innerproduct[\Omega_1]{\bm{v}_\alpha}{\rho \partial_t\bm{e}_\alpha} &=  +\innerproduct[\Omega_1]{\bm{v}_\alpha}{\Div\bm{E}_\beta}, \\
        \innerproduct[\Omega_1]{\bm{V}_\beta}{\bm{C}\partial_t\bm{E}_\beta} &= -\innerproduct[\Omega_1]{\Div\bm{V}_\beta}{\bm{e}_\alpha} + \boundary[\Gamma_1]{\mathcal{T}_{\beta}\bm{V}_\beta}{\bm{u}_{\partial,1 }} + \boundary[\Gamma_{\rm int}]{\mathcal{T}_{\beta}\bm{V}_\beta}{\bm{u}_{\partial, 1}^{\Gamma_{\rm int}}},
    \end{alignedat}
    \end{split}
\label{eq:Elasticity-weak-D}
\end{equation}
where $\mathcal{T}_\beta \bm{S} = \bm{S} \cdot \bm{n}|_{\partial\Omega_1}$ is the normal trace of a tensor (the traction). For the subdomain $\Omega_2$ where the Neumann boundary condition is natural, the functional setting is the following: seek $\bm{E}_\alpha\in H^1(\Omega_2; \mathbb{R}^2)$ and $\bm{E}_\beta\in H^{\rot\Rot}(\Omega_2; \mathbb{S})$ to satisfy $\forall\bm{v}_\alpha\in H^1(\Omega_2; \bbR^2)$ and $\forall\bm{V}_\beta\in H^{\rot\Rot}(\Omega_2; \mathbb{S})$
\begin{equation}
    \begin{split}
    \begin{alignedat}{2}
        \innerproduct[\Omega_2]{\bm{v}_\alpha}{\rho \partial_t\bm{e}_\alpha} &= -\innerproduct[\Omega_2]{\Grad\bm{v}_\alpha}{\bm{E}_\beta} + \boundary[\Gamma_2]{\mathcal{T}_{\alpha}\bm{v}_\alpha}{{u}_{\partial,2}} + \boundary[\Gamma_{\rm int}]{\mathcal{T}_{\alpha}\bm{v}_\alpha}{{u}_{\partial,2}^{\Gamma_{\rm int}}}, \\
        \innerproduct[\Omega_2]{\bm{V}_\beta}{\bm{C}\partial_t\bm{E}_\beta} &= +\innerproduct[\Omega_2]{\bm{V}_\beta}{\Grad\bm{e}_\alpha},
    \end{alignedat}
    \end{split}
\label{eq:Elasticity-weak-N}
\end{equation}
where $\mathcal{T}_\alpha \bm{u} = \bm{u}|_{\partial\Omega_2}$ is the Dirichlet trace. The space $H^{\rot\Rot}(\Omega_2; \mathbb{S})$ is the space of $\rot\Rot$ conforming symmetric tensor, where the $\rot\Rot$ operator (also known as incompatibility operator in mechanics) is given by 
$$
\rot\Rot \bm{S} = \partial_{xx} S_{yy} + \partial_{yy} S_{xx} - 2 \partial_{xy}S_{xy},
$$
or it ca be interpreted equivalently as the double divergence of a rotated second order tensor
$$
\rot\Rot\bm{S} = \div\Div(\bm{JSJ}^\top), \qquad \bm{J}:=\begin{pmatrix}
    0 & 1 \\
    -1 & 0
\end{pmatrix}.
$$
  
\paragraph{Finite element spaces} The mesh consists of a structured triangular mesh. Discontinuous Galerkin of order $1$ (DG$_1$) and conforming Arnold Winther elements \cite{arnold2002elasticity} of degree $3$, denoted by AW$_3$, are used for $\bm{e}_{\alpha, 1}$ and $\bm{E}_{\beta, 1}$ respectively on the $\Omega_1$ subdomain. Continuous Galerkin of order (CG$_2$) element for $\bm{a}_{\alpha, 2}$, and Discontinuous Galerkin of order $1$ (DG$_1$) for $\bm{e}_{\beta, 2}$ on the $\Omega_2$ subdomain.  The justification for this choice comes from the elasticity complex. However, given the fact that $H^{\rot\Rot}(\bbS)$ are yet not available in finite element libraries, discontinuous Galerkin finite elements are used as on simplicial meshes it holds
$$
\Grad \mathrm{CG}_k(\bbR^2) \subseteq \mathrm{DG}_{k-1}(\bbS),
$$
leading to an exact discrete subcomplex. The commuting diagram for the complexes and corresponding subcomplexes is as follows

\begin{minipage}[t]{0.45\textwidth}
\begin{center}
\begin{tikzcd}
H^{\Div}(\bbS) \arrow[r, "\Div"] \arrow[d, "\Pi"] & L^2(\bbR^2) \arrow[d, "\Pi"] \\
\mathrm{AW}_3 \arrow[r, "\Div"] & \mathrm{DG}_{1}(\bbR^2)
\end{tikzcd}  
\end{center}
\end{minipage} 
\begin{minipage}[t]{0.45\textwidth}
\begin{center}
\begin{tikzcd}
H^{1}(\bbR^2) \arrow[r, "\Grad"] \arrow[d, "\Pi"] & H^{\rot\Rot}(\bbS) \arrow[d, "\Pi"] \\
\mathrm{CG}_2(\bbR^2) \arrow[r, "\Grad"] & \mathrm{DG}_1(\bbS)
\end{tikzcd}  
\end{center}
\end{minipage}

\begin{figure}[htb]
    \centering
    \includegraphics[width=0.3\linewidth]{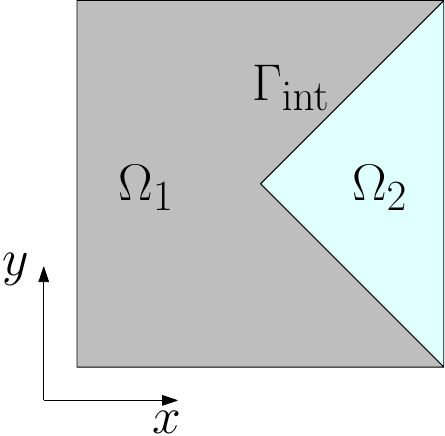}
    \caption{Decomposition of the domain for the elasticity problem.}
    \label{fig:elasticity_mesh}
\end{figure}

\begin{figure}[htbp]
	\centering
	\begin{subfigure}{0.32\textwidth}
		\centering
		\includegraphics[width=\textwidth]{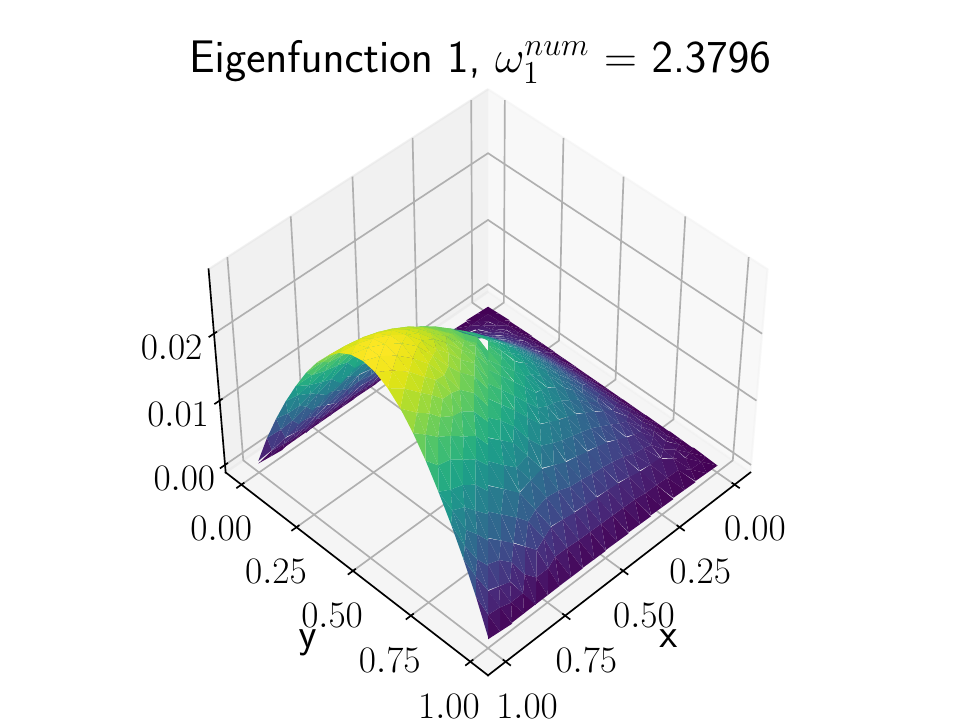}
	\end{subfigure}
	\begin{subfigure}{0.32\textwidth}
		\centering
		\includegraphics[width=\textwidth]{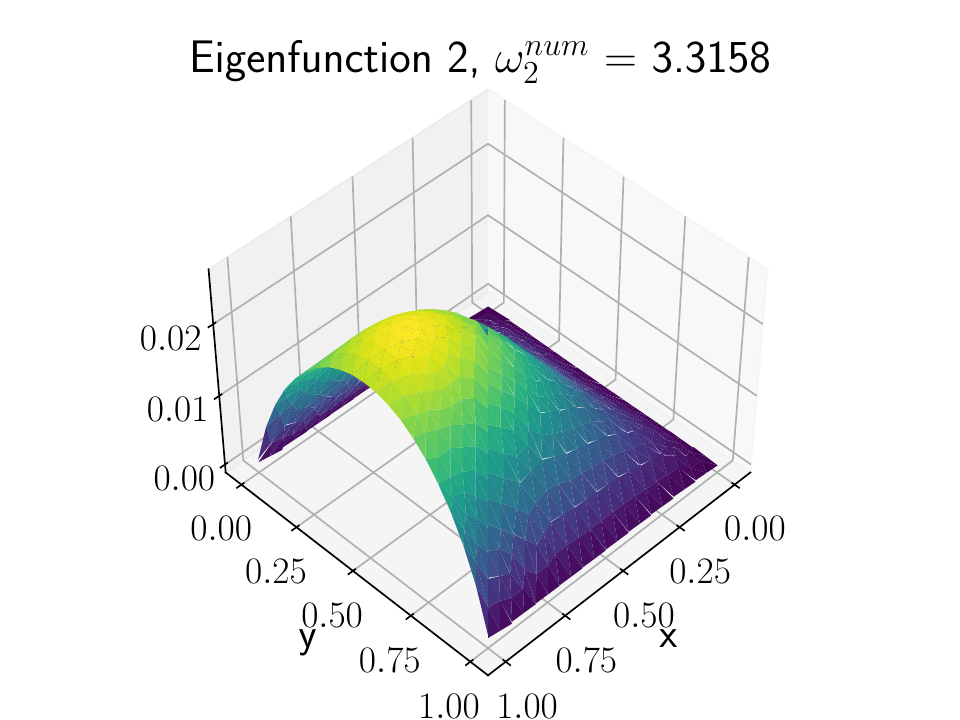}
	\end{subfigure}
        \begin{subfigure}{0.32\textwidth}
		\centering
		\includegraphics[width=\textwidth]{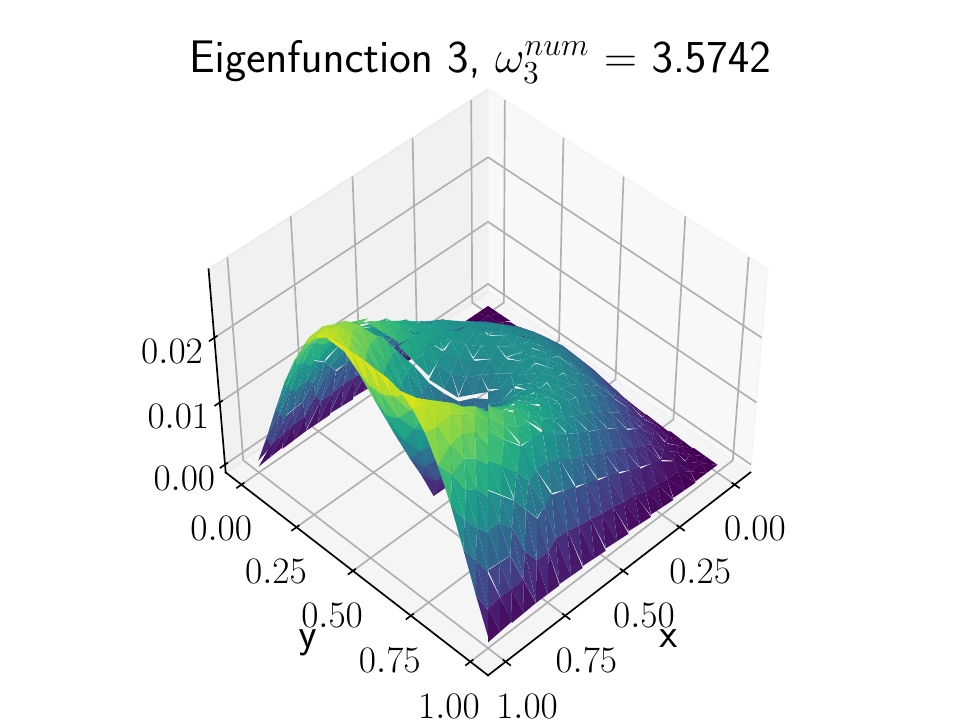}
	\end{subfigure}
	\begin{subfigure}{0.32\textwidth}
		\centering
		\includegraphics[width=\textwidth]{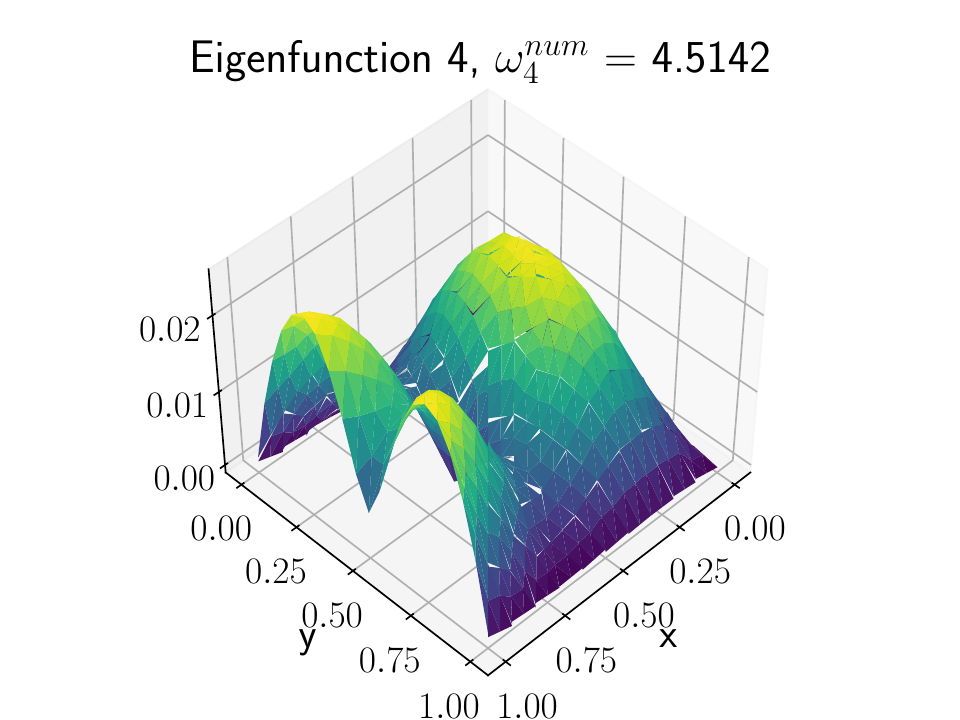}
	\end{subfigure}
        \begin{subfigure}{0.32\textwidth}
		\centering
		\includegraphics[width=\textwidth]{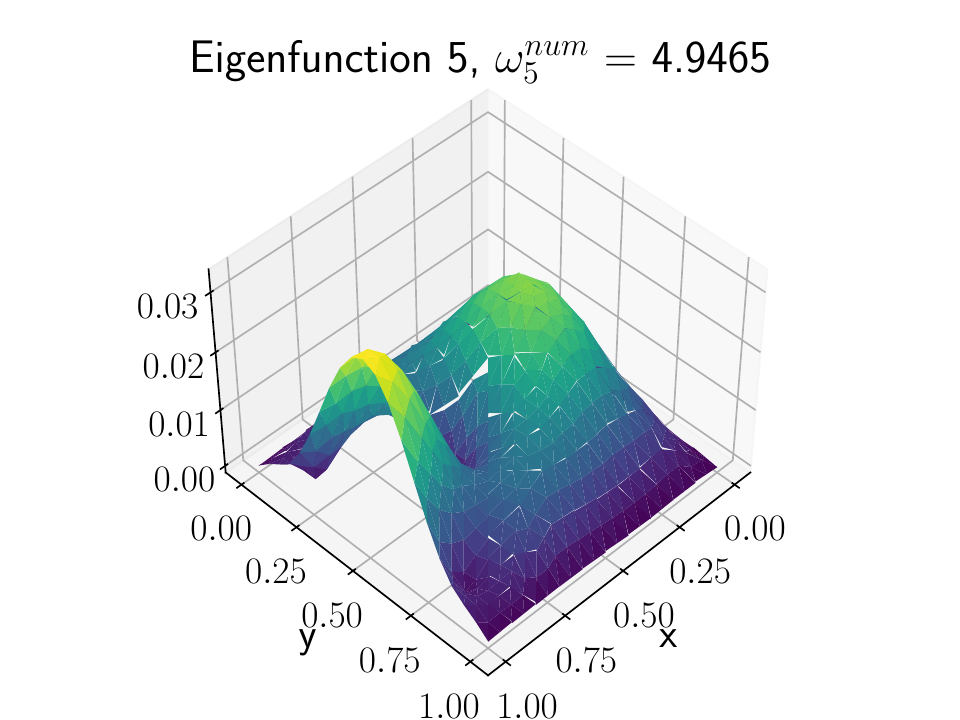}
	\end{subfigure}
	\begin{subfigure}{0.32\textwidth}
		\centering
		\includegraphics[width=\textwidth]{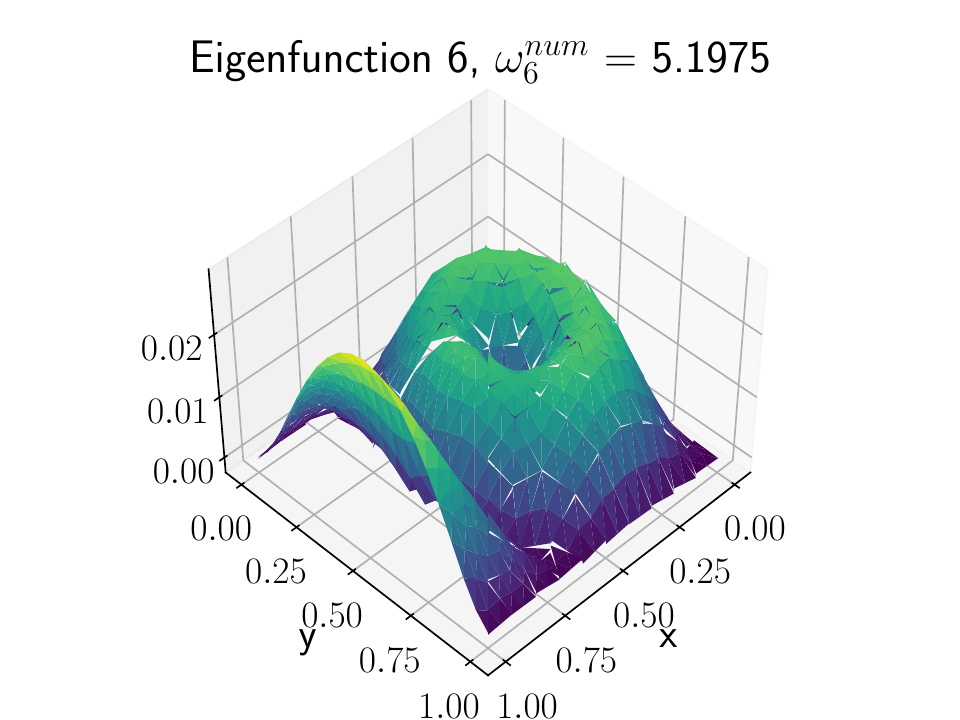}
	\end{subfigure}
	\caption{Numerical eigenvectors in terms of $\sqrt{\bm{e}_\alpha\cdot\bm{e}_\alpha}$ for the elastodynamics problem under mixed boundary conditions}
	\label{fig:elasticity_eigenvectors}
\end{figure}

The solution is found on union of meshes, that is $\mathfrak{T}_h=\mathfrak{T}_h^{\Omega_1}\cup\mathfrak{T}_h^{\Omega_2}$, with finite dimensional spaces for the $\Omega_1$ subdomain given by
\begin{equation}
    \begin{split}
         V_{\alpha, 1} &= \{\bm{u}_h\in L^2(\Omega_1; \bbR^2)|\; \forall T\in\mathfrak{T}_h^{\Omega_1}, \; \bm{u}_h|_T\in \mathrm{DG}_1(\bbR^2)\}, \\
         V_{\beta, 1} &= \{\bm{S}_h\in H^{\Div}(\Omega_1; \bbS)|\; \forall T\in\mathfrak{T}_h^{\Omega_1}, \; \bm{S}_h|_T\in \mathrm{AW}_3\}, 
    \end{split}
\end{equation}
where $T$ now denotes a triangular mesh element of $\mathfrak{T}_h$. For the $\Omega_2$ subdomain the mixed finite element spaces are
\begin{equation}
    \begin{split}
        V_{\alpha, 2} &= \{\bm{u}_h\in H^1(\Omega_2; \bbR^2)|\; \forall T\in\mathfrak{T}_h^{\Omega_2},\; \bm{u}_h|_T\in \mathrm{CG}_2(\bbR^2)\}, \\
        V_{\beta, 2} &= \{\bm{S}_h\in L^2(\Omega_2; \bbS)|\; \forall T\in\mathfrak{T}_h^{\Omega_2}, \bm{S}_h|_T\in \mathrm{DG}_1(\bbS)\}.
    \end{split}
\end{equation}
The finite dimensional system for the $\Omega_1$ subdomain becomes
\begin{equation}
    \begin{split}
        \begin{bmatrix}
            \mathbf{M}_{\alpha, 1} & 0 \\
            0 & \mathbf{M}_{\beta, 1}
        \end{bmatrix}\odv{}{t}
        \begin{pmatrix}
            \mathbf{e}_{\alpha,1} \\
            \mathbf{e}_{\beta,1}
        \end{pmatrix} &= 
        \begin{bmatrix}
            0 & \mathbf{D}_{\Div} \\
            -\mathbf{D}_{\Div}^\top & 0 
        \end{bmatrix}
        \begin{pmatrix}
            \mathbf{e}_{\alpha,1} \\
            \mathbf{e}_{\beta,1}
        \end{pmatrix} + 
        \begin{bmatrix}
            0 & 0 \\
            \mathbf{B}_{\beta}^{\Gamma_1} & \mathbf{B}_{\beta}^{\Gamma_{\rm int}}
        \end{bmatrix}
        \begin{pmatrix}
            \mathbf{u}_{\partial,1} \\
            \mathbf{u}_{\partial,1}^{\rm int}
        \end{pmatrix}, \\
        \begin{pmatrix}
            \mathbf{y}_{\partial,1} \\
            \mathbf{y}_{\partial,1}^{\rm int}
        \end{pmatrix} &= 
        \begin{bmatrix}
            0 & \mathbf{T}_\beta^{\Gamma_1} \\
            0 & \mathbf{T}_\beta^{\Gamma_{\rm int}}
        \end{bmatrix}
        \begin{pmatrix}
            \mathbf{e}_{\alpha,1} \\
            \mathbf{e}_{\beta,1}
        \end{pmatrix},
    \end{split}
\end{equation}
while for the $\Omega_2$ subdomain it becomes
\begin{equation}
    \begin{split}
        \begin{bmatrix}
            \mathbf{M}_{\alpha, 2} & 0 \\
            0 & \mathbf{M}_{\beta, 2}
        \end{bmatrix}\odv{}{t}
        \begin{pmatrix}
            \mathbf{e}_{\alpha,2} \\
            \mathbf{e}_{\beta,2}
        \end{pmatrix} &= 
        \begin{bmatrix}
            0 & -\mathbf{D}_{\Grad} \\
            \mathbf{D}_{\Grad}^\top & 0 
        \end{bmatrix}
        \begin{pmatrix}
            \mathbf{e}_{\alpha,2} \\
            \mathbf{e}_{\beta,2}
        \end{pmatrix} + 
        \begin{bmatrix}
            \mathbf{B}_\alpha^{\Gamma_2} & \mathbf{B}_\alpha^{\Gamma_{\rm int}} \\
            0 & 0
        \end{bmatrix}
        \begin{pmatrix}
            \mathbf{u}_{\partial,2} \\
            \mathbf{u}_{\partial,2}^{\rm int}
        \end{pmatrix}, \\
        \begin{pmatrix}
            \mathbf{y}_{\partial,2} \\
            \mathbf{y}_{\partial,2}^{\rm int}
        \end{pmatrix} &= 
        \begin{bmatrix}
            \mathbf{T}_\alpha^{\Gamma_2} & 0 \\
            \mathbf{T}_\alpha^{\Gamma_{\rm int}} & 0
        \end{bmatrix}
        \begin{pmatrix}
            \mathbf{e}_{\alpha,2} \\
            \mathbf{e}_{\beta,2}
        \end{pmatrix}.
        \end{split}
\end{equation}
\paragraph{Modal analysis} 
For this example we consider only a spectral analysis. The simulations take place on a unit square $\Omega = [0, 1]^2$, decomposed by an interface placed diagonally between the lower left and upper right vertex as shown in Fig. \ref{fig:elasticity_mesh}.
The numerical eigenvalues are obtained via the generalized eigenproblem
$$
i \omega^{num}_n \mathbf{M}\bm{\psi}_n = \mathbf{J}\bm{\psi}_n,
$$
where $i= \sqrt{-1}$ is the imaginary unit.
For the discretization $10$ finite elements per side are considered. The physical and geometrical parameters are 
$$
L=1 \; \mathrm{[m]}, \qquad \rho = 2700 \; \mathrm{[Kg/m^3]}, \qquad E = 70 \; \mathrm{[GPa]}, \qquad \nu=0.3.
$$
The resulting normalized eigenvalues, given by
$$
\widehat{\omega} = \omega \; L\sqrt{\frac{\rho}{E}}, \qquad
$$
are shown in Table \ref{tab:eigenvalues_elasticity} whereas the eigenvectors for variable $\bm{e}_\alpha$ (corresponding to the velocity) are plotted using the magnitude of the associated vector field in Fig. \ref{fig:elasticity_eigenvectors}. The obtained eigenfrequencies have the same accuracy with respect to a classical finite element discretization using quadratic Lagrange elements CG$_2$.

\begin{table}[htbp]
\centering
\begin{tabular}{c c c}
\hline
Mode & Proposed & Classical \\
\hline
1 & 2.3795 & 2.3803 \\
2 & 3.3158 & 3.3166 \\
3 & 3.5742 & 3.5751 \\
4 & 4.5142 & 4.5156 \\
5 & 4.9465 & 4.9468 \\
6 & 5.1975 & 5.1980 \\
\hline
\end{tabular}
\caption{Numerical eigenvalues of the elastodynamics problem using the proposed method and a standard finite element discretization.}
\label{tab:eigenvalues_elasticity}
\end{table}

%% file: Text/pH_Mindlin_plate.tex
\subsection{Mindlin plate}\label{sec:mindlin_plate}
We now consider the Mindlin plate problem in a unit square domain
\begin{equation}\label{eq:mindlin_primal}
\begin{aligned}
\rho h\partial_{tt}{w}-\div\bm{q} &=0, \qquad \Omega = [0, 1]^ 2, \\  
\rho J\partial_{tt}\bm{\bm{\theta}}-(\Div\bm{M} + \bm{q})&=0.
\end{aligned}
\end{equation}
Here $\rho$ is the density, $h$ the thickness, $J:=h^3/12$, $w$ is the vertical displacement, $\bm{\theta}$ the cross section rotation, $\bm{q}$ is shear force resulting and $\bm{M}$ the bending moment tensor. The shear force and bending moment are relative to the kinematic variables $w, \; \bm{\theta}$ via the following relations
$$
\bm{q} = K_{\rm sh} (\grad w - \bm{\theta}), \qquad \bm{M} = \bm{K}_b \Grad \bm{\theta},
$$
where $K_{\rm sh}$ is the shear rigidity and $\bm{K}_b$ the bending stiffness tensor.
For isotropic materials these parameters take the following expressions
$$
K_{\rm sh} = k G h, \qquad 
\bm{K}_b(\circ) = \frac{Eh^3}{12(1 - \nu^2)} ((1 - \nu)\circ + \nu\mathrm{tr}(\circ)\bm{I}_2)
$$
where $k$ is shear correction factor (that depends on the considered boundary conditions), $G = E/(2(1+\nu))$ the shear modulus, $E$ the Young modulus and $\nu$ the Poisson ratio. The discretization of the elastodynamics problem starts via the change of variables
\begin{equation}
    \bm{e}_{\alpha} = \begin{pmatrix}
        \partial_t w \\
        \partial_t \bm{\theta} \\
    \end{pmatrix} = \begin{pmatrix}
        v \\
        \bm{\omega} \\
    \end{pmatrix}, \qquad
    \bm{E}_{\beta} = \begin{pmatrix}
    K_{\rm sh}(\grad w - \bm{\theta})\\
    \bm{K}_b \Grad \bm{\theta} \\
    \end{pmatrix} = 
    \begin{pmatrix}
    \bm{q}\\
    \bm{M} \\
    \end{pmatrix},
\end{equation}
where uppercase is used to denote tensor variables. The system equivalent to \eqref{eq:energy_co-energy-system} is then written as
\begin{equation}
    \begin{bmatrix}
        \rho h & 0 & 0 & 0 \\
        0 & \rho J & 0 & 0 \\
        0 & 0 &  C_{\rm sh} & 0 \\
        0 & 0 & 0 & \bm{C}_b
    \end{bmatrix}
    \begin{pmatrix}
        \partial_t v \\
        \partial_t\bm{\omega} \\
        \partial_t \bm{q} \\
        \partial_t \bm{M} \\
    \end{pmatrix} = 
    \begin{bmatrix}
        0 & 0 & \div & 0 \\
        0 & 0 & \bm{I}_2 & \Div \\
        \grad & -\bm{I}_2 & 0 & 0 \\
        0 & \Grad & 0 & 0 \\
    \end{bmatrix}
    \begin{pmatrix}
        v \\
        \bm{\omega} \\
        \bm{q} \\
        \bm{M} \\
    \end{pmatrix},
\end{equation}
where $C_{\rm sh}:=K^{-1}_{\rm sh}$ is the shear compliance and $\bm{C}_b := \bm{K}_b^{-1}$ is the bending compliance tensor. The differential operator $\mathcal{L}$  and its formal adjoint $\mathcal{L}^*$ for the Mindlin plate problem are given by
$$
\mathcal{L}= \begin{bmatrix}
     \grad & -\bm{I}_2  \\
     0 & \Grad
\end{bmatrix}, \qquad 
\mathcal{L}^*= - \begin{bmatrix}
     \div & 0 \\
    \bm{I}_2 & \Div 
\end{bmatrix}
$$ 
The discretization is obtained by multiplying by the test functions and applying integration by parts as in \eqref{eq:weak_conforming_L} or \eqref{eq:weak_conforming_L*}. The resulting weak formulation for $\Omega_1$ is to find $v \in L^2(\Omega_1), \; \bm{\omega}\in L^2(\Omega_1; \bbR^2)$  and $\bm{q}\in H^{\div}(\Omega_1), \; \bm{M}\in H^{\Div}(\Omega_1; \bbS)$, such that $\forall{\psi}_v\in L^2(\Omega_1), \; \forall\bm{\psi}_\omega\in L^2(\Omega_1; \bbR^2), \; \forall\bm{\psi}_q\in H^{\Div}(\Omega_1)$ and $\forall\bm{\Psi}_M\in H^{\div}(\Omega_1; \bbS)$ it holds
\begin{equation}
    \begin{split}
    \begin{alignedat}{2}
        \innerproduct[\Omega_1]{\psi_v}{\rho h \partial_t v} &= + \innerproduct[\Omega_1]{\psi_v}{\div\bm{q}}, \\
        \innerproduct[\Omega_1]{\bm\psi_\omega}{\rho J \partial_t \bm{\omega}} &= + \innerproduct[\Omega_1]{\bm\psi_\omega}{\Div\bm{M} + \bm{q}}, \\
        \innerproduct[\Omega_1]{\bm{\psi}_q}{C_{\rm sh}\partial_t\bm{q}} &=- \innerproduct[\Omega_1]{\div\bm{\psi}_q}{v} - \innerproduct[\Omega_1]{\bm{\psi}_q}{\bm{\omega}} + \boundary[\Gamma_1]{\bm{\psi}_q \cdot \bm{n}}{v} + \boundary[\Gamma_{\rm int}]{\bm{\psi}_q \cdot \bm{n}}{v}, \\
        \innerproduct[\Omega_1]{\bm{\Psi}_M}{\bm{C}_b\partial_t\bm{M}} &=- \innerproduct[\Omega_1]{\Div\bm{\Psi}_M}{\bm{\omega}} + \boundary[\Gamma_1]{\bm{\Psi}_M \cdot \bm{n}}{\bm{\omega}} + \boundary[\Gamma_{\rm int}]{\bm{\Psi}_M \cdot \bm{n}}{\bm{\omega}}, \\
    \end{alignedat}
    \end{split}
\label{eq:Mindlin-weak-D}
\end{equation}
where the trace operator $\mathcal{T}_\beta$ is the normal trace applied to a consisting of a tensor and a vector
$$
\mathcal{T}_\beta \begin{pmatrix}
    \bm{q} \\
    \bm{M} \\
\end{pmatrix} = 
\begin{pmatrix}
    \bm{q}\cdot \bm{n}\vert_{\partial\Omega_1} \\
    \bm{M}\cdot \bm{n}\vert_{\partial\Omega_1} \\
\end{pmatrix}
$$
For the subdomain $\Omega_2$ where the Neumann boundary condition is natural, the functional setting is the following: seek $v \in H^1(\Omega_2; \mathbb{R}^2), \;  \bm{\omega}\in H^1(\Omega_2; \mathbb{R}^2)$ and $\bm{M}\in H^{\rot\Rot}(\Omega_2; \mathbb{S}), \; \bm{q}\in H^{\rot}(\Omega_2)$ to satisfy $\forall\psi_v \in H^1(\Omega_2), \; \forall\bm\psi_\omega \in H^1(\Omega_2; \bbR^2), \; \forall\bm{\Psi}_M\in H^{\rot\Rot}(\Omega_2; \mathbb{S}), \; \forall\bm{q}\in H^{\rot}(\Omega_2)$
\begin{equation}
    \begin{split}
    \begin{alignedat}{2}
        \innerproduct[\Omega_1]{\psi_v}{\rho h \partial_t v} &= - \innerproduct[\Omega_1]{\grad\psi_v}{\bm{q}} + \boundary[\Gamma_1]{\psi_v}{\bm{q} \cdot \bm{n}} + \boundary[\Gamma_{\rm int}]{\psi_v}{\bm{q} \cdot \bm{n}}, \\
        \innerproduct[\Omega_1]{\bm\psi_\omega}{\rho J \partial_t \bm{\omega}} &= - \innerproduct[\Omega_1]{\Grad\bm\psi_\omega}{\bm{M}} + \innerproduct[\Omega_1]{\bm\psi_\omega}{\bm{q}} + \boundary[\Gamma_1]{\bm{\psi}_\omega}{\bm{M}\cdot \bm{n}} + \boundary[\Gamma_{\rm int}]{\bm{\psi}_\omega}{\bm{M}\cdot \bm{n}}, \\
        \innerproduct[\Omega_1]{\bm{\psi}_q}{C_{\rm sh }\partial_t\bm{q}} &= +\innerproduct[\Omega_1]{\bm{\psi}_q}{\grad v - \bm{\omega}},\\
        \innerproduct[\Omega_1]{\bm{\Psi}_M}{\bm{C}_b\partial_t\bm{M}} &=+\innerproduct[\Omega_1]{\bm{\Psi}_M}{\Grad \bm{\omega}}, \\
    \end{alignedat}
    \end{split}
\label{eq:Mindlin-weak-N}
\end{equation}
where $\mathcal{T}_\alpha$ is the Dirichlet trace applied to a tuple consisting of a scalar and a vector
$$
\mathcal{T}_\alpha \begin{pmatrix}
    v \\
    \bm{\omega}
\end{pmatrix} = 
\begin{pmatrix}
    v\vert_{\partial\Omega_2}\\
    \bm{\omega}\vert_{\partial\Omega_2}
\end{pmatrix}.
$$
\paragraph{Finite element spaces}
The mesh consists of a structured triangular mesh. On the $\Omega_1$ subdomain, Discontinuous Galerkin of order $1$ (DG$_1$)  are used for $v$ and $\bm{\omega}$, Raviart-Thomas elements of degree two RT$_2$ are used for $\bm{q}$ and conforming Arnold Winther elements \cite{arnold2002elasticity} of degree three are used for $\bm{M}$. On the $\Omega_2$ subdomain Continuous Galerkin of order (CG$_2$) element for $v, \; \bm{\omega}$, and Discontinuous Galerkin of order $1$ (DG$_1$) for $\bm{q}, \; \bm{M}$.  The justification for this choice comes from the fact that the Mindlin plate combines the wave equation with 2D elasticity. Therefore the finite element subcomplex are the same as the ones used in Sec. \ref{sec:Wave} and \ref{sec:elasticity}. The solution is found on union of meshes, that is $\mathfrak{T}_h=\mathfrak{T}_h^{\Omega_1}\cup\mathfrak{T}_h^{\Omega_2}$, with finite dimensional spaces for the $\Omega_1$ subdomain given by
\begin{equation}
    \begin{split}
         V_{\alpha, 1} &= \{\bm{u}_h\in L^2(\Omega_1) \times L^2(\Omega_1; \bbR^2)|\; \forall T\in\mathfrak{T}_h^{\Omega_1}, \; \bm{u}_h|_T\in \mathrm{DG}_1 \times \mathrm{DG}_1(\bbR^2)\}, \\
         V_{\beta, 1} &= \{\bm{S}_h\in H^{\div}(\Omega_1) \times H^{\Div}(\Omega_1; \bbS)|\; \forall T\in\mathfrak{T}_h^{\Omega_1}, \; \bm{S}_h|_T\in \mathrm{RT}_2 \times \mathrm{AW}_3\}, 
    \end{split}
\end{equation}
where $T$ now denotes a triangular mesh element of $\mathfrak{T}_h$. For the $\Omega_2$ subdomain the mixed finite element spaces are
\begin{equation}
    \begin{split}
        V_{\alpha, 2} &= \{\bm{u}_h\in H^1(\Omega_2) \times H^1(\Omega_2; \bbR^2)|\; \forall T\in\mathfrak{T}_h^{\Omega_2},\; \bm{u}_h|_T\in \mathrm{CG}_2 \times \mathrm{CG}_2(\bbR^2)\}, \\
        V_{\beta, 2} &= \{\bm{S}_h\in L^2(\Omega_2; \bbR^2) \times L^2(\Omega_2; \bbS)|\; \forall T\in\mathfrak{T}_h^{\Omega_2}, \bm{S}_h|_T\in \mathrm{DG}_1(\bbR^2) \times \mathrm{DG}_1(\bbS)\}.
    \end{split}
\end{equation}
Introducing the finite element approximation, the following ODE is obtained for subdomain $\Omega_1$
\begin{equation}\label{eq:discrete_clamped_mindlin}
\mathrm{Diag}
\begin{bmatrix}
\mathbf{M}_{v, 1} \\
\mathbf{M}_{w, 1} \\
\mathbf{M}_{m, 1} \\
\mathbf{M}_{q, 1} \\
\end{bmatrix}
    \odv{}{t}
    \begin{pmatrix}
        \mathbf{v}_1 \\
        \mathbf{w}_1 \\
        \mathbf{m}_1 \\
        \mathbf{q}_1
    \end{pmatrix} = 
    \begin{bmatrix}
        0 & 0 & \mathbf{D}_{\div} & 0 \\
        0 & 0 & \mathbf{P} & \mathbf{D}_{\Div} \\
        -\mathbf{D}_{\div}^\top & -\mathbf{P}^\top & 0 & 0 \\
        0 & -\mathbf{D}_{\Div}^\top & 0 & 0
    \end{bmatrix}
    \begin{pmatrix}
        \mathbf{v}_1 \\
        \mathbf{w}_1 \\
        \mathbf{m}_1 \\
        \mathbf{q}_1
    \end{pmatrix} + \begin{bmatrix}
        0 & 0 \\
        0 & 0 \\
         \mathbf{B}_{\bm{q} \cdot \bm{n}} & 0 \\
        0 & \mathbf{B}_{\bm{M}\cdot \bm{n}} \\
    \end{bmatrix}\begin{pmatrix}
        \mathbf{u}_v \\
        \mathbf{u}_\omega
    \end{pmatrix}.
\end{equation}
The discrete system for the domain $\Omega_2$ is given by
\begin{equation}\label{eq:discrete_free_mindlin}
\mathrm{Diag}
\begin{bmatrix}
\mathbf{M}_{v, 2} \\
\mathbf{M}_{w, 2} \\
\mathbf{M}_{m, 2} \\
\mathbf{M}_{q, 2} \\
\end{bmatrix}
    \odv{}{t}
    \begin{pmatrix}
        \mathbf{v}_2 \\
        \mathbf{w}_2 \\
        \mathbf{n}_2 \\
        \mathbf{m}_2
    \end{pmatrix} = 
    \begin{bmatrix}
        0 & 0 & -\mathbf{D}_{\grad}^\top & 0 \\
        0 & 0 & \mathbf{P} & -\mathbf{D}_{\Grad}^\top \\
        \mathbf{D}_{\grad} & -\mathbf{P}^\top & 0 & 0 \\
        0 & \mathbf{D}_{\Grad} & 0 & 0
    \end{bmatrix}
    \begin{pmatrix}
        \mathbf{v}_2 \\
        \mathbf{w}_2 \\
        \mathbf{n}_2 \\
        \mathbf{m}_2
    \end{pmatrix} + \begin{bmatrix}
        \mathbf{B}_{v} & 0 \\
        0 & \mathbf{B}_{\bm\omega} \\
        0 & 0 \\
        0 & 0 \\
    \end{bmatrix}\begin{pmatrix}
        \mathbf{u}_q \\
        \mathbf{u}_M
    \end{pmatrix}.
\end{equation}

\begin{figure}[htbp]
	\centering
	\begin{subfigure}{0.32\textwidth}
		\centering
		\includegraphics[width=\textwidth]{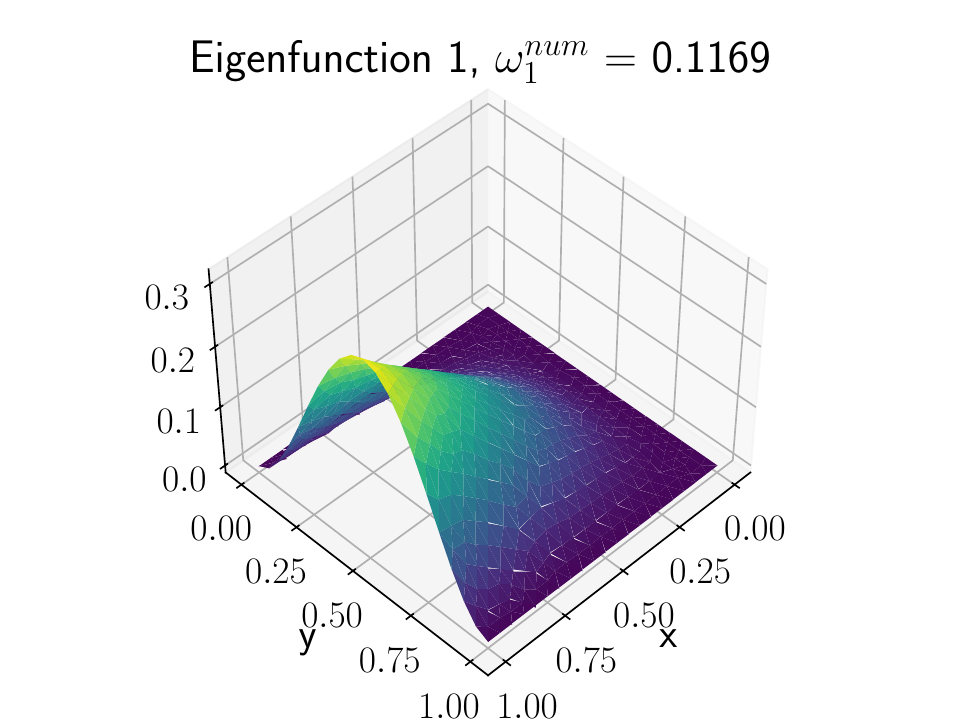}
	\end{subfigure}
	\begin{subfigure}{0.32\textwidth}
		\centering
		\includegraphics[width=\textwidth]{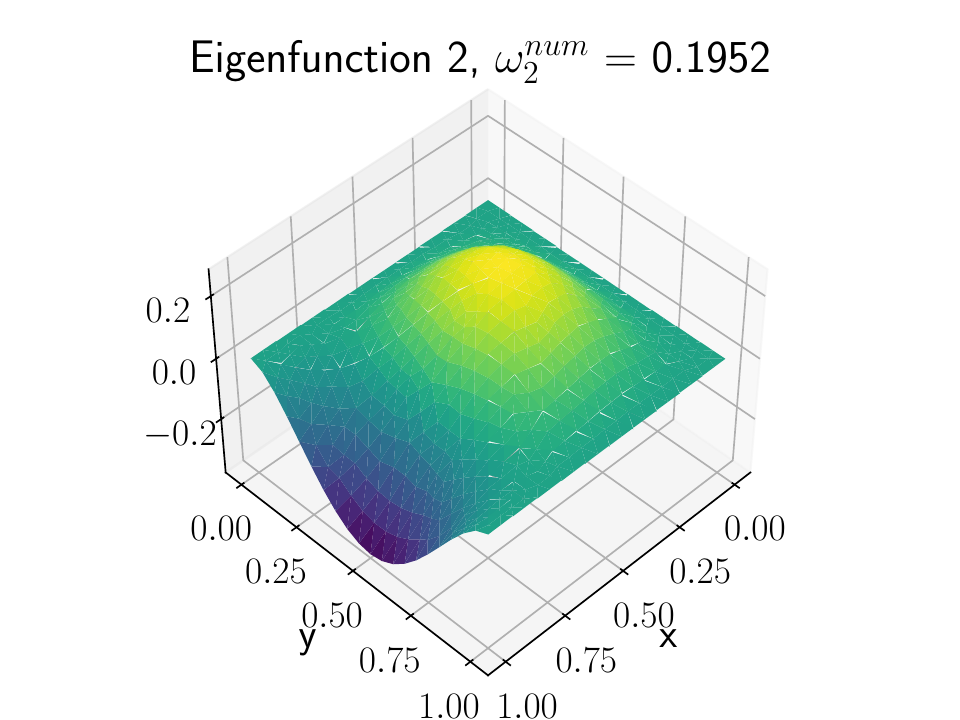}
	\end{subfigure}
        \begin{subfigure}{0.32\textwidth}
		\centering
		\includegraphics[width=\textwidth]{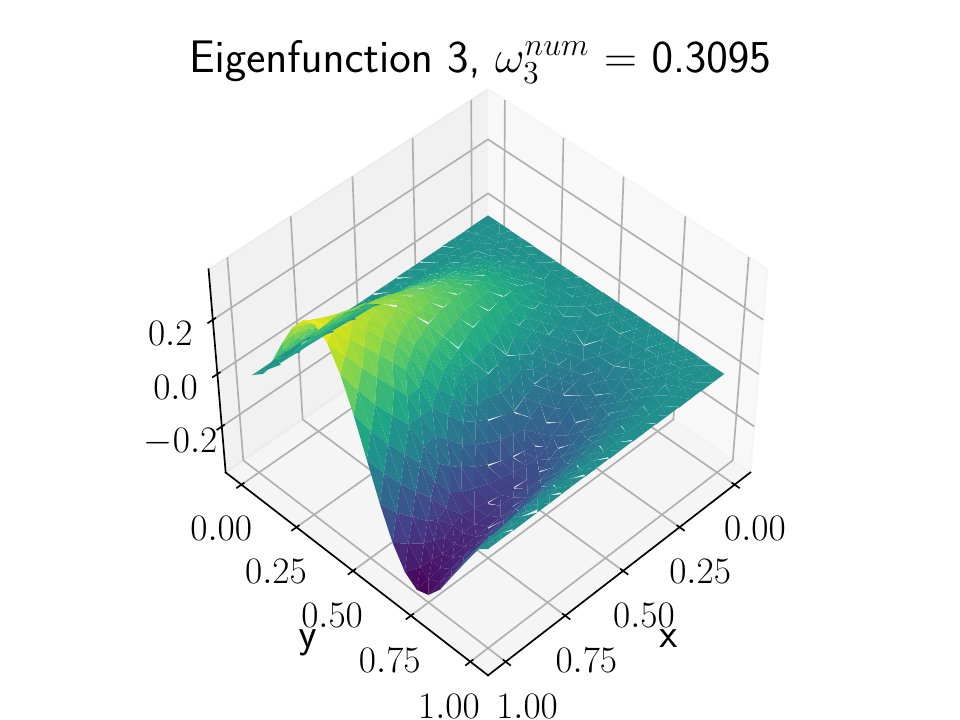}
	\end{subfigure}
	\begin{subfigure}{0.32\textwidth}
		\centering
		\includegraphics[width=\textwidth]{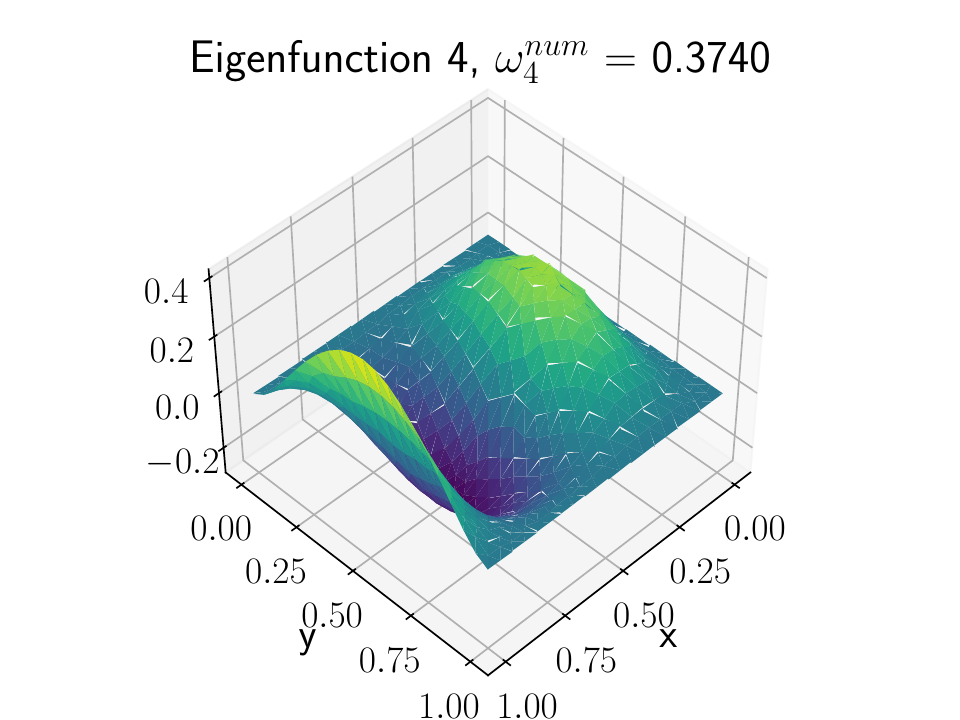}
	\end{subfigure}
        \begin{subfigure}{0.32\textwidth}
		\centering
		\includegraphics[width=\textwidth]{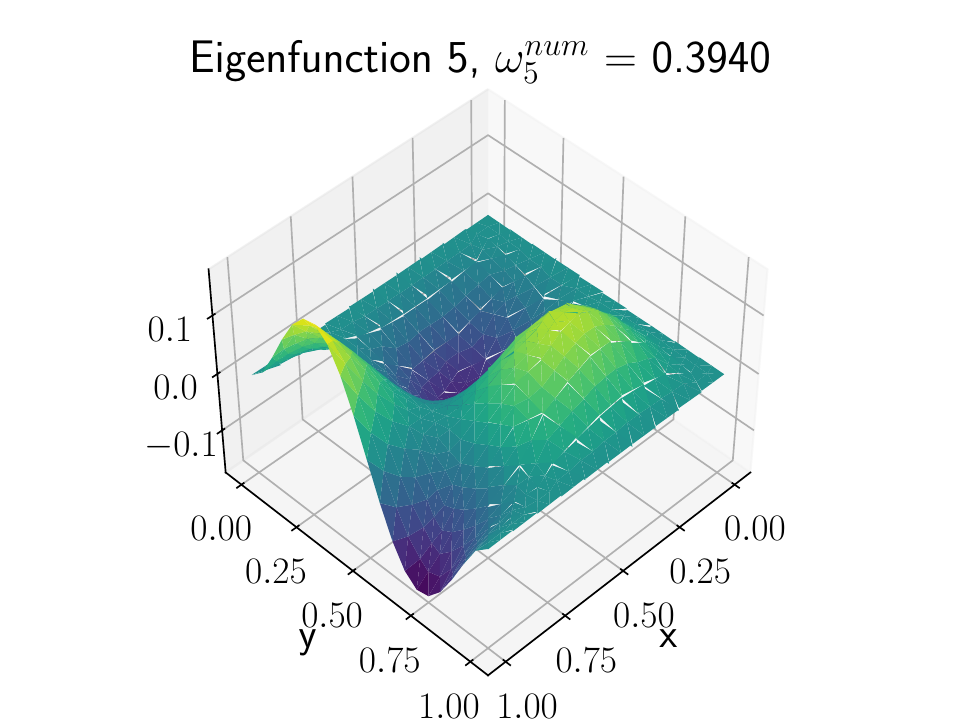}
	\end{subfigure}
	\begin{subfigure}{0.32\textwidth}
		\centering
		\includegraphics[width=\textwidth]{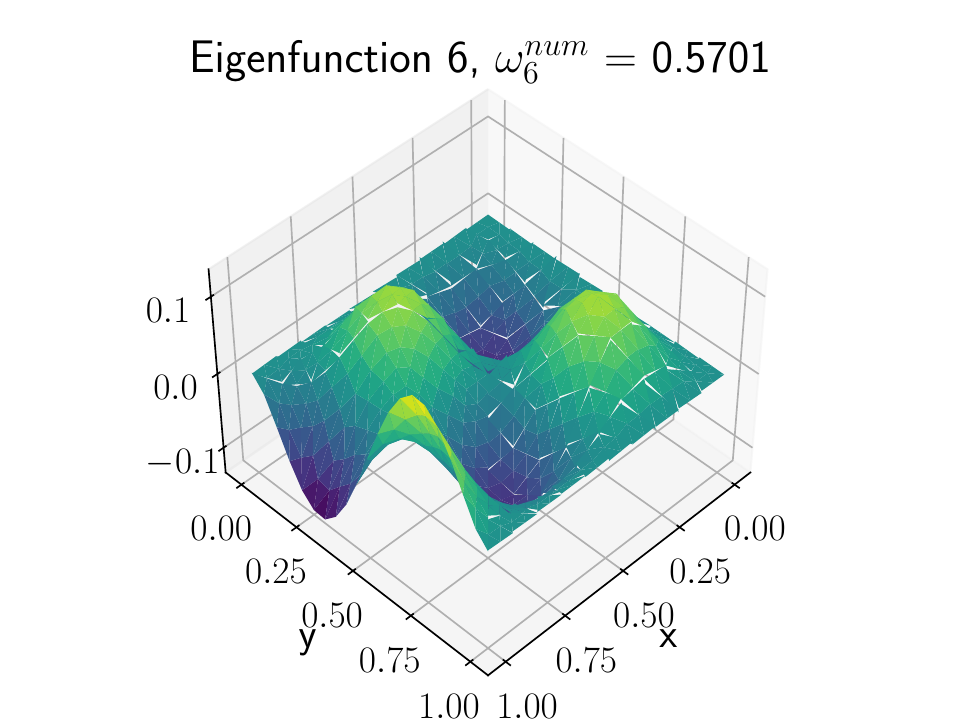}
	\end{subfigure}
	\caption{Numerical eigenvectors in terms of the vertical velocity $v$ for the Mindlin plate. The proposed approach and a standard primal discretization are compared with the results of \cite{dawe1980mindlin}.}
	\label{fig:eigenvectors_mindlin}
\end{figure}

\paragraph{Modal analysis} 
We consider an analogous setting to the linear elastodynamics example, on a unit square $\Omega = [0, 1]^2$, decomposed by an interface placed diagonally between the lower left and upper right vertex as shown in Fig. \ref{fig:elasticity_mesh}. The numerical eigenvalues are obtained via the generalized eigenproblem
$$
i \omega^{num}_n \mathbf{M}\bm{\psi}_n = \mathbf{J}\bm{\psi}_n,
$$
where $i= \sqrt{-1}$ is the imaginary unit. For the discretization $10$ finite elements per side are considered. The physical and geometrical paramters are  parameters are those of aluminum
$$
L=1 \; \mathrm{[m]}, \qquad h=0.01 \; \mathrm{[m]}, \qquad \rho = 2700 \; \mathrm{[Kg/m^3]}, \qquad E = 70 \; \mathrm{[GPa]}, \qquad \nu=0.3, \qquad k=0.8601.
$$
We consider a small thickness to show the robustness of the proposed methodology against shear locking phenomena. The resulting normalized eigenvalues, given by
$$
\widehat{\omega} = \omega \; L\sqrt{\frac{\rho}{G}}, \qquad \qquad G=\frac{E}{2(1+\nu)},
$$
are shown in Table \ref{tab:eigenvalues_mindlin} whereas the eigenvectors for variable $v$ (corresponding to the velocity) are plotted in Fig. \ref{fig:eigenvectors_mindlin}. The obtained eigenfrequencies are compared against a classical primal discretization using quadratic elements for $w$ and $\bm{\theta}$, with respect to reference \cite{dawe1980mindlin}, where an analytical approach is used. It can be noticed that the proposed discretization achieves better results and the mixed discretization avoids shear locking phenomena introducing different variables for $\bm{\omega}$ and $\bm{q}$.
\begin{table}[htbp]
\centering
\begin{tabular}{c c c c c}
\hline
Mode 
& Proposed 
& Classical 
& Reference \cite{dawe1980mindlin}
& Rel. Error (\%) \\
\cline{5-5}
& & & & Proposed \quad Classical \\
\hline
1 & 0.1168 & 0.1192 & 0.1171 & 0.256 \quad 1.793 \\
2 & 0.1951 & 0.2008 & 0.1951 & 0.000 \quad 2.922 \\
3 & 0.3094 & 0.3208 & 0.3093 & 0.032 \quad 3.718 \\
4 & 0.3739 & 0.3876 & 0.3740 & 0.027 \quad 3.636 \\
5 & 0.3940 & 0.4137 & 0.3931 & 0.229 \quad 5.240 \\
6 & 0.5700 & 0.6050 & 0.5695 & 0.088 \quad 6.234 \\
\hline
\end{tabular}
\caption{Numerical eigenvalues of the Mindlin plate problem using the proposed method and a standard finite element discretization.}
\label{tab:eigenvalues_mindlin}
\end{table}

%% file: Text/Conclusion.tex
\section{Conclusion}\label{sec:concl}

\revthree{
In this contribution a numerical strategy to impose mixed boundary conditions in port-Hamiltonian systems has been detailed. To this aim, a primal-dual formulation leveraging the machinery of Hilbert complexes has been used. The domain decomposition can be reinterpreted under the lens of Discontinuous Galerkin method: the interconnection of the two subdomains corresponds to a physically motivated choice of numerical fluxes. Indeed the natural boundary condition for each subdomain corresponds to the output of the dual formulation. The time integration can be performed using methods capable of preserving the Poisson structure of the system. A natural choice in the linear case is given by the implicit midpoint method. Integrators capable of decoupling the two subdomains, like the St\"ormer-Verlet scheme, are also of interest as they reduce the computational burden. The methodology can be extended to nonlinear problems but is limited to the case where the nonlinearity does not interfere with the differential operators. \\
The definition of an interface between two subdomains may represent a bottleneck in applications where the boundary subpartitions present an intricate topology. However, this contribution represents a proof of concept that shows that primal-dual mixed formulations can be used simultaneously. Problems with mixed boundary conditions can then be represented as ordinary differential equations rather than differential algebraic ones. This presents advantages that go beyond simulation purposes, as the removal of algebraic constraints is beneficial also in the context of numerical optimization.\\
An interesting perspective would be to push forward the method and perform the interconnection on each finite element. This would lead to a completely discontinuous Galerkin method where each finite element exchanges information with the adjacent elements via a feedback interconnection. Furthermore, it can be integrated as an actual domain decomposition approach to avoid the computational cost of solving monotonically large systems arising from partial differential equations. This idea may be combined with model reduction approaches to reduce each subdomain before performing the interconnection. The presented idea may also find application in static problems, very much in the same spirit of hybrid and discontinuous methods.
}